\documentclass[12pt]{article}

\usepackage{amssymb}
\usepackage{amsmath}
\usepackage{amscd}
\usepackage{ifpdf}

\usepackage[dvips]{graphicx}
\usepackage[ruled]{algorithm}
\usepackage{algorithmic}
\usepackage{palatino, url, multicol}

\numberwithin{equation}{section}

\begin{document}

\title{STORAGE ALLOCATION UNDER PROCESSOR SHARING II: FURTHER ASYMPTOTIC RESULTS}

\author{
Eunju Sohn\thanks{ Department of Mathematics, Statistics, and
Computer Science, University of Illinois at Chicago, 851 South
Morgan (M/C 249), Chicago, IL 60607-7045, USA. {\em Email:}
esohn3@math.uic.edu.} \and \and Charles Knessl\thanks{ Department of
Mathematics, Statistics, and Computer Science, University of
Illinois at Chicago, 851 South Morgan (M/C 249), Chicago, IL
60607-7045, USA. {\em Email:} knessl@uic.edu.\
\newline\indent\indent{\bf Acknowledgement:} This work was partly supported by NSF grant DMS 05-03745 and NSA grant H 98230-08-1-0102.
}}
\date{ }
\maketitle

\begin{abstract}
\noindent  We consider a processor sharing storage allocation model,
which has $m$ primary holding spaces and infinitely many secondary
ones, and a single processor servicing the stored items (customers). All of the
spaces are numbered and ordered. An arriving customer takes the
lowest available space. We define the traffic intensity $\rho$ to be
$\lambda / \mu$ where $\lambda$ is the customers' arrival rate and $\mu$ is
the service rate of the processor. We study the joint probability
distribution of the numbers of occupied primary and secondary spaces. We study the problem in two asymptotic limits: (1) $m\rightarrow \infty$ with a fixed $\rho <1$, and (2) $\rho \uparrow 1, \; m \rightarrow \infty$ with $m(1-\rho)= O(1)$.
\end{abstract}

\setlength{\baselineskip}{24pt}

\section{Introduction}

\indent We consider the following storage allocation model. Suppose that near a restaurant there are $m$ primary parking spaces and across
the street there are infinitely many additional ones. However, the
restaurant has only one waiter who serves all of the customers. All of
the parking spaces are numbered and ordered; the one with rank $= 1$
is closest to the restaurant and the primary spaces are numbered $\{1, 2, 3, ....., m\}$.
  We assume the following: (1) customers arrive
according to a Poisson process with rate $\lambda$, (2) the waiter works at rate $\mu$, (3) an arriving car
parks in the lowest-numbered available space, and (4) if there are $N$
customers in the restaurant, the waiter serves each customer at the
rate $\mu / N$. This corresponds to a processor sharing (PS) service discipline.

\indent Dynamic storage allocation and the fragmentation of computer
memory are among the many applications of this model.  We define
$N_1$ to be the number of occupied primary spaces and $N_2$ to be the number of occupied
secondary spaces.  Then we define $\textbf{S}$ to be the set of the indices of the occupied spaces, and the "wasted spaces" $W$ are defined as the difference between the largest
index of the occupied spaces (Max $\textbf{S}$) and the total
number of occupied spaces ($|\textbf{S}|=N_1+N_2$).  Coffman, Flatto, and Leighton \cite{CFL}
showed that for the processor-sharing model
$$E[W] = \Theta \left( \sqrt{\frac{1}{1-\rho}\log \left(\frac{1}{1-\rho}\right)} \right), \; \; \rho \uparrow 1$$
where $E[W]$ is the expected value of the wasted spaces.  Here $E[W] = \Theta(f(\rho))$ means that there exist positive constants $c, c'$ such that $c'f(\rho)\leq E[W] \leq c f(\rho)$. Also when $\rho \rightarrow 1$ (the heavy traffic case) Coffman and Mitrani \cite{CM}
obtained upper and lower bounds on $E[W]$ in the form
$$\frac{1}{2}\sqrt{\frac{\pi}{1-\rho}} \leq E[W]\leq\frac{1}{1-\rho}\left(\frac{\pi^2}{6}-1\right).$$

\indent A related model, the $M/M/\infty$ queue with ranked servers, has been studied by many authors \cite{A}, \cite{CK}, \cite{CL}, \cite{K}, \cite{P}, \cite{N}. This differs from the current model in that if there are a total of $N=N_1+N_2$ spaces occupied, the total service rate is $\mu N$, as each customer in the restaurant is served at rate $\mu$. For this model various asymptotic studies appear in \cite{A}, \cite{CL}, \cite{N}. In particular, Aldous \cite{A} showed that the mean number of the wasted spaces is $E[W] \sim \sqrt{2 \rho \log\log \rho}\; $ as $\; \rho = \lambda/\mu \rightarrow \infty$.

\indent A simple derivation of the exact joint distribution of finding $N_1$ (resp., $N_2$) occupied primary (resp., secondary) spaces appears in \cite{S} and detailed asymptotic results for this joint distribution appear in \cite{KB}, \cite{KC}, while the distribution of Max $\textbf{S}$ is analyzed in \cite{KA}. In \cite{KB} Knessl showed how to obtain asymptotic results for the infinite server model directly from the basic difference equation. Since the present processor sharing model does not seem amenable to exact solution, we shall employ such a direct asymptotic approach here.

\indent In this paper, we study the joint probability distribution of the numbers of occupied spaces in the PS model, letting $\pi(k, r) =  Prob[N_1=k, N_2=r]$ in the steady state. In part I \cite{KSS} we obtained exact solutions for $m=1$ and $m=2$, and developed a semi-numerical semi-analytic method for general $m$.  We also derived asymptotic results in the heavy traffic case $\rho\uparrow 1$, but with $m=O(1)$.  Here we shall obtain asymptotic results for $m\rightarrow \infty$, for the cases $0<\rho<1$ and $1-\rho=O(m^{-1})$.

\indent The paper is organized as follows. In section \textbf{2} we
state the problem and obtain the basic difference equations for $\pi(k,r)$. In section \textbf{3} we summarize our main results. In section \textbf{4} we consider $m \rightarrow \infty$ with a fixed $\rho$ ($0<\rho<1$) and study $\pi(k,r)$ for various ranges of $(k,r)$.  In section \textbf{5} we consider the double limit $m \rightarrow \infty$ and $\rho \uparrow 1$, with $m(1-\rho) = O(1)$. Some numerical studies and comparisons appear in section \textbf{6}.

\section{Statement of the problem}
We let $N_1 (t)$ (resp., $N_2 (t)$) denote the number of primary
(resp., secondary) spaces occupied at time $t$. The joint steady
state distribution function is
\begin{equation*}
\pi (k, r) = \pi (k, r; m)=\lim_{t\rightarrow\infty} Prob[ N_1 (t) = k,  N_2(t) = r ], \; 0\le k \le m,\; r \ge 0.
\end{equation*}

Let $\rho = \lambda/\mu$ be the traffic intensity and we assume the stability condition $\rho < 1$.
 The pair $(N_1, N_2)$ forms a Markov chain whose transition rates are sketched in \textbf{Fig. 1}. The state space is the lattice strip $\{(k,r): 0\leq k \leq m, \;  r\geq 0\}$ and the balance equations are
\begin{multline}\label{bl} (1I_{[k+r>0]}+\rho)\; \pi(k,r) = \rho \; \pi (k-1,
r)I_{[k\geq1]}+\frac{k+1}{k+r+1}\; \pi (k+1, r)I_{[k < m]}\\
 +\frac{r+1}{k+r+1}\; \pi (k, r+1)+\rho\; \pi (m, r-1)I_{[k=m,\; r\geq1]}.
\end{multline}
Here $I$ is an indicator function. The normalization condition is
\begin{equation} \label{n}
\sum_{r=0}^\infty \sum_{k=0}^m \pi(k,r) = 1.
\end{equation}

\indent From our viewpoint we will need to consider explicitly the boundary conditions inherent in (\ref{bl}), so we rewrite the main equation as
\begin{multline} \label{bl2}
(1+\rho)\; \pi(k,r) = \rho \; \pi (k-1,
r)+\frac{k+1}{k+r+1}\; \pi (k+1, r)
 +\frac{r+1}{k+r+1}\; \pi (k, r+1), \\ \; 0\leq k < m, \; r\geq 0, \; k+r > 0.
\end{multline}
and the boundary condition at $k=m$ is
\begin{equation} \label{bcm}
(1+\rho)\; \pi(m,r)\\ = \rho \; \pi (m-1,
r)+\frac{r+1}{m+r+1}\; \pi (m, r+1)+\rho\; \pi (m, r-1), \; r\geq1.
\end{equation}
There are also the two corner conditions
\begin{equation}
\rho \pi(0,0) = \pi(1,0)+\pi(0,1) \label{26}
\end{equation}
and
\begin{equation} \label{27}
(1+\rho)\pi(m,0) = \rho\pi(m-1,0)+\frac{1}{m+1} \pi(m,1).
\end{equation}
In (\ref{bl2}) when $k=0$ we interpret $\pi(-1,r)$ as 0. The boundary condition at $k=m$ in (\ref{bcm}) can be replaced by the
artificial boundary condition
\begin{equation}
\frac{m+1}{m+r+1}\pi(m+1, r) =
\rho\;\pi(m, r-1). \label{ab}
\end{equation}
This is obtained by extending (\ref{bl2}) to hold also at $k=m$ and comparing this to (\ref{bcm}).

\indent We note that the total number $N_1+N_2$ behaves as the number of customers in the $M/M/1-PS$ queue, which is well known to follow a geometric distribution. Thus we have
\begin{equation} \label{29}
\sum_{k+r=N} \pi(k,r) = (1-\rho)\rho^N, \; N \geq 0
\end{equation}
and we can rewrite this as
\begin{eqnarray}
\sum_{r=0}^N \pi(N-r, r) &=& (1-\rho)\rho^N, \; \; 0\leq N \leq m, \\
\sum_{r=N-m}^N \pi(N-r,r) &=& \sum_{k=0}^m \pi(k, N-k) = (1-\rho)\rho^N, \; N \geq m. \label{mm1}
\end{eqnarray}
These identities will provide a useful check on the calculations that follow.
\begin{center}
\textbf{Fig. 1}   A sketch of the transition rates.
\includegraphics[angle=0, width=1.0\textwidth]{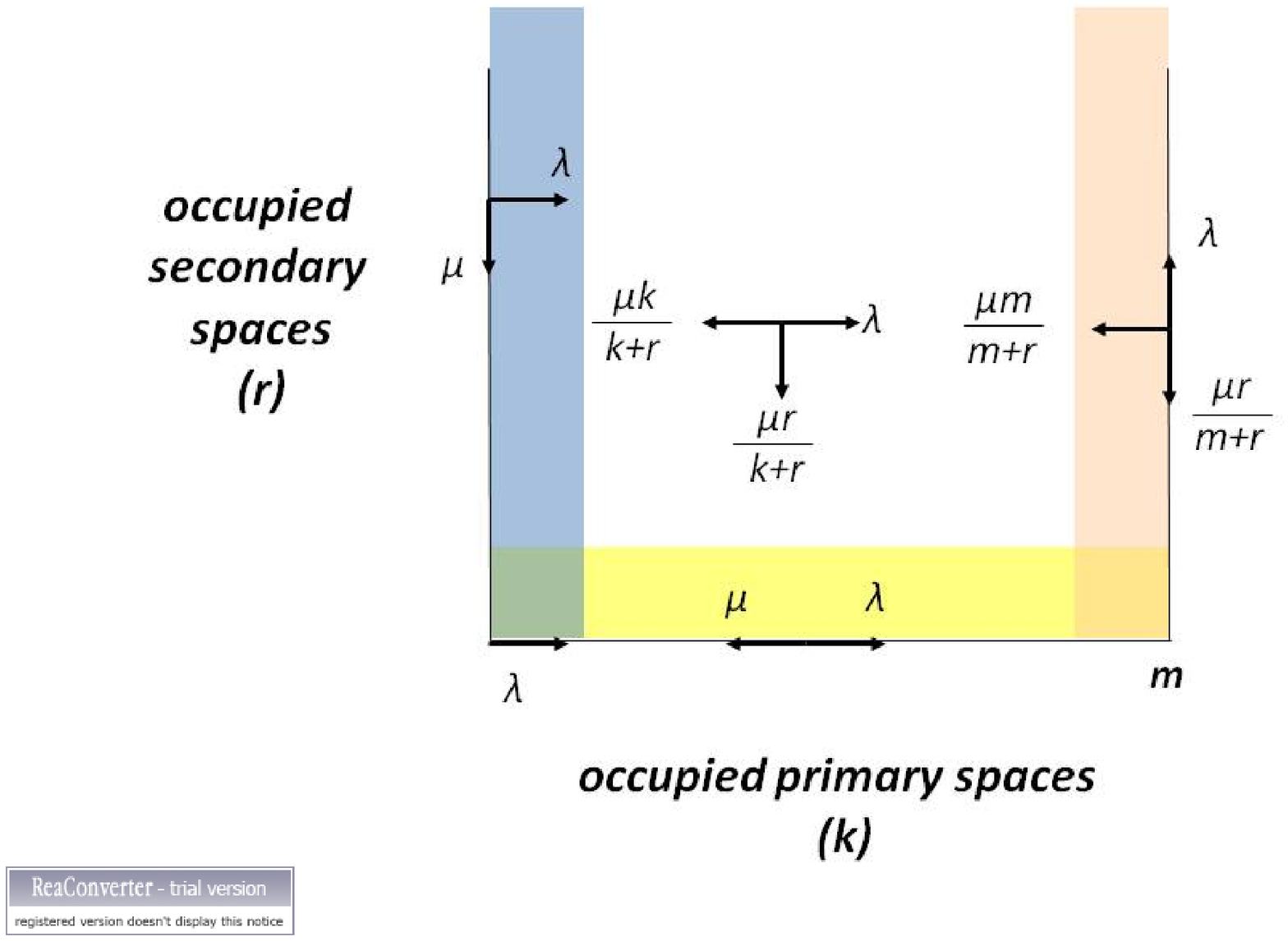}
\end{center}

\section{Summary}
Since the analysis shall become quite involved and technical, we collect here some of the main results.

We take $m\rightarrow \infty$, first assuming that $\rho<1$.  We use the scaled variables $x=k/m$ and $y=r/m$, so that $x$ is the fraction of primary spaces that are utilized.  The main approximation we obtain is then given by
\begin{equation}\label{30}
\pi(k,r) \sim K(x, y) e^{m\phi(x, y)}; \; \; 0<x\leq 1, \; \; y>0
\end{equation}
where
\begin{eqnarray}
\phi(x, y) &=& x\phi_x+y\phi_y+\log \rho-\log(s+1),  \label{phi}\\
\phi_x &=& \log\left[\frac{s+1-\rho-\rho s\;
e^{(1-\rho)t}}{s+1-\rho-s \; e^{(1-\rho)t}}\right],  \label{phix}\\
\phi_y &=& \log\left[\frac{\rho(s+1-\rho-\rho s\;
e^{(1-\rho)t})}{\rho(s+1-\rho)-\rho s\;
e^{(1-\rho)t}+(1-\rho)(s+1-\rho)e^t}\right], \label{phiy}
\end{eqnarray}
and $(x, y)$ are related to $(s,t)$ via the mapping
\begin{eqnarray}
x &=& x(s,t) = \frac{1}{(s+1)(1-\rho)^2}\left[s+1-\rho-s \; e^{(1-\rho)t}\right]\nonumber\\
&\times& \left[(s+1-\rho)\; e^{-(1-\rho)t}-s\rho^2-\rho(1-\rho)s e^{-t}\right],  \label{x}\\
y &=& y(s,t) = \frac{s}{(s+1)(1-\rho)}\left[(1-\rho)(s+1-\rho)-\rho s \;
e^{-\rho t}+\rho(s+1-\rho) \; e^{-t}\right]. \label{y}\nonumber\\
\end{eqnarray}
Then $K(x, y)$ is given by
\begin{eqnarray}\label{769}
K(x,y)&=& \sqrt{\frac{(1-\rho)(s+1-\rho)^3 [s+1-\rho-s\rho \;
e^{(1-\rho)t}]^3 }{(s+1-\rho-s \;
e^{(1-\rho)t})[\rho(s+1-\rho)-\rho s\;
e^{(1-\rho)t}+(1-\rho)(s+1-\rho)e^t]}}\nonumber\\
&\times& \frac{1}{(s+1)^2} \sqrt{\frac{e^{\rho t}}{\mid \jmath\mid}},
\end{eqnarray}
where $\jmath = x_ty_s-x_sy_t$ is the Jacobian of the transformation above.  This approximation to $\pi(k,r)$ is explicit in terms of $(s,t)$ but implicit in terms of $(x,y)$.  However, for $x\approx 1$ it becomes much more explicit, with
\begin{equation}\label{727}
\pi(k,r) \sim (1-\rho)\left(\frac{y+1-\rho}{y+1}\right)(y+1)^{k-m}\rho^{m(y+1)}.
\end{equation}
The above holds for $y>0$ and $k=m-O(1)$, which corresponds to all but a few occupied primary spaces.

Different expansions must be constructed in various boundary and corner regions of the strip $\{(x,y): 0\leq x\leq 1, \; \; y\geq 0\}$.  For $k=O(1)$ ($x=O(m^{-1})$),  which means that only a few primary spaces are occupied,
\begin{equation}\label{3}
\pi(k,r) \sim m^{k+1/2}e^{m\phi(0,y)}S_0(y) \frac{y^k}{k!}\left[1+\rho-e^{\phi_y(0,y)}\right]^k
\end{equation}
where
\begin{eqnarray}
\phi(0,y) &=& y\log\rho
+\frac{y}{1-\rho}\log\left[\frac{y+\sqrt{(y-1+\rho)^2+4y}-(1-\rho)}{y+\sqrt{(y-1+\rho)^2+4y}+(1-\rho)}\right] \nonumber\\
&-& \log\left[ y+\sqrt{(y-1+\rho)^2+4y}+1+\rho\right]+\log2+\log\rho, \label{760}\\
\phi_y(0,y) &=&
\log\rho+\frac{1}{1-\rho}\log\left[\frac{y+\sqrt{(y-1+\rho)^2+4y}-(1-\rho)}{y+\sqrt{(y-1+\rho)^2+4y}+(1-\rho)}\right], \label{761}\nonumber\\
 \\
S_0(y) &=& \frac{\sqrt{2\pi} (1-\rho)(s_*+1-\rho)^{3/2}\sqrt{s_*}}{(s_*+1)\sqrt{(s_*+1)^2-\rho}}, \\
s_*(y) &=& \frac{1}{2}\left[ y-1+\rho+\sqrt{(y-1+\rho)^2+4y} \right].  \label{734}
\end{eqnarray}
For $r=O(1)$  and $0<x<1$, which corresponds to having a few secondary spaces occupied and a fraction of the primary ones, we obtain
\begin{eqnarray}
\pi(k,r) &\sim& \frac{(1-\rho)^2 \rho^{m+r} x^{r/(1-\rho)}}{[1-\rho+\rho x^{1/(1-\rho)}]^{r+1}}, \; \; r\geq 1,  \label{738} \\
\pi(k,0) -(1-\rho)\rho^k &\sim&  -\frac{(1-\rho)\rho^{m+1} x^{1/(1-\rho)}}{1-\rho+\rho x^{1/(1-\rho)}}. \label{744}
\end{eqnarray}
Here we wrote the result for $r=0$ so as to estimate the deviation of $\pi(k,0)$ from the geometric distribution $(1-\rho)\rho^k = (1-\rho)\rho^{mx}.$

Near the corner $(x,y) = (0,0)$ we use the original discrete variables $(k,r)$ to find that
\begin{eqnarray}\label{763}
\pi(k,r) &\sim& \rho^{m+r} m^{-r/(1-\rho)}(1-\rho)^{2-r/(1-\rho)} \;
\Gamma\left(1+\frac{r}{1-\rho}\right) \nonumber\\
&\times& \frac{1}{2\pi i}\oint \frac{z^{-k-1}(1-\rho z)^{\rho r/(1-\rho)-1} }{(1-z)^{1+r/(1-\rho)}}dz, \; \; r\geq 1
\end{eqnarray}
and
\begin{eqnarray}
\pi(k,0)-(1-\rho)\rho^k &\sim& -\rho^m m^{-1/(1-\rho)}\sum_{j=0}^{k-1}\left(\frac{1-\rho^{k-j}}{1-\rho}\right)\frac{P(j,1)}{j+1}, \\
P(k,1) &=& \rho(1-\rho)^{2-1/(1-\rho)}\Gamma\left(\frac{2-\rho}{1-\rho}\right) \nonumber\\
&\times& \frac{1}{2\pi i}\oint \frac{z^{-k-1}(1-\rho z)^{\rho/(1-\rho)-1}}{(1-z)^{1+1/(1-\rho)}} dz.
\end{eqnarray}
Here the integrals are over a small loop about $z=0$, and these contour integrals may be expressed in terms of hypergeometric functions.  Near the other corner $(x,y)=(1,0)$ we use the variables $n=m-k$ and $r$, and obtain
\begin{eqnarray}
\pi(k,r) &\sim& (1-\rho)^2 \rho^{m+r}, \; \; r\geq1  \label{742}\\
\pi(k,0) &\sim&  (1-\rho)(\rho^{-n}-\rho)\rho^m, \; \; r=0. \label{749}
\end{eqnarray}
We note that for $\rho<1$ and $m\rightarrow \infty$ most of the probability mass occurs in the range $k=O(1)$ and $r=0$, and $\pi(k,r)$ is exponentially small in all of the other ranges.

Defining the marginal distribution by
\begin{equation}\label{m}
{\cal M} (k) = \sum_{r=0}^\infty \pi(k,r), \; \; 0\leq k\leq m
\end{equation}
and
\begin{equation}\label{nr}
{\cal N} (r) = \sum_{k=0}^m \pi(k,r), \; \; r\geq 0,
\end{equation}
we can easily obtain their expansions from the results for $\pi(k,r)$.  For the distribution of the number of occupied primary spaces we have
\begin{eqnarray}
{\cal M}(k) - (1-\rho)\rho^k  &\sim& \frac{\rho^{m+1}x^{\rho/(1-\rho)}}{(1-\rho)m}, \; \; 0<x\leq 1, \label{799}\\
{\cal M}(k) - (1-\rho)\rho^k  &\sim& \rho^m m^{-1/(1-\rho)}[P(k,0)+P(k,1)], \; \; k=O(1), \nonumber\\
\end{eqnarray}
where $P(k,0)$ is given by (\ref{765}).  For the distribution of the number of occupied secondary spaces we obtain
\begin{eqnarray}
{\cal N}(r) &\sim & (1-\rho)\rho^{m(y+1)}\left(\frac{y+1-\rho}{y}\right), \; \; y>0, \label{7b}\\
{\cal N}(r) &\sim& m (1-\rho)^3\rho^{m+r} \int_0^1 \frac{u^{r-\rho}}{(1-\rho+\rho u)^{r+1}}\; du, \; \; r\geq 1, \label{7a}\\
1-{\cal N}(0) &\sim& m(1-\rho)^2\rho^{m+1}\int_0^1 \frac{u^{1-\rho}}{1-\rho+\rho u} \; du.
\end{eqnarray}
In particular the mean number of occupied secondary spaces is
\begin{equation}
\sum_{r=1}^\infty r{\cal N}(r) \sim \left(\frac{1-\rho}{2-\rho}\right)m\rho^{m+1}.
\end{equation}

Finally we consider the double limit where $m\rightarrow \infty$ and $\rho \uparrow 1$.  We introduce the parameter $a=m(1-\rho)=O(1)$.  On the ($x,y$) scale we find that
\begin{equation}
\pi(k,r) \sim m^{-1} {\cal K}(x,y)e^{m\Psi(x,y)} \; ; \; 0<x\leq1, \; \; y>0
\end{equation}
where
\begin{equation}\label{Psi}
\Psi = x\Psi_x+y\Psi_y-\log(s+1)
\end{equation}
with
\begin{eqnarray}
\Psi_x &=& \log\left(\frac{1+s-st}{1-st}\right),\label{Psix}\\
\Psi_y &=& \log\left(\frac{1+s-st}{1-st+s e^t}\right),  \label{Psiy} \\
x(s,t)&=& \frac{1}{s+1}(1-st)(1+2s-st - s e^{-t}), \label{817}\\
y(s,t) &=& \frac{s}{s+1}\left[s+(1-st)e^{-t}\right]. \label{818}
\end{eqnarray}
and
\begin{eqnarray}\label{833}
{\cal K} &=&  \frac{as^{3/2}(s-st+1)e^{t/2}}{\sqrt{s(s+1)(1-st)(se^t-st+1)}\sqrt{2s(s+2)-[s(s+2)t^2-2t+s(s+2)-1]e^{-t}}} \nonumber\\
&\times& \exp \left[ast -\frac{as^2t^2}{2(s+1)}-a(s+1)\right].
\end{eqnarray}
For $k=m-O(1)$ ($x=1-O(m^{-1})$) the expression simplifies to
\begin{equation}
\pi(k, r) \sim \frac{1}{m}\frac{ay}{y+1} (y+1)^{k-m}e^{-a(y+1)}, \; \; y>0
\end{equation}
which shows that the total probability mass in the range $k=m-O(1)$ and $y>0$ is asymptotically $e^{-a}$.  The remaining mass occurs in the range $r=0$ and $0<x<1$, as the results below show that $\sum_{k=0}^m \pi(k,0) \sim 1-e^{-a}$.

For $k=O(1)$ ($x=O(m^{-1})$) we find that
\begin{equation}
\pi(k, r) \sim m^{k-1/2}{\cal S}_0(y) \frac{y^k}{k!}\left[2-e^{\Psi_y(0,y)}\right]^k e^{m\Psi(0,y)}
\end{equation}
where
\begin{eqnarray}
\Psi(0,y) &=& -\frac{y+\sqrt{y^2+4y}}{y+\sqrt{y^2+4y}+2} + \log\left(\frac{2}{y+\sqrt{y^2+4y}+2}\right), \label{831}\\
\Psi_y(0,y) &=& -\frac{2}{y+\sqrt{y^2+4y}},  \label{867} \\
{\cal S}_0(y) &=&\frac{a\sqrt{2\pi}\left(y+\sqrt{y^2+4y}\right)^{3/2}}{\left(y+\sqrt{y^2+4y}+2\right)\left(y+\sqrt{y^2+4y}+4\right)^{1/2}} \nonumber\\
&\times& \exp\left[-a\left(\frac{y+\sqrt{y^2+4y}}{2}+\frac{1}{y+\sqrt{y^2+4y}+2}\right)\right]. \label{868}
\end{eqnarray}

If $y>0$ the analysis and results for $1-\rho = O(m^{-1})$ are similar to the case $\rho<1$.  However this is not the case if $y=o(1)$.  We give below various results for $r=O(1)$ ($y=O(m^{-1})$) and $r=O(\sqrt{m})$ ($y=O(m^{-1/2})$).  First, for $r=O(1)$ and $0<x<1$ we have
\begin{equation}\label{861}
\pi(k,r) \sim \frac{a}{2}\; m^{r/2-1}r^{-r/2}x^{r/2+1/4}e^{-a(x+1)/2}e^{-2\sqrt{mr}(1-\sqrt{x})},  \; \; r\geq 1,
\end{equation}
and
\begin{equation}\label{876}
\pi(k, 0) -(1-\rho)\rho^k \sim -\frac{a}{2}\frac{x^{3/4}}{\sqrt{m}}e^{-a(x+1)/2}e^{-2\sqrt{m}(1-\sqrt{x})}.
\end{equation}
Note also that with this scaling of $\rho$
\begin{equation}
(1-\rho)\rho^k = \frac{a e^{-ax}}{m}\left[1-\frac{a^2 x}{2m}+O(m^{-2})\right].
\end{equation}
When $k, r = O(1)$ we have, for $r\geq 1$,
\begin{equation}\label{866}
\pi(k,r) \sim  a\sqrt{\pi}\; m^{-5/4}r^{3/4}e^{-(a+r)/2}e^{-2\sqrt{mr}} \frac{1}{2\pi i}\oint \frac{z^{-k-1}}{(1-z)^{r+2}}\; e^{r/(1-z)} dz,
\end{equation}
and
\begin{equation}
\pi(k,0) - (1-\rho)\rho^k \sim  -a\sqrt{\pi}e^{-(a+1)/2} m^{-5/4}e^{-2\sqrt{m}}\; \frac{1}{2\pi i}\oint \frac{z^{-k}e^{1/(1-z)}}{(1-z)^3}\; dz.
\end{equation}
The contour integral may be expressed in terms of a confluent hypergeometric function.

The analysis near the corner $(x, y)=(1, 0)$ is much more complicated in the case $1-\rho = O(m^{-1})$ than when $\rho<1$.  We have obtained some partial results in the corner range, and there are several nested corner layers that must be considered.  First when $1-x = O(m^{-1/2})$ and $y=O(m^{-1/2})$ we let $k=m-\sqrt{m}\xi$ and $r=\sqrt{m}R$.  We find that, for $\xi, R>0$,
\begin{equation}
\pi(k,r) \sim m^{-3/2}\Omega(\xi, R) = m^{-3/2}R^{-1}{\cal D}(\xi, R)
\end{equation}
where ${\cal D}(\xi, R)$ is expressible in terms of ${\cal D}(0, R)$ via the integral
\begin{eqnarray}\label{8g}
{\cal D}(\xi, R) &=&\frac{1}{2\sqrt{\pi}}\int_R^\infty \exp\left[-\frac{(\chi-R+\xi)^2}{4\log(\chi/R)}\right]\nonumber\\
&\times& \left\{\frac{{\cal D}(0, \chi)}{\chi^2\sqrt{\log(\chi/R)}}+\frac{(\chi-R){\cal D}(0, \chi)}{2\chi[\log(\chi/R)]^{3/2}}-\frac{{\cal D}_\chi(0,\chi)}{\chi\sqrt{\log(\chi/R)}}\right\} d\chi\nonumber \\
&+& \frac{\xi}{2\sqrt{\pi}}\int_R^\infty
\frac{{\cal D}(0, \chi)}{2\chi[\log(\chi/R)]^{3/2}}\exp\left[-\frac{(\chi-R+\xi)^2}{4\log(\chi/R)}\right]d\chi
\end{eqnarray}
and ${\cal D}(0, R)$ satisfies the integral equation
\begin{multline}\label{8av}
{\cal D}(0, R) = \frac{1}{\sqrt{\pi}}\int_R^\infty \left\{\frac{{\cal D}(0, \chi)}{\chi^2\sqrt{\log(\chi/R)}}+\frac{(\chi-R){\cal D}(0, \chi)}{2\chi[\log(\chi/R)]^{3/2}}-\frac{{\cal D}_\chi(0,\chi)}{\chi\sqrt{\log(\chi/R)}}\right\} \\
\times \exp\left[-\frac{(\chi-R)^2}{4\log(\chi/R)}\right] d\chi.
\end{multline}
We establish several asymptotic results for ${\cal D}(\xi, R)$ in subsection \textbf{5.4}, and also characterize ${\cal D}$ as the solution to a heat equation with a moving boundary.  In subsection \textbf{5.5} we consider the scale $k=m-O(m^{3/4})$ and $r=O(\sqrt{m})$ and obtain explicit expressions for $\pi(k, r)$,  but these follow simply by expanding the result on the ($x, y$) scale as $(x, y)\rightarrow (1, 0)$ along certain parabolic paths.  In subsection \textbf{5.6} we consider $x\sim 1$ and $r=O(1)$, and find that with
$k=m-\bar{\xi}(\log m)\sqrt{m}$  and $r\geq 1$
\begin{eqnarray}
\pi(k, r) &\sim& \frac{1}{2}ae^{-a}m^{r/2-1}m^{-\bar{\xi}\sqrt{r}}r^{-r/2}, \; \; \bar{\xi} > \sqrt{r}, \\
\pi(k, r) &\sim& \frac{ae^{-a}}{\sqrt{2\pi}}\frac{m^{-\bar{\xi}^2/2}}{m\sqrt{\log m}}
\frac{\Gamma(r-\bar{\xi}^2)}{r!}\; \bar{\xi} ^{2+\bar{\xi}^2}e^{-\bar{\xi}^2}, \; \; 0<\bar{\xi}<\sqrt{r},\\
\pi(k,r) &\sim& \frac{a e^{-a}}{2\sqrt{2\pi}}m^{-r/2-1}r^{-r/2}e^{-\circledast\sqrt{r\log m}} \int_{-\infty}^\circledast e^{-u^2/2}du, \nonumber\\&  & \; \; \;  \; \; \circledast = (\bar{\xi}-\sqrt{r})\sqrt{\log m} = O(1).
\end{eqnarray}
When $r=0$ we obtain
\begin{eqnarray}
\pi(k, 0)-(1-\rho)\rho^k &\sim& -\frac{1}{2}ae^{-a}m^{-\bar{\xi}-1/2}, \; \; \bar{\xi} >1, \\
\pi(k, 0)-\frac{a}{m}e^{-ax} &\sim& \frac{ae^{-a}}{\sqrt{2\pi}}\frac{m^{-\bar{\xi}^2/2}}{m\sqrt{\log m}}\Gamma(1+\bar{\xi}^2)\Gamma(-\bar{\xi}^2)\bar{\xi}^{1+\bar{\xi}^2}, \; \; 0<\bar{\xi}<1,\nonumber\\
\\
\pi(k, 0) - \frac{a}{m}e^{-ax} &\sim& -\frac{ae^{-a}}{2\sqrt{2\pi}}m^{-3/2}e^{-\circledast\sqrt{\log m}}\int_{-\infty}^\circledast e^{-u^2/2} du, \; \; \circledast=O(1). \nonumber\\
\end{eqnarray}
Note that $x=1-\bar{\xi}(\log m)m^{-1/2}$ so that $e^{-ax} \sim e^{-a}\left[1+a\bar{\xi}(\log m)/\sqrt{m}\right]$.

In subsection \textbf{5.6} we shall also discuss the scales $k=m-O(\sqrt{m\log m})$, $\; k=m-O(\sqrt{m})$, and $k=m-O(1)$, but we have not been able to resolve completely all of the corner layer(s) near $(x, y)=(1, 0)$.  Furthermore, the analysis also suggests that we may get different asymptotics for the limit $m\rightarrow \infty$ with $1-\rho = O(m^{-1/2})$, which we do not consider here.

\section{Asymptotic expansions for $m \rightarrow \infty$ with fixed $0< \rho < 1$}
In this section we examine the case when $m \rightarrow \infty$ and $\rho$ is fixed ($0< \rho < 1$).  We set
\begin{equation}
\delta = \frac{1}{m}, \; \; x = \frac{k}{m} = \delta k, \; \;  y = \frac{r}{m} = \delta r
\end{equation}
and then let
\begin{equation}\label{705}
\pi(k, r) = K(x,y;\delta) \exp\left[\frac{1}{\delta} \phi(x,y)\right].
\end{equation}
For $0 < x < 1$, and $\; y > 0$,  (\ref{bl2}) becomes
\begin{multline} \label{71}
(1+\rho)(x+y+\delta) \; K(x, y ;\delta) \exp\left[\frac{1}{\delta} \phi(x,y)\right] \\
=\rho (x+y+\delta)\; K(x-\delta, y ;\delta) \exp\left[\frac{1}{\delta} \phi(x-\delta, y)\right] \\
+(x+\delta)\; K(x+\delta, y ;\delta) \exp\left[\frac{1}{\delta} \phi(x+\delta,y)\right] \\
+(y+\delta)\; K(x, y+\delta ; \delta) \exp\left[\frac{1}{\delta} \phi(x,y+\delta)\right].
\end{multline}
The boundary condition (\ref{ab}) at $x=1$ ($k=m$) becomes
\begin{equation}\label{72}
\rho K(1, y-\delta ; \delta)\exp\left[\frac{1}{\delta}\phi(1, y-\delta)\right]=\frac{1+\delta}{1+y+\delta} \; K(1+\delta, y ; \delta)\exp\left[\frac{1}{\delta}\phi(1+\delta,y)\right].
\end{equation}

We assume that $K(x, y ; \delta)$ has an asymptotic expansion in powers of $\delta$, with
\begin{equation}\label{k}
K(x, y ; \delta) = K(x, y)+\delta K^{(1)}(x, y) + O(\delta^2).
\end{equation}
To compute the leading term for $\pi(k, r)$ we must compute $\phi(x, y)$ and the leading term in (\ref{k}).
We divide (\ref{71}) by $K(x, y)\exp[\phi(x, y)/\delta]$ and let $\delta \rightarrow 0$.  Then we similarly divide (\ref{72}) by $K(1,y)\exp[\phi(1, y)/\delta]$ and let $\delta \rightarrow 0$. This leads to the following 'eikonal' equation for $\phi$, for $0<x<1$ and $y>0$,
\begin{equation}
(1+\rho)(x+y)-\rho(x+y)e^{-\phi_x}-xe^{\phi_x}-ye^{\phi_y} = 0, \label{73}
\end{equation}
and at $x=1$ we have the boundary condition
\begin{equation}
(1+y)e^{-\phi_y} = \frac{e^{\phi_x}}{\rho}. \label{74}
\end{equation}
We solve (\ref{73}) and (\ref{74}) for $\phi(x,y)$ by using the method of
characteristics.  We write (\ref{73}) as $F(x, y, \phi, \phi_x, \phi_y)=0$ where
\begin{equation}
F \equiv (1+\rho)(x+y)-\rho(x+y)e^{-\phi_x}-xe^{\phi_x}-ye^{\phi_y}. \label{73F}
\end{equation}
The characteristic equations for this nonlinear PDE are \cite{CH}
\begin{eqnarray}
\dot{x}&=& \frac{dx}{dt} = \frac{\partial F}{\partial\phi_x}
=\rho(x+y)e^{-\phi_x}-xe^{\phi_x}, \label{75}\\
\dot{y} &=& \frac{dy}{dt} = \frac{\partial F}{\partial\phi_y} = -ye^{\phi_y}, \label{76}\\
\dot{\phi}&=& \frac{d\phi}{dt} = \phi_x \dot{x}+\phi_y \dot{y}, \label{77}\\
\dot{\phi_x} &=& -\frac{\partial F}{\partial x} =
-(1+\rho)+\rho \; e^{-\phi_x}+e^{\phi_x}, \label{78}\\
\dot{\phi_y} &=& -\frac{\partial F}{\partial y} =
-(1+\rho)+\rho \; e^{-\phi_x}+e^{\phi_y}. \label{79}
\end{eqnarray}
To solve this system we use rays starting from $(x, y)=(1, s)$ at $t=0$.  Thus the 'initial manifold' is $x=1$ where the boundary condition in (\ref{74}) applies, and all the rays start from $x=1$ at $t=0$. We view $x$ and $y$ as functions of $s$ and $t$.
We set $\phi_x \equiv \log \Upsilon(s,t)$ and, from (\ref{73F}), (\ref{75}), (\ref{76}) and (\ref{78}), obtain
\begin{eqnarray}
\frac{\dot{x}+\dot{y}}{x+y} = \frac{2\rho}{\Upsilon}-(1+\rho), \label{710} \\
\frac{\partial\Upsilon}{\partial t} = (\Upsilon-\rho)(\Upsilon-1). \label{711}
\end{eqnarray}
Solving (\ref{711}) gives
\begin{equation} \label{712}
\Upsilon= \frac{C_0(s) e^{(\rho-1)t}-\rho}{C_0(s) e^{(\rho-1)t}-1}
\end{equation}
and using (\ref{712}) we solve (\ref{710}) to obtain
\begin{equation}
x+y=C_1(s)[C_0(s) e^{(\rho-1)t}-\rho]^2 e^{(1-\rho)t}. \label{713}
\end{equation}
The functions $C_0(s)$ and $C_1(s)$ will be determined shortly.

From (\ref{712}) and (\ref{713}), (\ref{75}) can be written as
\begin{equation}
\frac{\partial x}{\partial t}+\left[\frac{C_0(s) e^{(\rho-1)t}-\rho}{C_0(s) e^{(\rho-1)t}-1}\right] x =\rho C_1(s)[C_0(s) e^{(\rho-1)t}-\rho] [C_0(s) -e^{(1-\rho)t}]. \label{714}
\end{equation}
Solving (\ref{714}) we obtain $x(s,t)$ as
\begin{equation}\label{702}
x(s,t) = [C_0(s) C_1(s)-\rho^2 C_1(s) e^{(1-\rho)t} + C_2(s) e^{-\rho t}][C_0(s) e^{(\rho-1)t}-1],
\end{equation}
and then $y(s, t)$ follows from (\ref{713}) as
\begin{equation}\label{703}
y(s,t) = C_0(s) C_1(s)(1-\rho)^2 + C_2(s)[e^{-\rho t}-C_0(s) e^{-t}].
\end{equation}
Applying the initial conditions ($x(s,0)=1, \;  y(s,0)=s$) we find that
\begin{eqnarray}
C_1(s) &=& \frac{s+1}{[C_0(s)-\rho]^2}, \\
C_2(s) &=& \frac{1}{C_0(s)-1}-\frac{C_0(s)-\rho^2}{[C_0(s)-\rho]^2}(s+1). \label{704}
\end{eqnarray}

To determine $C_0(s)$ we set $\phi_x=A(s)$ and $\phi_y=B(s)$ at $t=0$, and from (\ref{73}) and (\ref{74}) obtain the following equations along the initial manifold
\begin{eqnarray}
(1+\rho)(1+s)&=&\rho(1+s)e^{-A}+e^{A}+se^{B}, \label{70}\\
\rho(1+s)e^{-A} &=& e^{B}. \label{701}
\end{eqnarray}
Solving (\ref{70}) and (\ref{701}) leads to
\begin{equation*}
A(s) = \log(s+1), \; \; B(s)=\log\rho.
\end{equation*}
Since $e^{\phi_y}= -\dot{y}/y$ from (\ref{76}), $\phi_y(1,s)=\log[-\dot{y}(s,0)/s]=B(s)=\log\rho$.  Thus we obtain $C_0(s)= (s+1-\rho)/s$ and then, from (\ref{702})-(\ref{704}), obtain (\ref{x}) and (\ref{y}).
Solving (\ref{78}) and (\ref{79}) then yields (\ref{phix}) and (\ref{phiy}).
Finally, we obtain $\phi(x,y)$ by integrating (\ref{77}) with respect to $t$. Using integration by parts, we can write
\begin{equation}
\phi = \int (\phi_x \dot{x} + \phi_y \dot{y}) dt = \phi_x x+\phi_y y - \int (\dot{\phi_x}x+\dot{\phi_y}y) dt. \label{723}
\end{equation}
The integrand in the last integral in (\ref{723})  is $-F = 0$.  We can also check from (\ref{phix}) - (\ref{y}) that indeed $\dot{\phi_x}/\dot{\phi_y} = -y/x$.  Therefore, $\phi(x,y) = \phi_x x+\phi_y y+f(s)$ where $f(s)$ is some function of $s$.  At $t=0$,  $ \; \phi(1,s) = \log(s+1)+s\log\rho+f(s)$ and since $y=s$ at $t=0$,
\begin{equation*}
\frac{\partial \phi(1,s)}{\partial s} = \frac{1}{s+1} + \log \rho +f'(s) = \phi_y(1,s) = \log \rho.
\end{equation*}
Thus $f(s) = -\log (s+1)+ constant$.  We will later show, by asymptotic matching, that near the corner $(x, y)=(1, 0)$ the solution $\pi(k,r)$ must be $O(\rho^m)=O\left(\exp[(\log \rho)/\delta]\right)$. Thus $\phi(1,0) = \log\rho$ and we then obtain (\ref{phi}).
The value of the $constant$ can also be ultimately found by normalization or by using (\ref{29}). Note, however, that our expansion (\ref{705}) applies only in the domain $0<x\leq1, \; y>0$.  We shall construct different expansions for $x\approx0$ and $y\approx0$, and also near the two corners $(x,y)=(0,0)$ and $(1,0)$.

We can plot the rays, which are given by (\ref{x}) and (\ref{y}) in parametric forms. Each fixed value of $s$ corresponds to a particular ray parameterized by $t$.  We sketch the rays in \textbf{Fig. 2} and \textbf{Fig. 3}.  In \textbf{Fig. 2} we plot the rays for $s>0$ and $t$ in the range  $0 \leq t \leq t_{max}\equiv \log [1+(1-\rho)/s]/(1-\rho)$.  At $t=0$ we have $x=1$ and at $t=t_{max}(s)$ the rays reach $x=0$.  We also note that, in view of (\ref{phix}), $\phi_x$ develops a singularity at $t=t_{max}(s)$.  However, using (\ref{x}) and (\ref{y}) we can continue the rays for $t>t_{max}$, and this continuation is sketched in \textbf{Fig. 3}.  We see that each ray reaches a minimum value of $x$, where it has a cusp.  The cusp occurs when $\dot{x} =\dot{y}=0$, which happens at $t=t_c(s)\equiv \log[(s+1-\rho)/(\rho s)]/(1-\rho) > t_{max}(s)$.  From (\ref{x}) and (\ref{y}) we also see that when $t=t_c, \; x+y =0$, so that the line $y=-x$ is the locus of the cusp points.  For $t>t_c$ we again have $\dot{x}>0$ and the rays eventually re-enter the domain $x>0$, and then $x$ increases past $x=1$.  The full rays are sketched in \textbf{Fig. 3}, up to the time they return to $x=1$.  The figure clearly shows the cusps and their locus.  However, since the cusps occur outside of the domain of interest $\{0\leq x\leq 1, \; y\geq 0\}$  they do not affect the solution very much, and we only consider the rays for $0\leq t\leq t_{max}$.  The origin $(0,0)$ is the only cusp point in the domain, and this will affect significantly the asymptotic structure of $\pi(k,r)$ for the scale $k, r=O(1)$.  We analyze this range in subsection \textbf{4.4}.

\begin{center}
\includegraphics[angle=0, width=0.6\textwidth]{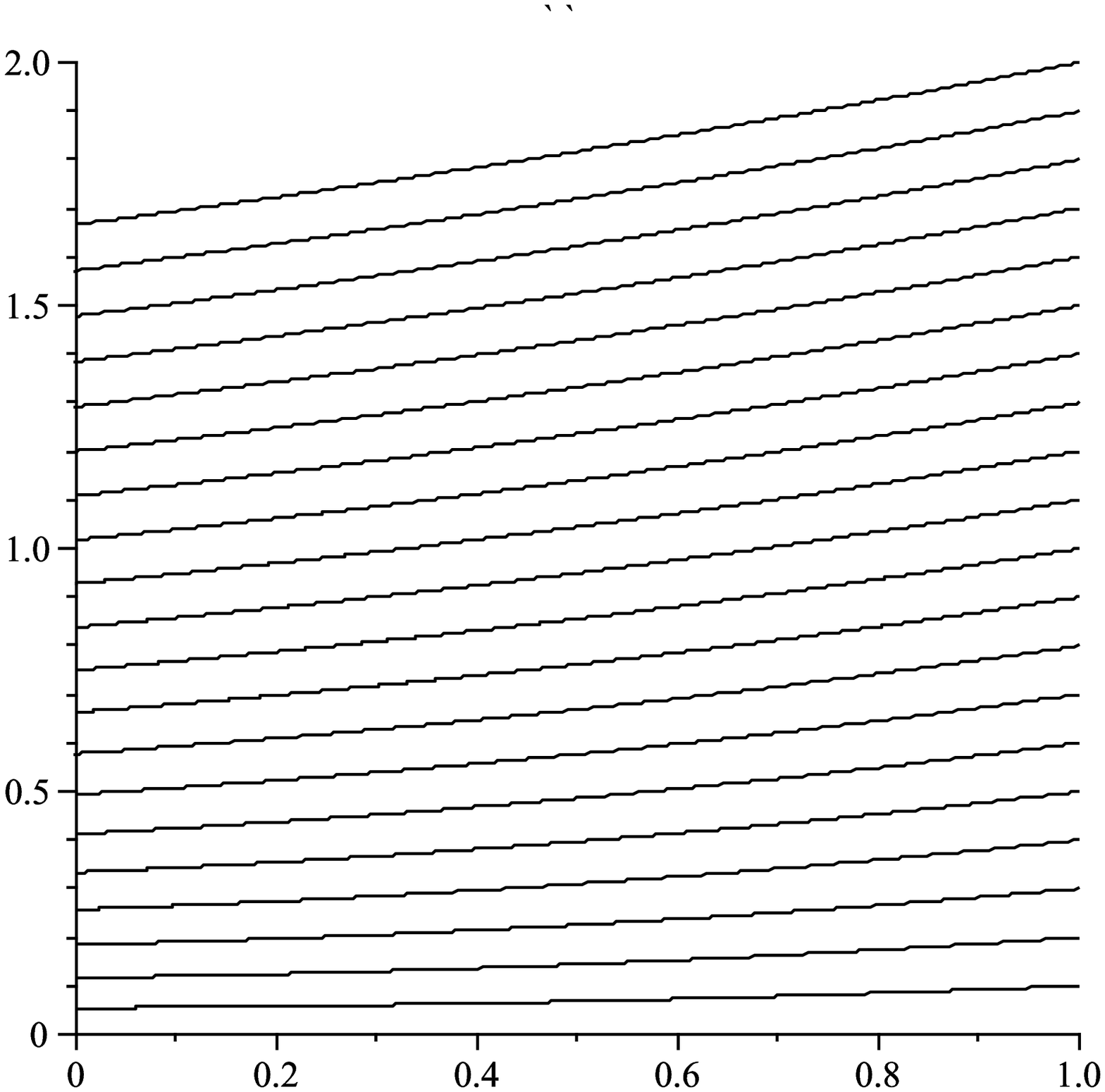}\\
\textbf{Fig. 2}   A sketch of the rays from $x=1$.
\end{center}

\begin{center}
\includegraphics[angle=0, width=0.6\textwidth]{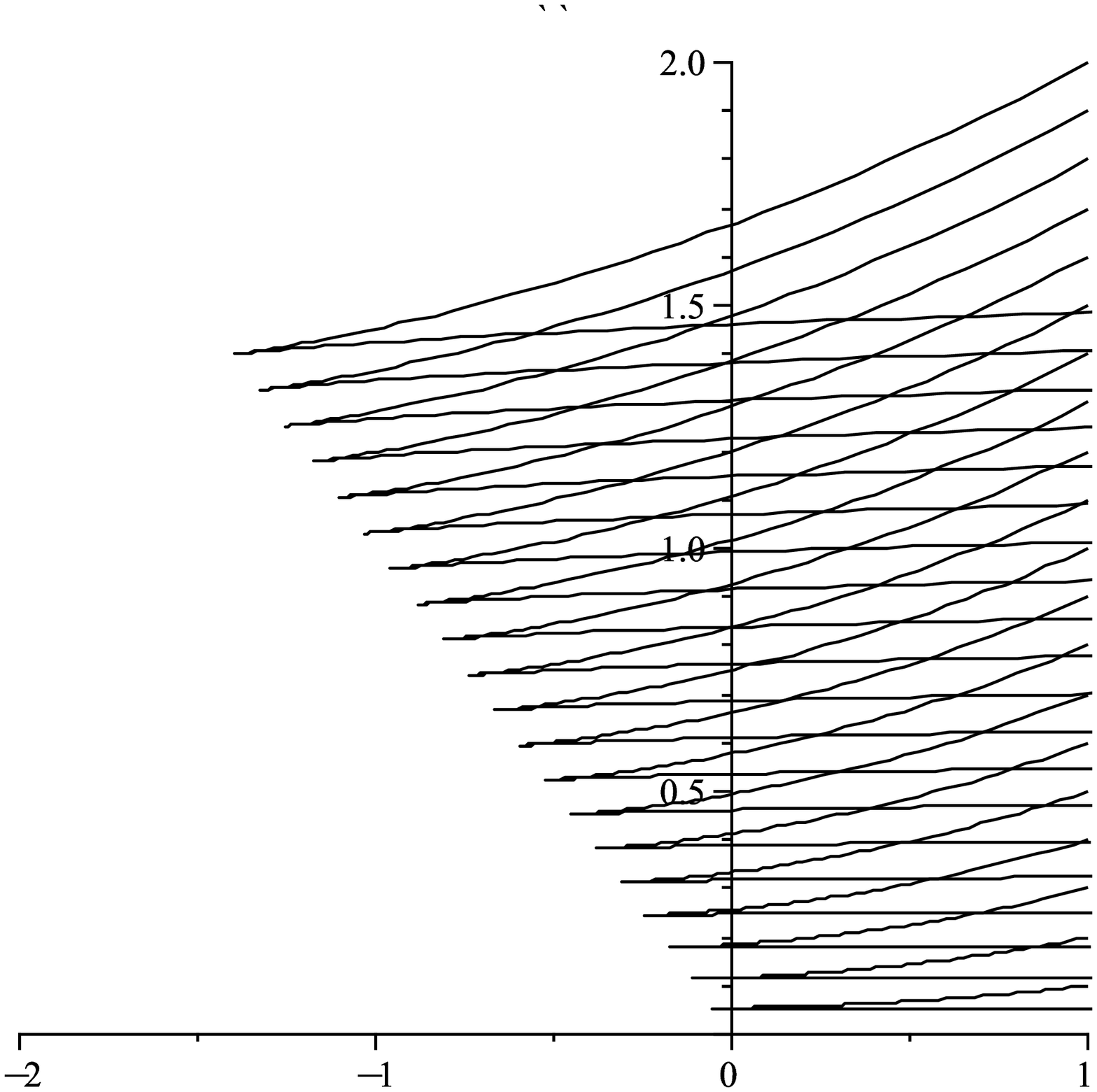}\\
\textbf{Fig. 3}   A sketch of the extended rays from $x=1$.
\end{center}

We next calculate $K(x,y)$ in (\ref{k}).  Expanding (\ref{71}) as $\delta\rightarrow 0$ and using the fact that $\phi$ satisfies (\ref{73}) we obtain the following 'transport' equation.
\begin{multline}
\left[\left(\frac{x+y}{2}\phi_{xx}+1\right)\rho e^{-\phi_x} + \left(\frac{x}{2} \phi_{xx} + 1\right)e^{\phi_x} +\left(\frac{y}{2} \phi_{yy}+1\right)e^{\phi_y} -(1+\rho)\right] K \\
 = \left[\rho(x+y) e^{-\phi_x}-x e^{\phi_x}\right]K_x-y e^{\phi_y} K_y. \label{719}
\end{multline}
We note that the right hand side is $\dot{x}K_x+\dot{y}K_y = \dot{K}$, which is the directional derivative of $K$ along a ray.  The factor that multiplies $K$ in the left hand side of (\ref{719}) may be rewritten as
\begin{multline}
-\frac{1}{2}\frac{\partial}{\partial x} [\rho(x+y)e^{-\phi_x}-x e^{\phi_x}] +\frac{1}{2}\frac{\partial}{\partial y}(y e^{\phi_y})-(1+\rho)+\frac{1}{2}(e^{\phi_x}+e^{\phi_y})+\frac{3}{2}\rho e^{-\phi_x}\\
=-\frac{1}{2}\frac{\partial \dot{x}}{\partial x} -\frac{1}{2}\frac{\partial \dot{y}}{\partial y}-(1+\rho)+\frac{1}{2}(e^{\phi_x}+e^{\phi_y})+\frac{3}{2}\rho e^{-\phi_x}. \label{715}
\end{multline}
Defining $\jmath(s,t)$ to be $x_t y_s-x_s y_t$, which is the Jacobian of the transformation from $(x,y)$ to $(s, t)$ coordinates, we find that
\begin{eqnarray}
\frac{\partial \dot{x}}{\partial x} &=& \frac{\partial x_t}{\partial x} = x_{tt}t_x+x_{ts}s_x=\frac{x_{tt}y_s-x_{ts}y_t}{\jmath}, \label{716}\\
\frac{\partial \dot{y}}{\partial y} &=& \frac{\partial y_t}{\partial y} = y_{tt}t_y+y_{ts}s_y=\frac{y_{ts}x_t-y_{tt}x_s}{\jmath}, \label{717}
\end{eqnarray}
and
\begin{equation}
\frac{\partial \jmath}{\partial t} = \dot{\jmath} = x_{tt}y_s+y_{ts}x_t-x_{ts}y_t-y_{tt}x_s. \label{718}
\end{equation}
Thus, from (\ref{78}), (\ref{79}), and (\ref{715}) - (\ref{718}), we can write (\ref{719}) as
\begin{equation}\label{784}
\frac{\dot{K}}{K} = -\frac{1}{2}\frac{\dot{\jmath}}{\jmath}+\frac{1}{2}(\dot{\phi_x}+\dot{\phi_y})+\frac{1}{2}\rho e^{-\phi_x}.
\end{equation}
Integrating (\ref{784}) using (\ref{phix}) for $\phi_x$ we obtain
\begin{equation} \label{726}
K(x,y)= K_0(s)\left[ \frac{|(s+1-\rho)e^{-(1-\rho)t}-s\rho|}{|\jmath|}\right]^{1/2} \exp\left(\frac{\phi_x+\phi_y+t}{2}\right)
\end{equation}
where $K_0(s)$ is a function of $s$.

To determine $K_0(s)$ we use the boundary condition at $x=1$.  From (\ref{72}) we obtained (\ref{74}) to leading order in $\delta$, and at the next order we obtain
\begin{multline}
\rho e^{-\phi_y}\left[\frac{(y+1)}{2}\phi_{yy} K(1,s)+K(1,s)-(y+1)K_y(1,s)\right] \\
= e^{\phi_x}\left[\left(\frac{\phi_{xx}}{2}+1\right)K(1,s)+K_x(1,s)\right]. \label{720}
\end{multline}
Since $\phi_{xx}(1,s)=s/(\rho-s-1)$ and $\phi_{yy}(1,s)= 0$, (\ref{720}) becomes
\begin{equation}
(s+1)[K_x(1,s)+K_y(1,s)]+\frac{s(2\rho-1-s)}{2(\rho-s-1)}K(1,s) = 0. \label{721}
\end{equation}
Using the relations
\begin{eqnarray*}
K_x &=& K_t t_x + K_s s_x = K_t \frac{y_s}{\jmath}-K_s\frac{y_t}{\jmath}= \frac{K_t+\rho s K_s}{\rho-s-1},\\
K_y &=& K_t t_y + K_s s_y =-K_t \frac{x_s}{\jmath}+K_s \frac{x_t}{\jmath} = K_s,\\
\jmath(s,0)&=& \rho-1-s,
\end{eqnarray*}
(\ref{721}) becomes
\begin{equation}
\frac{K_t}{K}+(\rho-1) (s+1)\frac{K_s}{K}+\frac{\rho s}{1+s}-\frac{s}{2} = 0, \;  \; t=0.\label{722}
\end{equation}
Evaluating (\ref{719}) at $t=0$ yields
\begin{equation} \label{724}
\frac{K_t}{K} = \frac{s(\rho+s+1)}{2(\rho-s-1)}+\frac{\rho}{s+1} +s
\end{equation}
and thus, from (\ref{722}) and (\ref{724}), we obtain the following ODE for $K$ at $t=0$ with respect to $s$.
\begin{equation}\label{785}
\frac{K_s}{K} = \frac{1}{s+1-\rho}-\frac{1}{s+1}.
\end{equation}
Integrating (\ref{785}) gives
\begin{equation} \label{725}
K(1,s) = C^*\left(1-\frac{\rho}{s+1}\right)
\end{equation}
where $C^*$ is a constant.  Equating (\ref{726}) at $t=0$ with (\ref{725}) we obtain $K_0(s)$ as
\begin{equation*}
K_0(s) = C^*\frac{(s+1-\rho)^{3/2}}{(s+1)^2}\frac{1}{\sqrt{\rho(1-\rho)}}.
\end{equation*}

To obtain $C^*$ we use the fact that $N_1+N_2$ follows the geometric distribution in (\ref{29}).  Expanding $\phi$ and $K$ near $x=1$, which corresponds to $t=0$, the approximation $\pi(k,r)\sim K(x, y)\exp[\phi(x, y)/\delta]$ simplifies to
\begin{eqnarray}\label{751}
\pi(k,r)\sim C^*\left(1-\frac{\rho}{y+1}\right)\exp\left(\frac{y+1}{\delta}\log\rho\right)\exp\left[\frac{x-1}{\delta}\log(y+1)\right].
\end{eqnarray}
Here we used $K(x, y)\sim K(1, y)$ and $\phi(x, y) = \phi(1, y)+\phi_x (1, y)(x-1)+\cdot\cdot\cdot\cdot\cdot.$  If $x+y$ is constant, we set $x+y=z>1$ and use (\ref{29}) with $N=z/\delta=O(\delta^{-1})$.  We sum over $(k,r) = (m, N-m), (m-1, N-m+1), (m-2, N-m+2), \cdot\cdot\cdot$, which corresponds to $(x, y)=(1, z-1), (1-\delta, z-1+\delta), (1-2\delta, z-1+2\delta), \cdot\cdot\cdot.$  Hence asymptotically for $\delta = 1/m \rightarrow 0\;$  (\ref{29}) becomes
\begin{eqnarray*}
\sum_{k+r=z/\delta}\pi(k,r)&=&(1-\rho)\rho^{z/\delta} \\
&\sim & C^*\left[1-\frac{\rho}{z+O(\delta)}\right]\sum_{j=0}^\infty \exp\left[\left(\frac{z}{\delta}+j\right)\log\rho\right]\exp\left[-j\log(z+O(\delta))\right] \\
&\sim & C^*\left(1-\frac{\rho}{z}\right)\rho^{z/\delta}\sum_{j=0}^\infty \left(\frac{\rho}{z}\right)^j = C^* \rho^{z/\delta}.
\end{eqnarray*}
Therefore, $C^*=1-\rho$ and thus we obtain (\ref{769}).
We also obtain $\pi(k,r)$ near $x=1$ from (\ref{751}) as (\ref{727}), which applies for $k=m-O(1)$ and $y>0$.

The ray solution $\pi(k, r) \sim K(x,y) \exp[m \phi(x,y)]$ is valid in the interior of the state space and also near $x=1$.  However, the expansion breaks down near the boundaries $x=0$ and $y=0$, and near the corner points $(x, y)=(1, 0)$ and $(0,0)$.  We proceed to analyze  these regions separately.

\subsection{Boundary layer near $x=0$}
We observed from (\ref{phix}) that $\phi_x$ has a logarithmic singularity as $x\rightarrow 0$. We consider the scale $k=O(1)$ and $y>0$, which corresponds to $x=O(\delta)$ and $r=O(\delta^{-1})= O(m)$.  We set
\begin{equation*}
\pi(k,r) = \delta^{\nu-k}\exp\left[\frac{1}{\delta}\phi(0,y)\right]S_k(y; \delta),
\end{equation*}
and rewrite (\ref{bl2}) as
\begin{multline}\label{728}
(1+\rho)[y+\delta(k+1)]\; e^{\phi(0,y)/\delta}S_k(y; \delta) = \delta\rho [y+\delta(k+1)]e^{\phi(0,y)/\delta}S_{k-1}(y; \delta)\\
+(k+1)\; e^{\phi(0,y)/\delta}S_{k+1}(y; \delta)
 +(y+\delta)\; e^{\phi(0,y+\delta)/\delta}S_k(y+\delta; \delta)
\end{multline}
for $k>0$ and $y > 0$.  The boundary condition at $k=0$ is
\begin{equation}\label{729}
(1+\rho)(y+\delta)\; e^{\phi(0,y)/\delta}S_0(y; \delta)= e^{\phi(0,y)/\delta}S_1(y; \delta) + (y+\delta)e^{\phi(0,y+\delta)/\delta}S_0(y+\delta; \delta).
\end{equation}
We expand $S_k(y; \delta)$ in the form
$$S_k(y; \delta) = S_k(y)+\delta S_k^{(1)}(y)+O(\delta^2).$$
In (\ref{728}) and (\ref{729}) we let $\delta\rightarrow 0$  and obtain the following equation for the leading term $S_k(y)$
\begin{eqnarray} \label{731}
\left[1+\rho-e^{\phi_y(0,y)}\right]y &=& (k+1)\; \frac{S_{k+1}(y)}{S_k(y)},  \;  \;  k \geq 0.
\end{eqnarray}
The general solution to (\ref{731}) is
\begin{equation*}
S_k(y) = S_0(y) \frac{y^k}{k!}\left[1+\rho-e^{\phi_y(0,y)}\right]^k
\end{equation*}
and thus
\begin{equation} \label{732}
\pi(k,r) \sim \delta^{\nu-k}e^{\phi(0,y)/\delta}S_0(y) \frac{y^k}{k!}\left[1+\rho-e^{\phi_y(0,y)}\right]^k.
\end{equation}
It remains to determine the constant $\nu$ and the function $S_0(y)$.

We next asymptotically match (\ref{732}) to the ray expansion $K(x, y)\exp[\phi(x, y)/\delta]$, in an intermediate limit where $x\rightarrow 0$ but $k=x/\delta \rightarrow \infty$. Using Stirling's formula for $k!$ the expansion in (\ref{732}) for $k\rightarrow \infty$ becomes
\begin{eqnarray}\label{795}
\pi(k,r) &\sim& \delta^{\nu}e^{\phi(0,y)/\delta}S_0(y) \frac{e^k}{\sqrt{2\pi k}}\exp[-k\log k -k \log\delta]\left[1+\rho-e^{\phi_y(0,y)}\right]^k y^k  \nonumber\\
&=& \delta^{\nu+1/2} \frac{S_0(y)}{\sqrt{2\pi x}}\exp\left[\frac{1}{\delta}\left\{\phi(0,y)+x+x\log \left(\frac{y}{x}\right)+x\log\left[1+\rho-e^{\phi_y(0,y)}\right]\right\}\right].\nonumber\\
\end{eqnarray}
Here we rewrote the result in terms of $x$.

We expand the ray solution as $x\rightarrow 0$.  Along $x=0$ we can explicitly invert the transformation from $(x, y)$ to ray coordinates.  When $x=0$ we have $t=t_{max}(s)$ and $s+1-\rho-s e^{(1-\rho)t}=0$.  Then from (\ref{y}) we find that $s$ and $y$ are related by $y=s(s+1-\rho)/(s+1).$ We thus define $s_*=s_*(y)$ by (\ref{734}).
Thus, a ray that starts from $(x, y)=(1, s_*)$ hits the $y-$axis at the point $(0, y)$.  Then we use (\ref{phi}) and (\ref{phiy}) and evaluate these expressions at $t=t_{max}$ and $s=s_*$, to obtain explicit expressions (\ref{760}) and (\ref{761}) for $\phi(0,y)$ and $\phi_y(0,y)$.
We note that $\phi_y(0,y)=\log \rho -t_{max}(s_*)$.  Also, from (\ref{phix}) and (\ref{x}) we obtain
\begin{equation}\label{788}
\phi_x = -\log x + \log \left(\frac{\left[s+1-\rho-\rho s \; e^{(1-\rho)t}\right]\left[(s+1-\rho)\; e^{-(1-\rho)t}-s\rho^2-\rho(1-\rho)s e^{-t}\right]}{(s+1)(1-\rho)^2}\right),
\end{equation}
which is exact for all $(x, y)$ and indicates the logarithmic singularity in $\phi_x$ as $x\rightarrow 0$.  By expanding (\ref{788}) for $s\rightarrow s_*$ and $t\rightarrow t_{max}$ we obtain
\begin{eqnarray}
\phi_x &=& -\log x + \log \left[\frac{s_*(s_*+1-\rho)(1+\rho-\rho e^{-t})}{(s_*+1)}\right]+o(1) \nonumber\\
&=&-\log x + \log y +\log\left[1+\rho-e^{\phi_y(0,y)}\right]+o(1),
\end{eqnarray}
which integrates to
\begin{equation}
\phi(x, y) = \phi(0, y)-x\log x +x+x\log y +x \log\left[1+\rho-e^{\phi_y(0,y)}\right]+o(x).
\end{equation}
The above precisely agrees with the exponential terms in (\ref{795}).  The expansions will thus asymptotically match provided that
$\nu = -1/2$ and
\begin{equation*}
S_0(y) = \lim_{x \rightarrow 0}\left[\sqrt{2\pi x} \; K(x,y)\right]
\end{equation*}
By expanding $\jmath=x_t y_s-x_s y_t$ and then $K$ in (\ref{769}) for $s\rightarrow s_*$ and $t\rightarrow t_{max}$, we obtain after some calculation
\begin{eqnarray}
S_0(y) &=& \frac{\sqrt{2\pi} (1-\rho)(s_*+1-\rho)^{3/2}\sqrt{s_*}}{(s_*+1)\sqrt{(s_*+1)^2-\rho}} \nonumber\\
&=&\frac{(1-\rho)\sqrt{2\pi y}\;
\left[(1+\rho)y+(1-\rho)\sqrt{(y-1+\rho)^2+4y}+(1-\rho)^2\right]^{1/2}}{\left[y^2+2(1+\rho)y+(y+1+\rho)\sqrt{(y-1+\rho)^2+4y}+(1-\rho)^2\right]^{1/2}}, \nonumber\\
 \label{759}
\end{eqnarray}
This completes the determination of $\nu$ and $S_0(y)$ in (\ref{732}), and verifies the matching between the $(x, y)$ and $(k, y)$ scales.

\subsection{Boundary layer near $y=0$}
We consider small values of $y$, and return to the original discrete
variable $r=y/\delta$, thus setting
\begin{equation}\label{7201}
\pi(k,r) = \Pi_r(x;\delta) \; \; r\geq 1
\end{equation}
and
\begin{equation}\label{7202}
\pi(k,0) = (1-\rho)\rho^k+\Pi_0(x;\delta).
\end{equation}
We assume for now that $0<x<1$, and will treat the corner regions,
where $x\approx 0$ or $x \approx 1$, separately.  Also we expect
that for large $m$ and $\rho<1$, the secondary servers will rarely
be needed.  Thus we expect that $\pi(k,r)$ will be mostly
concentrated along $\pi(k,0)$, with this function being roughly
geometric in $\rho$.  Thus we wrote $\pi(k,0)$ in a form so as to
estimate its deviation from a geometric distribution.  We also note
that $(1-\rho)\rho^k I_{\{r=0\}}$ is an exact solution to
(\ref{bl2}), and also satisfies the corner condition (\ref{26}) at
$(k, r)=(0, 0)$.  It fails to satisfy the full problem in (\ref{bl})
only due to the boundary condition at $k=m$.  With (\ref{7201}) and
(\ref{7202}), (\ref{bl2}) becomes
\begin{eqnarray}
(1+\rho)\; \Pi_r(x; \delta) &=& \rho \; \Pi_r (x-\delta; \delta)+\frac{x+\delta}{x+\delta(r+1)}\; \Pi_r (x+\delta; \delta) \nonumber\\
 &+&\frac{\delta(r+1)}{x+\delta(r+1)}\; \Pi_{r+1}(x; \delta),       \;    \;   r > 0. \label{735} \\
(1+\rho)\; \Pi_0(x; \delta) &=& \rho \; \Pi_0 (x-\delta; \delta)+
\Pi_0 (x+\delta; \delta)
 +\frac{\delta}{x+\delta}\; \Pi_1 (x; \delta), \;  \;  r=0. \label{736}\nonumber\\
\end{eqnarray}
We expand $\Pi_r(x; \delta)$ as
\begin{equation*}
\Pi_r(x; \delta) =  \Pi_r(x)+\delta \Pi_r^{(1)}(x)+O(\delta^2).
\end{equation*}
From (\ref{735}) and (\ref{736}) we obtain to leading order the following
equations
\begin{eqnarray}
x(1-\rho)\Pi_r'-r\Pi_r+(r+1)\Pi_{r+1} &=& 0, \label{739}\\
x(1-\rho)\Pi_0'+\Pi_1 &=& 0. \label{740}
\end{eqnarray}
Equations (\ref{739}) and (\ref{740}) express the 'current' value $\Pi_r(x)$ in terms of the 'future' value $\Pi_{r+1}(x)$. However, for sufficiently large $r$ we can use the ray expansion to compute $\pi(k,r)$ asymptotically.  Expanding $K(x, y)\exp[m\phi(x, y)]$ for $y\rightarrow 0$ and $0<x<1$ leads to
\begin{eqnarray}
\pi(k,r) &=&  K(x,0) \exp\{m[\phi(x,0)+y\phi_y(x,0)+O(y^2)]\} \nonumber\\
&\sim& \frac{(1-\rho)^2}{1-\rho+\rho e^{-t}}\;  \rho^m [e^{\phi_y(x,0)}]^r \nonumber\\
&=& \frac{(1-\rho)^2 \rho^{m+r} (e^{-t})^r}{[1-\rho+\rho e^{-t}]^{r+1}}. \label{737}
\end{eqnarray}
Here we used the fact that  $y=0$ implies that $s=0$, and thus (cf.  (\ref{x})) $x(0,t) = e^{(\rho-1)t}$ and $e^{-t} = x^{1/(1-\rho)}$. The ray that starts from the corner $(1, 0)$ is the line segment $y=0, \; 0\leq x \leq 1$ and this ray reaches $(x, 0)$ when $t=-(\log x) /(1-\rho)$. Therefore, from (\ref{737}), as $y\rightarrow 0 \;  \; \pi(k,r)$ becomes (\ref{738}).
We can easily check that (\ref{738}) satisfies (\ref{739}). Thus we have the leading term $\Pi_r(x)$ for $r\geq 1$.  To obtain $\Pi_0(x)$ we solve (\ref{740}) with $\Pi_1(x)$ computed from the right hand side of (\ref{738}) with $r=1$.  We thus obtain
\begin{equation*}
\Pi_0(x) =  \rho^m \left[\tilde{C}+\frac{(1-\rho)^2}{1-\rho+\rho x^{1/(1-\rho)}}\right]
\end{equation*}
where $\tilde{C}$ is a constant. Since $\pi(0,0)=1-\rho$, we expect that $ \; \Pi_0(0) = \rho^m [\tilde{C}+(1-\rho)]=0$, which implies that $\tilde{C}=-(1-\rho)$.  Therefore, for $r=0$ we have (\ref{744}).
In subsection \textbf{4.4} we will consider $\pi(k,r)$ for $r$ and $k=O(1)$, and will give a more precise argument, based on asymptotic matching, to show that $\Pi_0(0)=0$.

We know that for $k\leq m, \; \sum_{n=0}^m \pi(k-n, n) = (1-\rho)\rho^k$ exactly, and this implies that $\sum_{J=0}^{x/\delta} \Pi_J(x-\delta J; \delta)=0$.
Since $0<x<1$, asymptotically we should have $\sum_{J=0}^\infty \Pi_J(x)=0$, which is indeed true since
\begin{equation*}
\sum_{r=1}^\infty \Pi_r(x) = \sum_{r=1}^\infty \frac{(1-\rho)^2 \rho^{m+r} x^{r/(1-\rho)}}{[1-\rho+\rho x^{1/(1-\rho)}]^{r+1}} = \frac{(1-\rho)\rho^{m+1} x^{1/(1-\rho)}}{1-\rho+\rho x^{1/(1-\rho)}}=-\Pi_0(x).
\end{equation*}

\subsection{Corner layer near $(x,y)=(1,0)$}
We consider $(x, y)$ near the corner $(1, 0)$.  We again use the discrete variable $r$ and consider $k=m-O(1)$ (or $1-x=O(\delta)$).  Thus we define $n$ and $L$ by
\begin{equation*}
n=m-k, \;\;  \pi(k,r) = L(n,r; m).
\end{equation*}
The main balance equation (\ref{bl2}), written in terms of $(n, r)$, becomes
\begin{multline} \label{741}
(1+\rho) L(n,r; m) = \rho L(n+1,r; m) + \frac{m-n+1}{m-n+1+r} L(n-1,r; m) \\
+\frac{r+1}{m-n+1+r} L(n,r+1; m).
\end{multline}
Here we used $\pi(k\pm 1, r)=L(n\mp 1, r; m)$.  The artificial boundary condition (\ref{ab}) becomes
\begin{equation}\label{706}
\frac{m+1}{m+r+1}L(-1, r; m) =\rho\;L(0, r-1; m), \; \; r\geq 2,
\end{equation}
and the corner condition in (\ref{27}) is
\begin{equation}\label{707}
(1+\rho)L(0,0; m) = \rho L(1,0; m)+\frac{1}{m+1} L(0,1; m).
\end{equation}
By requiring (\ref{741}) to hold also at $n=0$, thus defining $L(-1, r; m)$, we can use (\ref{706}) (or (\ref{ab})) instead of (\ref{bl2}).  Then (\ref{707}) may be replaced by
\begin{equation}\label{708}
L(-1, 0; m)=0.
\end{equation}
The condition at the other corner $(0,0)$ corresponds asymptotically to $n=\infty$ and will play no role.

We expand $L$ for $m\rightarrow \infty$ (or $\delta\rightarrow 0$) with
\begin{equation*}
L(n, r; m) =  L(n, r)+\frac{1}{m} L^{(1)}(n, r)+O(m^{-2}).
\end{equation*}
To leading order we obtain from (\ref{741})
\begin{equation}\label{743}
(1+\rho) L(n,r)= \rho L(n+1,r) +  L(n-1,r), \; r \geq 0
\end{equation}
and (\ref{706}) and (\ref{708}) lead to
\begin{eqnarray}
L(-1,r) &=& \rho \; L(0,r-1), \; \; r\geq 2, \label{700}\\
L(-1, 0) &=& 0. \label{746}
\end{eqnarray}

We can infer the solution of (\ref{743})-(\ref{746}) be expanding the ray solution as $(x, y)\rightarrow (1,0)$.  This can be obtained by simply letting $x\rightarrow 1$ in (\ref{738}), which suggests that
\begin{equation}
L(n,r) = (1-\rho)^2 \rho^{m+r}, \; \; r\geq 1, \label{748}
\end{equation}
which is independent of $n$.  We can easily verify that (\ref{748}) satisfies (\ref{743}) and (\ref{700}).  To obtain $L(n, 0)$ we let $x\rightarrow 1$ in (\ref{744}) and recall that $k=m-n$.  Hence,
\begin{equation}\label{709}
L(n,0) = (1-\rho)(\rho^{-n}-\rho)\rho^m.
\end{equation}
We also note that (\ref{707}) implies asymptotically that
\begin{equation*}
(1+\rho)L(0,0)=\rho L(1,0),
\end{equation*}
which is indeed satisfied by (\ref{709}).  Thus we have obtained the leading term for $\pi(k,r)$ in the corner region in a very simple form, with (\ref{742}) and (\ref{749}).
We note that when $r=0$ and $k=m-O(1)$ even the leading term indicates a deviation from the geometric distribution $(1-\rho)\rho^k=(1-\rho)\rho^{m-n}$.

As a check we verify that (\ref{29}) is satisfied asymptotically.  In  terms of $n$ this identity becomes
\begin{equation}\label{730}
(1-\rho) \rho^N = \sum_{k+r=N}\pi(k,r)  = \sum_{m-n+r=N}L(n,r; m).
\end{equation}
We consider $N=m+\hat{N}$ with $m, N \rightarrow \infty$ but with $\hat{N}= O(1)$.  Then the corner range approximation can be used to approximate $\pi(k,r)$ in the sum in (\ref{730}).  We may have $\hat{N}<0$ or $\hat{N}\geq 0$.  If $\hat{N}<0, \; k+r= m+\hat{N}<m$ and the approximations in (\ref{742}) and (\ref{749}) yield
\begin{eqnarray*}
\sum_{r=n+\hat{N}}L(n,r; m) &=& \sum_{r=0}^{m+\hat{N}} L( r-\hat{N}, r; m) \sim \sum_{r=0}^\infty L(r-\hat{N}, r) \\
&=&  L(-\hat{N}, 0)+\sum_{r=1}^\infty L(r-\hat{N}, r)\\
&=& \rho^m \left[ (1-\rho)(\rho^{\hat{N}}-\rho)+ \sum_{r=1}^\infty (1-\rho)^2\rho^r \right] \\
&=& (1-\rho) \rho^N.
\end{eqnarray*}
so that (\ref{29}) holds for $\hat{N}<0$.
For $\hat{N}>0$ we have $N>m$ and then necessarily $r\geq 1$ so that
\begin{eqnarray*}
\sum_{r=n+\hat{N}}L(n,r; m) &=& \sum_{r=\hat{N}}^{m+\hat{N}} L( r-\hat{N}, r; m)\\
&\sim&  \sum_{r=\hat{N}}^\infty L(r-\hat{N}, r) = \rho^m (1-\rho)^2 \sum_{r=\hat{N}}^\infty \rho^r
\end{eqnarray*}
which again evaluates to $(1-\rho)\rho^{m+\hat{N}}$ and $m+\hat{N}=N$.  When $N=m$ ($\hat{N}=0$) either of the above calculations apply, since $\pi(m, 0)\sim L(0,0)$ can be computed by either of the formulas in (\ref{742}) and (\ref{749}).  This means that $\pi(k,0)$ can be computed by setting $r=0$ in the expression for $r\geq 1$, but only if $k=m$.  This observation was also used in the analysis in section \textbf{5}.

\subsection{Corner layer near $(x,y)=(0,0)$}
We consider near $(x,y)=(0,0)$ and go back to the original discrete variable $(k,r)$, with $k, r = O(1)$.  We also set
\begin{equation}
\pi(k,r) = \rho^m m^{-r/(1-\rho)} P(k,r; m), \; \; r>0, \label{781}
\end{equation}
and
\begin{equation}
\pi(k,0) = (1-\rho)\rho^k+\rho^m m^{-1/(1-\rho)}P(k,0; m), \; \; r=0. \label{753}
\end{equation}
We also require that this corner expansion match to the expansion in subsection \textbf{4.2}, in an intermediate limit where $k\rightarrow \infty$ but $x=k/m \rightarrow 0$.  By expanding (\ref{738}) and (\ref{744}) as $x\rightarrow 0$, we see that $x^{r/(1-\rho)}$ becomes $O(m^{-r/(1-\rho)})$ on the $k$-scale. We thus scaled $\pi(k,r)$ to be $O(\rho^m m^{-r/(1-\rho)})$ expecting that $P(k,r;m)$ will be $O(1)$.  Note also that $\pi(k,1)$ and $\pi(k, 0)-(1-\rho)\rho^k$ are scaled to be of the same order in $m$.
Then the balance equation (\ref{bl2}) becomes
\begin{eqnarray}\label{752}
(1+\rho) P(k,r; m) &=& \rho   P(k-1,r; m)+\frac{k+1}{k+r+1} \; P(k+1,r; m) \nonumber\\
 &+&\frac{r+1}{k+r+1}\; m^{-1/(1-\rho)} P(k,r+1 ; m),  \;  \;  k>0, \; \; r > 0, \nonumber\\
 \\
(1+\rho) P(k,0; m) &=& \rho P(k-1,0; m)+  P(k+1,0; m) \nonumber\\
&+& \frac{1}{k+1} \; P(k,1; m), \;  \;  k >0, \; r = 0. \label{754}
\end{eqnarray}
Expanding $P(k,r;m)$ as $P(k,r;m)=P(k,r)+o(1)$ we obtain from (\ref{752}) and (\ref{754}) the following equations for the leading term
\begin{eqnarray}\label{755}
(1+\rho)  P(k,r) &=& \rho P(k-1,r)+\frac{k+1}{k+r+1} \; P(k+1,r), \; \;  r > 0, \\
(1+\rho) P(k,0) &=& \rho P(k-1,0)+  P(k+1,0) +\frac{1}{k+1} \; P(k,1), \;  \;  k >0, \; r = 0.\nonumber\\ \label{747}
\end{eqnarray}
Here we define  $P(-1,r)=0$  so that (\ref{755}) holds for $k\geq0$.  However, we must now consider the corner condition (\ref{26}).
Since $\pi(0,0)=1-\rho$,  (\ref{753}) implies that $P(0,0)=0$ and thus (\ref{26}) asymptotically becomes
\begin{equation} \label{764}
0 = P(1,0)+P(0,1).
\end{equation}
This is also consistent with (\ref{747}) at $k=0$.  We shall first analyze (\ref{755}) for $r>0$, and then solve (\ref{747}) for $P(k,0)$.

For $r>0$ we use the generating function $\xi(z)= \xi(z;r)\equiv \sum_{k=0}^\infty P(k,r) z^k$.
Multiplying (\ref{755}) by $z^k$ and summing over  $k\geq 0$ gives
\begin{equation} \label{756}
\left[-\rho z^2+(1+\rho)z-1\right]\xi'(z) = \left[\rho (r+2)z-(1+\rho)(r+1)\right]\xi(z).
\end{equation}
The general solution to (\ref{756}) is
\begin{equation*}
\xi(z;r) = \xi_*(r)(1-z)^{-1-r/(1-\rho)}\left(\frac{\rho}{1-\rho z}\right)^{1-\rho r/(1-\rho)}
\end{equation*}
where $\xi_*(r)$ is a function of $r$.  In particular,
\begin{equation*}
\xi(0;r) = \xi_*(r)\rho^{1-\rho r/(1-\rho)} = P(0,r).
\end{equation*}
Inverting the generating function, we can write for $r>0$,
\begin{equation} \label{757}
P(k,r) = \xi_*(r) \frac{1}{2\pi i}\oint \frac{z^{-k-1}}{(1-z)^{1+r/(1-\rho)}}\left(\frac{\rho}{1-\rho z}\right)^{1-\rho r/(1-\rho)} dz,
\end{equation}
where the integral is a complex contour integral along a small loop around $z=0$.
We obtain $\xi_*(r)$ by matching. Letting $x\rightarrow 0$ ($x=k/m$) in (\ref{738}), and invoking the asymptotic matching condition shows that
\begin{equation}\label{745}
P(k,r) \sim (1-\rho)\left(\frac{\rho}{1-\rho}\right)^r k^{r/(1-\rho)}, \; \; k\rightarrow \infty.
\end{equation}
We now expand (\ref{757}) as $k\rightarrow \infty$ and verify that (\ref{745}) is satisfied.  This will also determine the function $\xi_*(r)$.
As $k\rightarrow \infty$ the behavior of the integral in (\ref{757}) is determined by the singularity closest to $z=0$, which occurs at $z=1$.  Thus, expanding the integrand in (\ref{757}) near $z=1$ we obtain
\begin{eqnarray}
P(k,r) &\sim& \xi_*(r)\left(\frac{\rho}{1-\rho}\right)^{1-\rho r/(1-\rho)}\frac{\Gamma\left(k+1+\frac{r}{1-\rho}\right)}{\Gamma\left(1+\frac{r}{1-\rho}\right)k\Gamma(k)} \nonumber\\
&=& \xi_*(r)\left(\frac{\rho}{1-\rho}\right)^{1-\rho
r/(1-\rho)}\frac{1}{k \; B\left(k,
1+\frac{r}{1-\rho}\right)}.\label{758}
\end{eqnarray}
Expanding the Beta function as $k\rightarrow \infty$ yields
\begin{equation*}
B\left(k, 1+\frac{r}{1-\rho}\right) \sim
\Gamma\left(1+\frac{r}{1-\rho}\right)k^{-1-r/(1-\rho)},
\end{equation*}
and thus (\ref{758}) simplifies to
\begin{equation*}
P(k,r)\sim \xi_*(r)\left(\frac{\rho}{1-\rho}\right)^{1-\rho
r/(1-\rho)}\frac{k^{r/(1-\rho)}}{\Gamma\left(1+\frac{r}{1-\rho}\right)}.
\end{equation*}
Therefore, (\ref{745}) is satisfied if
\begin{equation*}
\xi_*(r)=\frac{(1-\rho)^2}{\rho}\left(\frac{\rho}{1-\rho}\right)^{r/(1-\rho)}
\Gamma\left(1+\frac{r}{1-\rho}\right).
\end{equation*}

We show that (\ref{757}), when expanded for $r\rightarrow \infty$ asymptotically matches to the expansion in subsection \textbf{4.1} as $y\rightarrow 0$.
Letting $y\rightarrow 0$ in (\ref{760}), (\ref{761}), and (\ref{759}) we obtain
\begin{eqnarray*}
S_0(y)&\sim& (1-\rho)^{3/2}\sqrt{2\pi y}, \\
\phi(0,y) &=& (y+1)\log\rho
+\frac{y \log y}{1-\rho}-\frac{2y\log(1-\rho)}{1-\rho}-\frac{y}{1-\rho}+O(y^2),\\
\phi_y(0,y) &=& \frac{1}{1-\rho}\log y+\log\rho -\frac{2}{1-\rho}\log(1-\rho)+o(1).
\end{eqnarray*}
Thus as $y \rightarrow 0$ ($y=r/m$) the asymptotic behavior of (\ref{732}) is
\begin{equation}
\pi(k,r) \sim \sqrt{2\pi r} (1-\rho)^{3/2-2r/(1-\rho)}(1+\rho)^k\rho^{r+m}\frac{r^k}{k!}\left(\frac{r}{m}\right)^{r/(1-\rho)} e^{-r/(1-\rho)} \label{762}.
\end{equation}
This is an approximation to $\pi(k,r)$ that applies in the matching region where $k=O(1), \; \; r\rightarrow \infty, \; \;  r/m\rightarrow 0.$

Next we expand the corner approximation (\ref{781}), which is given by (\ref{763}),
as $r\rightarrow \infty$.  We expand the integrand in (\ref{763}) around $z=0$ by setting $z=u/r$ and evaluate the integral, using
\begin{equation*}
\frac{r^k}{2\pi i}\oint u^{-k-1} e^{(1+\rho)u} du = \frac{r^k}{k!}(1+\rho)^k.
\end{equation*}
Using Stirling's formula we can expand the Gamma function in (\ref{763}) as
\begin{equation}\label{768}
\Gamma\left(1+\frac{r}{1-\rho}\right) \sim \sqrt{\frac{2\pi r}{1-\rho}}e^{-r/(1-\rho)}\left(\frac{r}{1-\rho}\right)^{r/(1-\rho)}, \; \; r\rightarrow \infty.
\end{equation}
Therefore, (\ref{763}) as $r\rightarrow \infty$ agrees with (\ref{762}).

Next we verify that the corner approximation (\ref{763}) matches to the ray solution $K(x, y)\exp[m \phi(x, y)]$, in an intermediate limit where $k, r \rightarrow \infty$ but $x=k/m, \; y=r/m \rightarrow 0$.  We expand (\ref{763}) for $k$ and $r$ simultaneously large, writing the integrand as
\begin{eqnarray*}
\frac{\exp[H(z)]}{z(1-z)(1-\rho z)}&\equiv& \frac{1}{z(1-z)(1-\rho z)}\nonumber\\
&\times& \exp\left[-k \log z+\frac{\rho r}{1-\rho}\log(1-\rho z)-\frac{r}{1-\rho}\log(1-z)\right]
\end{eqnarray*}
and using the saddle point method.  The saddle point(s) satisfy $H'(z)=0$ so that
\begin{equation} \label{772}
\frac{k}{r} = \frac{z}{1-\rho}\left(-\frac{\rho^2}{1-\rho z}+\frac{1}{1-z}\right).
\end{equation}
(\ref{772}) defines the saddle as a function of $k/r$ and the saddle $z_0=z_0(k/r)$ is given explicitly by
\begin{equation}\label{786}
z_0 = \frac{1}{2 \rho}\left[1+\rho+\sqrt{1-2\rho\left(\frac{k-r}{k+r}\right)+\rho^2}\right].
\end{equation}
We note that $z_0 \rightarrow 0$ as $k/r \rightarrow 0$ and $z_0 \rightarrow 1$ as $k/r \rightarrow \infty$.
Evaluating the integral along the steepest descent path through the saddle point we obtain
\begin{eqnarray}\label{750}
\pi(k,r) &\sim& \rho^{m+r}m^{-r/(1-\rho)}(1-\rho)^{2-r/(1-\rho)}\; \Gamma\left(1+\frac{r}{1-\rho}\right) \nonumber\\
&\times& \frac{1}{\sqrt{2\pi  F''(z_0)}\;(1-z_0)(1-\rho z_0)\; z_0^{k+1}}\left[\frac{(1-\rho z_0)^{\rho}}{1-z_0}\right]^{r/(1-\rho)} \nonumber\\
&\sim& (1-\rho)^{3/2} \sqrt{\frac{r}{F''(z_0)}} \; \frac{1}{z_0(1-z_0)(1-\rho z_0)} \nonumber\\
&\times& \rho^{m+r}z_0^{-k}\left[\frac{1}{e}\frac{r}{m}\frac{(1-\rho z_0)^{\rho}}{(1-\rho)^2(1-z_0)}\right]^{r/(1-\rho)}.
\end{eqnarray}
In the ray solution we let $x \rightarrow 0$ and $y\rightarrow 0$, which corresponds to $s \rightarrow 0$ and $\; t\rightarrow \infty$, with $x/y$ fixed (corresponding to $s e^{(1-\rho)t} =O(1)$).  We find that in this limit
\begin{eqnarray*}
\frac{1-\rho-s e^{(1-\rho)t}}{1-\rho-\rho s e^{(1-\rho)t}} &\sim& z_0,\\
\phi_x &\sim& -\log z_0, \\
\phi_y &\sim& \log \rho + \frac{1}{1-\rho}[\log y -2\log(1-\rho)+\rho \log (1-\rho z_0)-\log(1-z_0)],\\
\log(1+s) &\sim&  s \sim  \frac{y}{1-\rho}.
\end{eqnarray*}
Therefore, as $(x, y)\rightarrow (0,0)$, in terms of $z_0=z_0(x/y)=z_0(k/r)$,
$$\exp[m\phi(x,y)] \sim \rho^{m+r}z_0^{-k}\left[\frac{1}{e}\frac{r}{m}\frac{(1-\rho z_0)^{\rho}}{(1-\rho)^2(1-z_0)}\right]^{r/(1-\rho)}.$$
This agrees precisely with the exponentially varying terms in (\ref{750}).  From (\ref{772}) and (\ref{786}) we also obtain
\begin{eqnarray*}
\frac{1}{(1-z_0)(1-\rho z_0)} &=& \frac{k+r}{r},\\
 z_0\sqrt{F''(z_0)} &=& \frac{(1+\rho)}{\sqrt{2\rho}}\frac{(k+r)}{\sqrt{r}}\sqrt{\nabla-\nabla^2}
\end{eqnarray*}
where
\begin{equation}\label{770}
\nabla \equiv \sqrt{1-\frac{4\rho k}{(1+\rho)^2(k+r)}} = \sqrt{1-\frac{4\rho}{(1+\rho)^2}\frac{x}{(x+y)}}.
\end{equation}
It follows that
\begin{equation}\label{787}
\sqrt{\frac{r}{F''(z_0)}} \; \frac{(1-\rho)^{3/2}}{z_0(1-z_0)(1-\rho z_0)} =   \frac{(1-\rho)^{3/2}}{(1+\rho)}\sqrt{\frac{2\rho}{\nabla-\nabla^2}}.
\end{equation}
By asymptotic matching, the above should agree with the expansion of $K(x, y)$ as $x, y\rightarrow 0$.
As $s\rightarrow 0$ and $t\rightarrow \infty$ with $s e^{(1-\rho)t} =O(1)$, from (\ref{769}) we obtain
\begin{equation}\label{771}
K(x,y) \sim  \frac{(1-\rho)[1-\rho-\rho s e^{(1-\rho)t}]}{[1-\rho-s e^{(1-\rho)t}]^{1/2}[1-\rho+\rho s e^{(1-\rho)t}]^{1/2}}.
\end{equation}
Since $x(s,t) \sim [1-\rho-s e^{(1-\rho)t}][(1-\rho) e^{-(1-\rho)t}-s\rho^2]/(1-\rho)^2$
and $y(s,t) \sim s(1-\rho)$ in this limit, (\ref{771}) agrees with (\ref{787}). This completes the matching verifications.

It remains to compute $P(k,0)$, by solving (\ref{747}) and (\ref{764}).  We use the generating function
$$\xi_0(z) \equiv \sum_{k=0}^\infty P(k,0) z^k,$$
multiply (\ref{747}) by $z^k$, and sum over $k\geq 0$, to obtain
\begin{equation*}
\xi_0(z) = \frac{1}{(1+\rho)z-\rho z^2-1}\sum_{k=0}^\infty \frac{z^{k+1}}{k+1}P(k,1),
\end{equation*}
and thus
\begin{eqnarray}
P(k,0) &=& \frac{1}{2\pi i}\oint z^{-k-1}\xi_0(z) dz \nonumber\\
&=& -\frac{1}{2\pi i} \oint \frac{z^{-k-1}}{(1- z)(1-\rho z)}\left[\sum_{j=0}^\infty \frac{z^{j+1}}{j+1}P(j,1)\right] dz. \label{765}
\end{eqnarray}
Evaluating the contour integral leads to
\begin{equation}
P(k,0) = -\sum_{j=0}^{k-1} \frac{1-\rho^{k-j}}{1-\rho}\;\frac{P(j,1)}{j+1}.
\end{equation}

Finally we verify the matching between (\ref{765}) and (\ref{744}).  To evaluate the integral in (\ref{765}) as $k\rightarrow \infty$,  we expand the integrand around $z=0$ and obtain
\begin{eqnarray}
P(k,0) &\sim& -\frac{1}{1-\rho}\left[\frac{P(k-1,1)}{k}+\frac{P(k-2,1)}{k-1}+\frac{P(k-3,1)}{k-2}+\cdot\cdot\cdot+P(0,1)\right] \nonumber\\
&\sim& -\frac{\rho}{1-\rho}\sum_{n=1}^k\frac{(n-1)^{1/(1-\rho)}}{n}\nonumber\\
&\sim& -\frac{\rho}{1-\rho} \int_1^k x^{1/(1-\rho)-1} dx \sim -\rho k^{1/(1-\rho)}. \label{766}
\end{eqnarray}
Therefore, from (\ref{753}) and (\ref{766}) we obtain
\begin{equation}\label{767}
\pi(k,0) \sim (1-\rho)\rho^k-  \rho^{m+1} \left(\frac{k}{m}\right)^{1/(1-\rho)}.
\end{equation}
This applies for $k\rightarrow \infty$ but with $k/m\rightarrow 0$.
If we let $x\rightarrow 0$ in (\ref{744}), we obtain precisely (\ref{767}), which verifies the matching. Thus, for $r=0$ the approximation to $\pi(k,r)$ must be computed by using (\ref{757}) with $r=1$ in (\ref{765}) and then (\ref{753}) with $P(k, 0; m)\sim P(k,0)$.

\subsection{Marginal distributions}

We define the marginal distributions by (\ref{m}) and (\ref{nr}).
First, we discuss ${\cal M}(k)$ for various ranges of $k$. Summing  (\ref{bl}) over $r$, we obtain
\begin{multline}\label{773}
\sum_{r=0}^\infty(1I_{[k+r>0]}+\rho)(k+r+1)\; \pi(k,r) = \rho
I_{[k\geq1]} \;\sum_{r=0}^\infty (k+r+1)\pi (k-1,
r)\\
+(k+1)\sum_{r=0}^\infty \pi (k+1, r) +\sum_{r=0}^\infty(r+1)\pi(k,
r+1).
\end{multline}
If we define ${\cal M}_1(k) = \sum_{r=1}^\infty r\;\pi(k,r)$, then (\ref{773})
can be written as
\begin{multline}\label{774}
\rho(k+1)[{\cal M}(k)-{\cal M}(k-1)]-(k+1)[{\cal M}(k+1)-{\cal M}(k)]\\
= \pi(0,0)I_{[k=0]}-\rho[{\cal M}_1(k)-{\cal M}_1(k-1)]
\end{multline}
where ${\cal M}_1(-1) = 0$ and ${\cal M}(-1)=0$. Setting ${\cal M}(k) =
(1-\rho)\rho^k+\tilde{M}(k)$, from (\ref{774}), we obtain
\begin{equation}\label{775}
[\tilde{M}(k+1)-\tilde{M}(k)] - \rho[\tilde{M}(k)- \tilde{M}(k-1)] = \frac{\rho}{(k+1)}[
{\cal M}_1(k)-{\cal M}_1(k-1)]
\end{equation} where $\tilde{M}(-1) =0$.

On the $x$ scale (corresponding to $k \rightarrow \infty, \; 0<x<1$), if we set $\tilde{M}(k) = \rho^m \bar{M}(x; \delta)$ and ${\cal M}_1(k) = \rho^m\bar{M}_1(x; \delta)$, (\ref{775}) becomes
\begin{multline}\label{776}
(x+\delta)[\bar{M}(x+\delta; \delta)-\bar{M}(x; \delta)]-\rho (x+\delta)[\bar{M}(x; \delta)- \bar{M}(x-\delta; \delta)]\\
=\delta \rho [\bar{M}_1(x; \delta)-\bar{M}_1(x-\delta; \delta)].
\end{multline}
Expanding $\bar{M}(x; \delta)$ and $\bar{M}_1(x; \delta)$ as
$$\bar{M}(x; \delta) = \delta \bar{M}(x)+O(\delta^2), \; \; \bar{M}_1(x; \delta) = \bar{M}_1(x)+O(\delta),$$
from (\ref{776}), we obtain to leading order
\begin{equation}\label{797}
\bar{M}'(x)= \frac{\rho}{1-\rho}\frac{\bar{M}_1'(x)}{x}.
\end{equation}
For $0<x<1$ $\; \pi(k,r)$ is maximal in the range $r=O(1)$ for $r\geq 1$.  Thus we use (\ref{738}) to evaluate the sum in (\ref{m}), to obtain
\begin{eqnarray}\label{777}
{\cal M}_1(k)&=&\rho^m \bar{M}_1(x ; \delta) \nonumber\\
&\sim&\rho^m \frac{(1-\rho)^2}{[1-\rho+\rho x^{1/(1-\rho)}]} \sum_{r=1}^\infty r \left\{\frac{\rho x^{1/(1-\rho)}}{[1-\rho+\rho x^{1/(1-\rho)}]}\right\}^r \nonumber\\
&=& \rho^{m+1} x^{1/(1-\rho)}.
\end{eqnarray}
This gives the mean number of occupied secondary spaces, if the number of occupied primary spaces is $k=mx$.

To leading order we have ${\cal M}(k) \sim (1-\rho)\rho^k \sim \pi(k, 0)$.
To obtain the correction (second) term for ${\cal M}(k)$, we need to solve (\ref{797}) and then use (\ref{777}), which yields
\begin{equation}
\bar{M}(x)=\frac{\rho}{1-\rho}\int_0^x \frac{1}{\varsigma}\bar{M}_1'(\varsigma)d\varsigma
= \frac{\rho}{1-\rho} \; x^{\rho/(1-\rho)}.
\end{equation}
Therefore, we obtain (\ref{799}).
Using (\ref{799}) we note that
\begin{equation}
\sum_{k=0}^m {\cal M}(k) \sim 1-\rho^{m+1}+\frac{\rho^{m+1}}{1-\rho} \int_0^1 x^{1/(1-\rho)-1}dx = 1.
\end{equation}
Thus, the marginal is properly normalized up to order $O(\rho^m)$.  Expression (\ref{799}) breaks down as $x\rightarrow 0$, since then  we must use the corner expansion for $\pi(k, r)$ to evaluate (\ref{m}).

For the $k$ scale ($x = \delta k, \; k=O(1)$) we must use the approximation to $\pi(k, r)$ valid for $k, r = O(1)$.  Thus from (\ref{781}) and (\ref{753}) we obtain
\begin{eqnarray}
\tilde{M}(k) &\sim& \rho^m m^{-1/(1-\rho)}[P(k,0)+P(k,1)] \label{782}, \\
{\cal M}_1(k) &\sim& \rho^m m^{-1/(1-\rho)} P(k,1) \label{783}
\end{eqnarray}
and we can easily check that (\ref{782}) and (\ref{783}) satisfy (\ref{775}).  Then $P(k, 1)$ and $P(k, 0)$ are given by (\ref{757}) and (\ref{765}).

We verify the matching between the $k$ and $x$ scales. For $k\rightarrow \infty$, we use (\ref{783}) and write (\ref{775}) asymptotically as
\begin{equation}\label{798}
(1-\rho)\tilde{M}'(k) \sim \rho^{m+1} m^{-1/(1-\rho)}\frac{1}{(k+1)}[ P(k,1) -P(k-1,1)].
\end{equation}
Also for $k\rightarrow \infty$, we use (\ref{745}) in (\ref{798}) and obtain
\begin{equation}
\tilde{M}'(k)\sim \frac{\rho^{m+2}}{(1-\rho)^2} m^{-1/(1-\rho)} k^{1/(1-\rho)-2},
\end{equation}
and thus
\begin{equation}
\tilde{M}(k) \sim \frac{\rho^{m+1}}{(1-\rho)} \; m^{-1/(1-\rho)} k^{1/(1-\rho)-1}.
\end{equation}
This matches to $\tilde{M}(k)$ on the $x$ scale, in view of (\ref{799}).

Near the boundary of $x=1$, to obtain ${\cal M}(k)$ we use the corner $(x,y)=(1,0)$ solution in subsection \textbf{4.3}.  From (\ref{742}) (or \ref{748})) and (\ref{749}) (or (\ref{709})), we obtain
\begin{equation*}
{\cal M}(k) \sim \sum_{r=1}^\infty (1-\rho)^2 \rho^{m+r} + (1-\rho)(\rho^{-n}-\rho)\rho^m = (1-\rho)\rho^{m-n} = (1-\rho)\rho^k.
\end{equation*}
Thus the leading term for the marginal ${\cal M}(k)$ is the geometric distribution $(1-\rho)\rho^k$, for all ranges of $0\leq k\leq m$.  The correction term is roughly $O(\rho^m)$ and contains the effects of the (rarely used) secondary servers.

Next we discuss ${\cal N}(r)$, the marginal distribution of the number of occupied secondary spaces, for various ranges of $y=\delta r$.  On the $y$ scale ($r = O(m)$), we use the fact that $\pi(k,r)$ is maximal in the range $k=m-O(1)$.  By summing (\ref{727}) over $k\leq m$ we obtain the distribution of $N_2$ as
\begin{equation*}
{\cal N}(r) = \sum_{n=0}^m \pi(m-n,r) \sim (1-\rho)\left(\frac{y+1-\rho}{y+1}\right)\rho^{r+m}\sum_{n=0}^m (y+1)^{-n},
\end{equation*}
which can be written as (\ref{7b}).
(\ref{7b}) applies only for $y>0$, and has a singularity as $y\rightarrow 0$.

On the $r$ scale ($r=O(1)$), we integrate (\ref{738}) over $0\leq x\leq 1$ and obtain (\ref{7a}) from
\begin{equation}
{\cal N}(r) \sim  m (1-\rho)^2 \rho^{m+r} \int_0^1 \frac{x^{r/(1-\rho)}}{[1-\rho+\rho x^{1/(1-\rho)}]^{r+1}} dx.
\end{equation}
Here we approximated the sum in (\ref{nr}) by an integral and changed variables using $u=x^{1/(1-\rho)}$.  To check the matching between the $r$ scale and the $y$ scale we let $r\rightarrow \infty$ in (\ref{7a}). We expand the integrand in (\ref{7a}) around $u=1$ since it is an increasing function as $r\rightarrow \infty$.  If we set $u=1-w/r$, (\ref{7a}) becomes
\begin{eqnarray}
{\cal N}(r) &\sim& m (1-\rho)^3\rho^{m+r} \int_0^r \frac{(1-w/r)^{r-\rho}}{(1-\rho w/r)^{r+1}}\frac{dw}{r}, \nonumber\\
&\sim& (1-\rho)^3\rho^{m+r} \; \frac{m}{r}\int_0^\infty e^{-(1-\rho)w} dw = (1-\rho)^2\rho^{m+r}\;\frac{m}{r}. \label{7c}
\end{eqnarray}
We can easily check that letting $y\rightarrow 0$ in (\ref{7b}) and using $y=r/m$, we also obtain (\ref{7c}).

We need a different expression for ${\cal N}(r)$ when $r=0$. From (\ref{744}), again summing over $0\leq k \leq m$ ( and integrating the correction term over $0\leq x \leq 1$) we obtain
\begin{eqnarray}
{\cal N}(0) &=& \sum_{k=0}^m\pi(k,0) \nonumber\\
&\sim& \sum_{k=0}^m (1-\rho)\rho^k-  m(1-\rho)\rho^{m+1}\int_0^1 \frac{x^{1/(1-\rho)}}{1-\rho+\rho x^{1/(1-\rho)}} dx \nonumber\\
&=& 1- m(1-\rho)^2\rho^{m+1}\int_0^1 \frac{u^{1-\rho}}{1-\rho+\rho u} du + O(\rho^{m}). \label{7d}
\end{eqnarray}
The normalization condition is satisfied since
\begin{eqnarray*}
\sum_{r=1}^\infty {\cal N}(r) &\sim& m (1-\rho)^3\rho^m \int_0^1 \frac{u^{-\rho}}{1-\rho+\rho u}\left[\sum_{r=1}^\infty\left(\frac{\rho u}{1-\rho+\rho u}\right)^r\right] du \\
&=& m (1-\rho)^2\rho^{m+1} \int_0^1 \frac{u^{1-\rho}}{1-\rho+\rho u} du \sim 1-{\cal N}(0).
\end{eqnarray*}
We have thus estimated the difference $1-{\cal N}(0)$, which gives the exponentially small probability of having to use the secondary spaces.  Note also that the mean number of occupied secondary spaces is, from (\ref{7a}),
\begin{eqnarray*}
\sum_{r=1}^\infty r{\cal N}(r) &\sim& m (1-\rho)^3\rho^m \int_0^1 \left[\sum_{r=1}^\infty \frac{r (\rho u)^r}{(1-\rho+\rho u)^r}\right] \frac{ u^{-\rho}}{(1-\rho+\rho u)}du \nonumber\\
&=& m(1-\rho)\rho^{m+1}\int_0^1 u^{1-\rho}du = \left(\frac{1-\rho}{2-\rho}\right) m \rho^{m+1}.
\end{eqnarray*}

\section{Asymptotic expansions for $m \rightarrow \infty$ with $\rho \uparrow 1$}

In this section, we study the heavy traffic case, in which $\rho \uparrow 1$, with $m \rightarrow \infty$.  We introduce the parameter $a$, with
\begin{equation}
\rho = 1-\frac{a}{m} = 1-a \delta,   \;   \; a= O(1).
\end{equation}
We shall again analyze (\ref{bl2}) for different ranges of $(k,r)$ (or $(x, y)$).  For some of the ranges the analysis will closely parallel that of section \textbf{4}, which had $\rho<1$.  However, when $y\approx 0$ the analysis is much different, and we will show that the scalings $r=O(\sqrt{m})$ ($y= O(\sqrt{\delta})$)  and $r=O(1)$ ($y=O(\delta)$) will lead to very different types of asymptotics.

We first consider $0<x \leq1$ and $y>0$,  and set
\begin{equation*}
\pi(k,r) = {\cal P}(x, y) ={\cal K}(x, y; \delta)\exp\left[\frac{1}{\delta} \Psi(x, y)\right].
\end{equation*}
Then  (\ref{bl2}) becomes
\begin{multline}
(2-a \delta)(x+y+\delta) \; {\cal K}(x, y ;  \delta) \exp\left[\frac{1}{\delta} \Psi(x,y)\right] \\
=(1-a \delta) (x+y+\delta)\; {\cal K}(x-\delta, y ; \delta) \exp\left[\frac{1}{\delta} \Psi(x-\delta, y)\right] \\
+(x+\delta)\; {\cal K}(x+\delta, y ; \delta) \exp\left[\frac{1}{\delta} \Psi(x+\delta,y)\right] \\
+(y+\delta)\; {\cal K}(x, y+\delta ; \delta) \exp\left[\frac{1}{\delta} \Psi(x,y+\delta)\right].  \label{8}
\end{multline}
The boundary condition (\ref{ab}) at $x=1$ ($k=m$) becomes
\begin{equation}\label{80}
(1-a\delta){\cal K}(1, y-\delta ; \delta)\exp\left[\frac{1}{\delta}\Psi(1, y-\delta)\right]=\frac{1+\delta}{1+y+\delta} \;{\cal K}(1+\delta, y ; \delta)\exp\left[\frac{1}{\delta}\Psi(1+\delta,y)\right].
\end{equation}
We assume that ${\cal K}(x, y ; \delta)$ has an asymptotic expansion in powers of $\delta$, with
\begin{equation}\label{819}
{\cal K}(x, y ; \delta) = \delta{\cal K}(x, y)+\delta^2 {\cal K}^{(1)}(x, y) + O(\delta^3).
\end{equation}
To compute $\Psi(x, y)$ we divide (\ref{8}) by $\delta {\cal K}(x, y)\exp[\Psi(x, y)/\delta]$ and let $\delta \rightarrow 0$. Then we similarly divide (\ref{80}) by
$\delta {\cal K}(1,y)\exp[\Psi(1, y)/\delta]$ and let $\delta \rightarrow 0$.  We thus obtain the following equation for $\Psi$, for $0<x<1$ and $y>0$,
\begin{equation}\label{801}
(x+y)(2-e^{-\Psi_x})-xe^{\Psi_x}-ye^{\Psi_y} = 0,
\end{equation}
and at $x=1$  the boundary condition
\begin{equation}\label{802}
(1+y)e^{-\Psi_y(1,y)} = e^{\Psi_x(1,y)}.
\end{equation}
We solve (\ref{801}) and (\ref{802}) for $\Psi(x,y)$ by using the method of
characteristics, as we did in section \textbf{4} again writing (\ref{801}) as ${\cal F}(x, y, \Psi, \Psi_x, \Psi_y)=0$ where
\begin{equation}
{\cal F} \equiv (x+y)(2-e^{-\Psi_x})-xe^{\Psi_x}-ye^{\Psi_y}. \label{803}
\end{equation}
The characteristic equations for this nonlinear PDE are \cite{CH}
\begin{eqnarray}
\dot{x}&=& \frac{dx}{dt} = \frac{\partial {\cal F}}{\partial\Psi_x}
=(x+y)e^{-\Psi_x}-xe^{\Psi_x}, \label{804}\\
\dot{y} &=& \frac{dy}{dt} = \frac{\partial {\cal F}}{\partial\Psi_y} = -ye^{\Psi_y}, \label{805}\\
\dot{\Psi}&=& \frac{d\Psi}{dt} = \Psi_x \dot{x}+\Psi_y \dot{y}, \label{806}\\
\dot{\Psi_x} &=& -\frac{\partial {\cal F}}{\partial x} =
-2+ e^{-\Psi_x}+e^{\Psi_x}, \label{807}\\
\dot{\Psi_y} &=& -\frac{\partial {\cal F}}{\partial y} =
-2+ e^{-\Psi_x}+e^{\Psi_y}. \label{808}
\end{eqnarray}
To solve this system we again use rays starting from $(x, y)=(1, s)$ at $t=0$.
We set $\Psi_x \equiv \log \hat{\Upsilon}(s,t)$ and, from (\ref{803}) - (\ref{805}) and (\ref{807}), obtain
\begin{eqnarray}
\frac{\dot{x}+\dot{y}}{x+y} &=& 2\left(\frac{1}{\hat{\Upsilon}}-1\right), \label{809} \\
\frac{\partial \hat{\Upsilon}}{\partial t} &=& (\hat{\Upsilon}-1)^2. \label{81}
\end{eqnarray}
The solutions to (\ref{81}) and (\ref{809}) are
\begin{eqnarray}
\hat{\Upsilon}&=& 1-\frac{1}{t+C_0(s)}, \label{811}\\
x+y &=& C_1(s)[t+C_0(s)-1]^2. \label{810}
\end{eqnarray}
Using (\ref{811}) and (\ref{810}) in (\ref{804}) we obtain
\begin{equation}
\frac{\partial x}{\partial t}+\left[1-\frac{1}{t+C_0(s)}\right] x = C_1(s)[t+C_0(s)][ t+C_0(s)-1]. \label{812}
\end{equation}
The general solution to (\ref{812}) is
\begin{equation}\label{813}
x(s,t) = [t+ C_0(s)]\{C_1(s)[t+C_0(s)-2] + C_2(s) e^{-t}\},
\end{equation}
Subtracting (\ref{813}) from (\ref{810}) we obtain
\begin{equation}
y(s, t) = C_1(s)-C_2(s)[t+C_0(s)]e^{-t}. \end{equation}
Applying the initial conditions ($x(s,0)=1, \;  y(s,0)=s$) we find that
\begin{eqnarray}
C_1(s) &=& \frac{s+1}{[C_0(s)-1]^2}, \nonumber\\
C_2(s) &=& \frac{1}{C_0(s)}-\frac{C_0(s)-2}{[C_0(s)-1]^2}(s+1). \label{814}
\end{eqnarray}

To determine $C_0(s)$ we set $\Psi_x=\alpha(s)$ and $\Psi_y=\beta(s)$ at $t=0$, and from (\ref{801}) and (\ref{802}) obtain
\begin{eqnarray}
2(1+s)&=&\rho(1+s)e^{-\alpha}+e^{\alpha}+se^{\beta}, \label{815}\\ (1+s)e^{-\alpha} &=& e^{\beta}. \label{816}
\end{eqnarray}
Solving (\ref{815}) and (\ref{816}) leads to
\begin{equation*}
\alpha(s) = \log(s+1), \; \; \beta(s)=0.
\end{equation*}
Therefore, $\Psi_x = \alpha(s)=\log(s+1)=\log \hat{\Upsilon} = \log[1-1/C_0(s)]$ at $t=0$.  This implies that $C_0(s) = -1/s$, which then gives $C_1(s) = C_2(s) =s^2/(s+1)$.  Thus we obtain $x(s,t)$ and $y(s, t)$ as (\ref{817}) and (\ref{818}).
From (\ref{811}) and (\ref{805}) we also obtain (\ref{Psix}) and (\ref{Psiy}).
In the same manner as in section \textbf{4}, we solve for $\Psi$ as
\begin{eqnarray*}
\Psi &=& \Psi_x x+\Psi_y y-\int (\dot{\Psi_x}x+\dot{\Psi_y}y) \; dt, \nonumber\\
&=& x\Psi_x+y\Psi_y+f(s).
\end{eqnarray*}
Here we used $\dot{\Psi_x}x+\dot{\Psi_y}y= 0$. At $t=0$, $\Psi(1,s) = \log(1+s)+f(s)$ and
$$\Psi_y(1, s) = \frac{\partial\Psi(1,s)}{\partial s} = \frac{1}{s+1}+f'(s) = 0.$$
This implies that $f(s) = -\log(s+1)+constant$.  We will later show that $\Psi$ must vanish at $x=1$, $\; y=0$ (i.e., $t=0$, $\; s=0$)  so that this $constant=0$. Therefore, we have (\ref{Psi}).

We plot the rays using (\ref{817}) and (\ref{818})  in \textbf{Fig. 4} and \textbf{Fig. 5}.  In \textbf{Fig. 4} we plot the rays for $s>0$ and $t$ in the range  $0 \leq t \leq t_{max}\equiv 1/s$. At $t=0$ we have $x=1$ and at $t=t_{max}(s)$ the rays reach $x=0$.  We note that $\Psi_x$ develops a singularity at $t=t_{max}(s)$.  We continue the rays for $t>t_{max}$ in \textbf{Fig. 5}.  There is a cusp when $\dot{x} =0$, which happens at $t= t_c(s) \equiv 1+1/s > t_{max}(s)$. At $t=t_c$, $x+y=0$ and thus the locus of the cusps is on the line $y=-x$.  By comparing \textbf{Fig. 2} and \textbf{Fig. 4} we see that the main difference in the rays occurs in their behavior near $y=0$, and indeed the asymptotics will be much different for $y\approx 0$ and $\rho \uparrow 1$, than what we found for $y\approx 0$ and $\rho<1$ in section \textbf{4}.  Near $y=0$ the ray expression breaks down and we will analyze these cases in subsections \textbf{5.2} - \textbf{5.5}.

\begin{center}
\includegraphics[angle=0, width=0.6\textwidth]{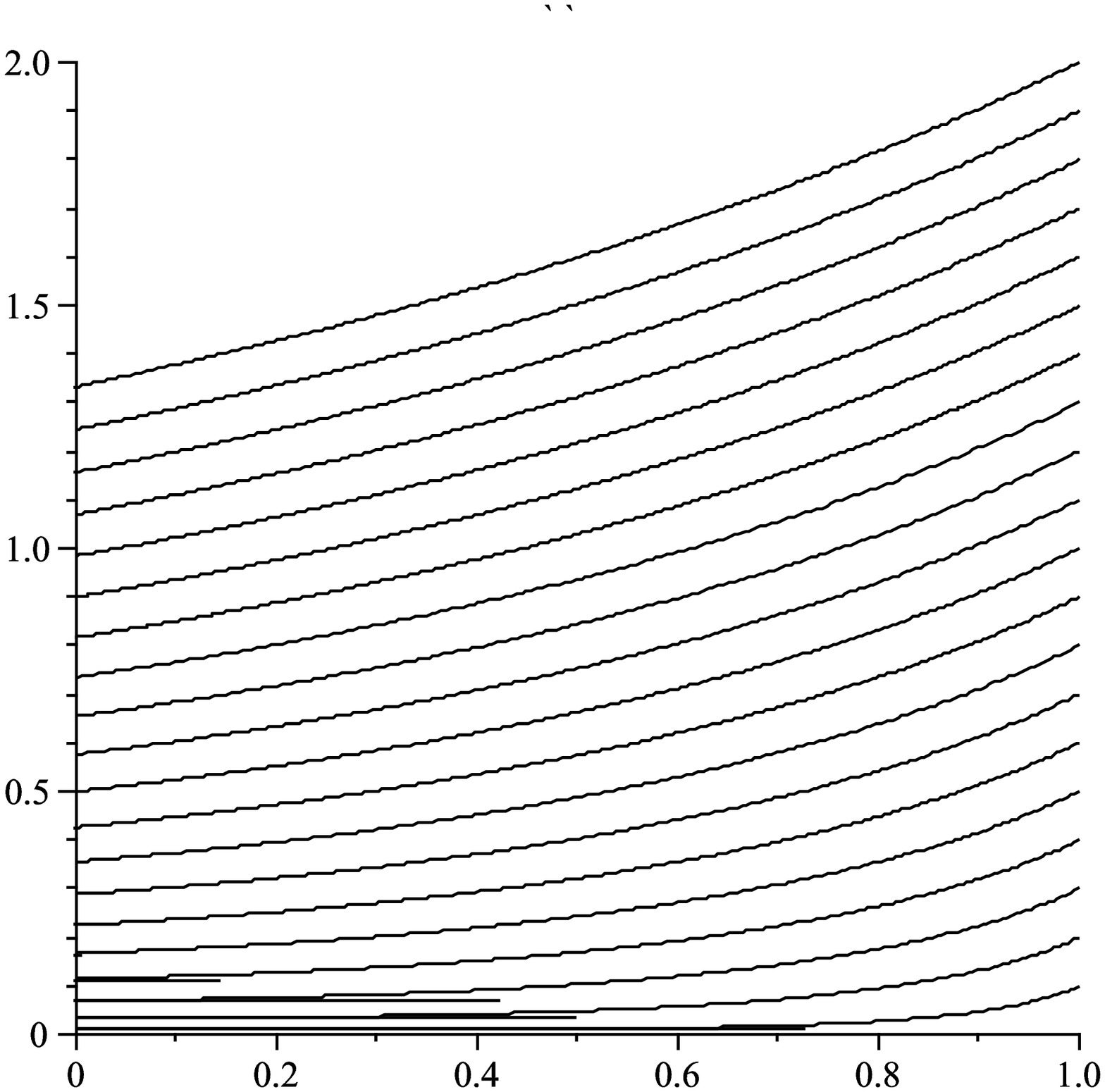}\\
\textbf{Fig. 4}   A sketch of the rays from $x=1$.
\end{center}

\begin{center}
\includegraphics[angle=0, width=0.6\textwidth]{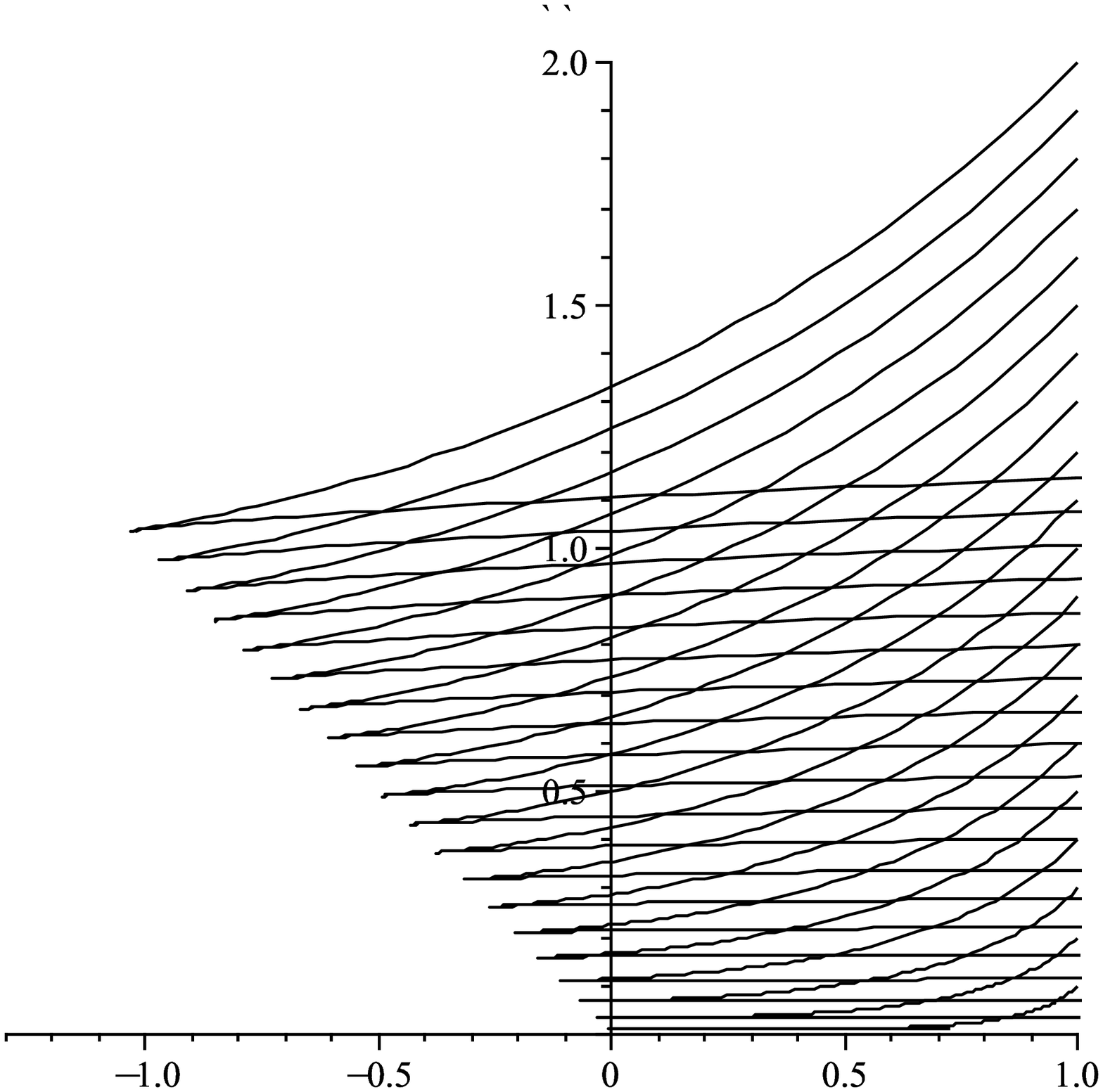}\\
\textbf{Fig. 5}   A sketch of the extended rays from $x=1$.
\end{center}

We next calculate ${\cal K}(x,y)$ in (\ref{819}).  In the same manner as in section \textbf{4} we obtain the following 'transport' equation:
\begin{multline}
\{ \left[\left(\frac{\Psi_{xx}}{2}-a\right)(x+y)+1\right] e^{-\Psi_x} + \left(\frac{x}{2} \Psi_{xx} + 1\right)e^{\Psi_x} +\left(\frac{y}{2} \Psi_{yy}+1\right)e^{\Psi_y} \\
-2+a(x+y)\; \} {\cal K}
=\left[(x+y) e^{-\Psi_x}-x e^{\Psi_x}\right]{\cal K}_x-y e^{\Psi_y}{\cal K}_y.\label{820}
\end{multline}
The right hand side is $\dot{x}{\cal K}_x+\dot{y}{\cal K}_y = \dot{{\cal K}}$ and thus, from (\ref{804}), (\ref{805}), (\ref{807}), (\ref{808}), and (\ref{716}) - (\ref{718}), we can write (\ref{820}) as
\begin{equation}\label{821}
\frac{\dot{{\cal K}}}{{\cal K}} = -\frac{1}{2}\frac{\dot{\jmath}}{\jmath}+\frac{1}{2}(\dot{\Psi_x}+\dot{\Psi_y})+a(x+y)(1-e^{-\Psi_x})+\frac{1}{2} e^{-\Psi_x}.
\end{equation}
Integrating (\ref{821}) using (\ref{Psix}) for $\Psi_x$ we obtain
\begin{equation} \label{822}
{\cal K} = {\cal K}_0(s)|\jmath|^{-1/2} \exp\left(\frac{\Psi_x+\Psi_y}{2}\right)|s-st+1|^{1/2}\exp\left[ast-\frac{as^2t^2}{2(s+1)}+\frac{t}{2}\right].
\end{equation}
where ${\cal K}_0(s)$ is a function of $s$.
Using (\ref{817}) and (\ref{818}), we obtain $\jmath$ as
\begin{equation}
\jmath = -\frac{s(s-st+1)}{(s+1)^3}\left\{2s(s+2)-[s(s+2)t^2-2t+s(s+2)-1]e^{-t}\right\}.
\end{equation}
Therefore,
\begin{eqnarray}
{\cal K} &=& \frac{{\cal K}_0(s)(s+1)^{3/2}(s-st+1)e^{t/2}}{\sqrt{s(1-st)(se^t-st+1)}\sqrt{2s(s+2)-[s(s+2)t^2-2t+s(s+2)-1]e^{-t}}} \nonumber\\
&\times& \exp\left[ast -\frac{as^2t^2}{2(s+1)}\right].  \label{823}
\end{eqnarray}

We determine ${\cal K}_0(s)$  by matching ${\cal P}(x, y)$ to the ray solution in section \textbf{4}.  Setting $\rho=1-a\delta$ in (\ref{phi}) and expanding the result for $\delta \rightarrow 0$  yields
\begin{eqnarray}\label{824}
\phi &=& \phi\Bigl |_{\rho=1}-a \delta \left(\frac{d\phi}{d\rho}\right)\Bigl |_{\rho=1} + O(\delta^2) \nonumber\\
&\sim& x\phi_x\Bigl |_{\rho=1}+y\phi_y\Bigl |_{\rho=1}-\log(s+1)-a\delta(\phi_s s_\rho+ \phi_t t_\rho+\phi_\rho)\Bigl |_{\rho=1}.
\end{eqnarray}
We can easily check that when $\rho=1$,  $\phi_x = \Psi_x$, $\phi_y = \Psi_y$, and (\ref{x}) and (\ref{y}) reduce to (\ref{817}) and (\ref{818}) respectively.
After some calculation we find that when $\rho = 1$,
\begin{eqnarray}
\phi_s s_\rho+ \phi_t t_\rho+\phi_\rho
&=& \phi_s s_\rho+ \phi_t t_\rho+x_\rho \phi_x+x\phi_{x\rho}+y_\rho\phi_y+y\phi_{y\rho}+1 \nonumber\\
&=& \frac{s^2t^2}{2(s+1)}-st+s+1.
\end{eqnarray}
Thus (\ref{824}) becomes, as $\delta\rightarrow 0$,
\begin{equation}\label{826}
\phi \sim x\Psi_x+y\Psi_y-\log(s+1)+a\delta\left[st -\frac{s^2t^2}{2(s+1)}-(s+1)\right],
\end{equation}
and expanding (\ref{769}) as $\delta\rightarrow 0$ gives
\begin{equation}\label{825}
K \sim \delta a |\jmath|^{-1/2}e^{t/2}\frac{s^{3/2}}{(s+1)^2}\frac{(s-st+1)^{3/2}}{\sqrt{se^t-st+1}}.
\end{equation}
Therefore, from (\ref{Psi}), (\ref{823}), (\ref{826}), and (\ref{825}), the matching suggests that
\begin{equation*}
{\cal K}_0(s) =  \frac{a s^{3/2}}{(s+1)^2}\; e^{-a(s+1)},
\end{equation*}
and thus we obtain (\ref{833}).

The ray solution ${\cal P}(x, y) \sim \delta {\cal K}(x, y)\exp[\Psi(x, y)/\delta]$ is valid in the interior of the state space and also near $x=1$.  However, the expansion breaks down near the boundaries $x=0$ and $y=0$, and near the corner points $(x, y)=(1, 0)$ and $(0,0)$.  We analyze these regions separately in the following subsections.

Near $x=1$ we can simplify the ray solution to
\begin{equation}\label{854}
\pi(k, r) = {\cal P}(x, y) \sim \frac{\delta ay}{y+1} \exp[m(x-1)\log(y+1)-a(y+1)],
\end{equation}
so that the exponent is $O(1)$ for $k=m-O(1)$ (i.e.,  $1-x =O(m^{-1})=O(\delta)$).  The total probability mass in this region is given by
\begin{eqnarray*}
{\cal M}_+ &=& Prob\left[N_1=m-O(1), \; N_2=O(m)\right] \\
&\sim& a \int_0^\infty \sum_{n=0}^\infty (y+1)^{-n} \frac{y}{y+1} e^{-a(y+1)}  dy = e^{-a}.
\end{eqnarray*}
In subsection \textbf{5.2} we shall consider the scale $0<x<1$ and $r=O(1)$, and obtain the mass in this range as
$${\cal M}_- = Prob\left[0\leq N_1 \leq m, \; N_2 = 0\right] \sim 1-e^{-a},$$
which complements ${\cal M}_+$.  Thus for $m\rightarrow \infty$ and in the heavy traffic limit, there will either be no secondary spaces occupied, or the number of occupied primary spaces will be close to the maximum $m$ and then there will be a large number ($O(m)$) of secondary spaces occupied.  From (\ref{854}) we can also infer the conditional limit laws
\begin{eqnarray*}
Prob\left[N_1=m-n|N_2 = my\right] &\sim& y(y+1)^{-n-1}, \; \; y>0, \\
Prob\left[N_2=my|N_1=m-n\right] &\sim& \frac{y(y+1)^{-n-1}e^{-ay}}{\mathcal{C}_n},
\end{eqnarray*}
with $$\mathcal{C}_n = \int_0^\infty \frac{u e^{-au}}{(u+1)^{n+1}}  du.$$
This may be written in terms of incomplete Gamma functions as
\begin{equation*}
\mathcal{C}_n = e^a a^n \left[\frac{1}{a}\; \Gamma(1-n, a)-\Gamma(-n, a)\right].
\end{equation*}

\subsection{Boundary layer near $x=0$}

We consider the scale $k=O(1)$ and $y>0$, which corresponds to $x=O(\delta)$ and $r=O(\delta^{-1})= O(m)$.  From (\ref{Psix}) and (\ref{817}) we see that $\Psi_x$ has a logarithmic singularity as $x\rightarrow 0$.  We set
\begin{equation}
\pi(k, r) = \delta^{\hat{\nu}-k}\exp\left[\frac{1}{\delta}\Psi(0,y)\right]{\cal S}_k(y; \delta),
\end{equation}
and rewrite (\ref{bl2}) as
\begin{multline}\label{827}
(2-a\delta)[y +(k+1)\delta]e^{\Psi(0,y)/\delta}{\cal S}_k(y; \delta) = \delta(1-a\delta)[y+(k+1)\delta]e^{\Psi(0,y)/\delta}{\cal S}_{k-1}(y; \delta)\\+(k+1)e^{\Psi(0,y)/\delta}{\cal S}_{k+1}(y; \delta)
 +(y+\delta) e^{\Psi(0,y+\delta)/\delta}{\cal S}_k(y+\delta; \delta), \\ \; k>0, \; y> 0.
\end{multline}
The boundary condition at $k=0$ can be written as
\begin{multline}\label{828}
(2-a\delta)(y+\delta)e^{\Psi(0,y)/\delta}{\cal S}_0(y; \delta) = e^{\Psi(0,y)/\delta}{\cal S}_1(y; \delta)+(y+\delta) e^{\Psi(0,y+\delta)/\delta}{\cal S}_0(y+\delta; \delta).
\end{multline}
We expand ${\cal S}_k(y; \delta)$ in the form
$${\cal S}_k(y; \delta) = {\cal S}_k(y)+\delta {\cal S}_k^{(1)}(y)+O(\delta^2).$$
and let $\delta\rightarrow 0$ in (\ref{827}) and (\ref{828}) to obtain the following equation for the leading term ${\cal S}_k(y)$
\begin{eqnarray} \label{829}
\left[2-e^{\Psi_y(0,y)}\right] y {\cal S}_k(y) = (k+1){\cal S}_{k+1}(y),  \; k>0, \; y> 0.
\end{eqnarray}
The general solution to (\ref{829}) is
\begin{equation*}
{\cal S}_k(y) = {\cal S}_0(y) \frac{y^k}{k!}\left[2-e^{\Psi_y(0,y)}\right]^k
\end{equation*}
and thus
\begin{equation} \label{830}
\pi(k, r) \sim \delta^{\hat{\nu}-k}{\cal S}_0(y) \frac{y^k}{k!}\left[2-e^{\Psi_y(0,y)}\right]^k\exp\left[\frac{1}{\delta}\Psi(0,y)\right].
\end{equation}
In the same way as in section \textbf{4} we determine the constant $\hat{\nu}$ and the function ${\cal S}_0(y)$ by matching (\ref{830}) asymptotically to the ray solution ${\cal P}(x, y) \sim \delta {\cal K}(x, y)e^{\Psi(x, y)/\delta}$ in an intermediate limit ($x\rightarrow 0$ but $k=x/\delta \rightarrow \infty$). Using Stirling's formula for $k!$ we can write the expansion in (\ref{830}) for $k\rightarrow \infty$ in terms of $x$, as
\begin{equation}\label{832}
\pi(k,r) \sim  \delta^{\hat{\nu}+1/2} \frac{{\cal S}_0(y)}{\sqrt{2\pi x}}\exp\left[\frac{1}{\delta}\left\{\Psi(0,y)+x+x\log \left(\frac{y}{x}\right)+x\log\left[2-e^{\Psi_y(0,y)}\right]\right\}\right].
\end{equation}

We expand the ray solution ${\cal P}(x, y)$ as $x\rightarrow 0$, inverting the transformation from $(x, y)$ to ray coordinates along $x=0$.  When $x=0$ we have $t=t_{max}(s)=1/s$.  Then from (\ref{818}) we find that $s$ and $y$ are related by $y=s^2/(s+1).$ We thus define $\hat{s}=\hat{ s}(y)$ by
\begin{equation}
\hat{s}(y) = \frac{1}{2}\left[ y+\sqrt{y^2+4y} \right].
\end{equation}
Note that this is the same as (\ref{734}) with $\rho=1$.
Thus, a ray that starts from $(x, y)=(1, \hat{s})$ hits the $y-$axis at the point $(0, y)$.  Using (\ref{Psi}) and (\ref{Psiy}) and evaluating these expressions at $t=t_{max}$ and $s=\hat{s}$ we obtain explicit expressions for $\Psi(0,y)$ and $\Psi_y(0,y)$:
\begin{eqnarray*}
\Psi(0,y) &=& -\frac{\hat{s}}{\hat{s}+1}-\log(\hat{s}+1), \\
\Psi_y(0,y) &=& -t_{max}= -\frac{1}{\hat{s}},
\end{eqnarray*}
which can be rewritten as (\ref{831}) and (\ref{867}) respectively.
Also, from (\ref{Psix}) and (\ref{817}) we obtain
\begin{equation}\label{855}
\Psi_x = -\log x + \log \left[\frac{(1+s- st)(1+2s-st-s e^{-t})}{s+1}\right],
\end{equation}
which is exact for all $(x, y)$ and indicates the logarithmic singularity in $\Psi_x$ as $x\rightarrow 0$.  By expanding (\ref{855}) for $s\rightarrow \hat{s}$ and $t\rightarrow t_{max}$ we obtain
\begin{eqnarray}
\Psi_x &=& -\log x + \log \left[\frac{\hat{s}^2(2-e^{-t})}{\hat{s}+1}\right]+o(1) \nonumber\\
&=&-\log x + \log y + \log\left[2-e^{\Psi_y(0,y)}\right]+o(1),
\end{eqnarray}
which integrates to
\begin{equation}
\Psi(x, y) = \Psi(0, y)-x\log x +x+x\log y +x \log\left[2-e^{\Psi_y(0,y)}\right]+o(x).
\end{equation}
The above agrees precisely with the exponential terms in (\ref{832}) and thus the expansions will asymptotically match provided that
$\hat{\nu} = 1/2$ and ${\cal S}_0(y) = \lim_{x \rightarrow 0}\left[\sqrt{2\pi x} \; {\cal K}(x,y)\right]$.
By expanding $\jmath=x_t y_s-x_s y_t$ and then ${\cal K}$ in (\ref{833}) for $s\rightarrow \hat{s}$ and $t\rightarrow t_{max}$, we obtain \begin{equation}
{\cal S}_0(y) = \frac {a \hat{s}^{3/2}\sqrt{2\pi}}{(\hat{s}+1)\sqrt{\hat{s}+2}}\exp\left[-a\left(\hat{s}+\frac{1}{2(\hat{s}+1)}\right)\right],
\end{equation}
which can also be written as (\ref{868}) in terms of $y$.
We have thus determined $\hat{\nu}$ and ${\cal S}_0(y)$ in (\ref{830}), and verified the matching between the $(x, y)$ and $(k, y)$ scales.

\subsection{Boundary layer near $y=0$}
We consider small values of $y$ and use the original discrete
variable $r=y/\delta$.  We set, for $0<x<1$,
\begin{equation}\label{834}
\pi(k,r) = m^{\nu_1+r/2}e^{-\sqrt{mr}{\cal H}(x)}{\cal Q}_r(x;\delta) \; \; r\geq 1
\end{equation}
and
\begin{equation}\label{835}
\pi(k,0) = (1-\rho)\rho^k+m^{\nu_1}\sqrt{m}e^{-\sqrt{m}{\cal H}(x)}{\cal Q}_0(x;\delta), \; \; r=0,
\end{equation}
where ${\cal H}(x)$ is a positive function of $x$.  The forms of the expansions in (\ref{834}) and (\ref{835}) are indicated by the behavior of the ray expansion as $y\rightarrow 0$, for $0<x<1$, which we shall discuss shortly.  We also note that in the present heavy traffic limit
\begin{equation}
(1-\rho)\rho^k = \delta a (1-\delta a)^{x/\delta}= \delta a e^{-ax}\left[1-\frac{a^2 x}{2}\delta +O(\delta^2)\right].
\end{equation}
With (\ref{834}) and
(\ref{835}), from (\ref{bl2}), we obtain
\begin{multline}
(2-a\delta)[x+\delta(r+1)]e^{-\sqrt{mr}{\cal H}(x)}{\cal Q}_r(x; \delta) = (1-a\delta)[x+\delta(r+1)]e^{-\sqrt{mr}{\cal H}(x-\delta)}{\cal Q}_r (x-\delta; \delta)\\
+(x+\delta)e^{-\sqrt{mr}{\cal H}(x+\delta)}{\cal Q}_r (x+\delta; \delta)
 +\delta(r+1)\sqrt{m}e^{-\sqrt{m(r+1)}{\cal H}(x)}{\cal Q}_{r+1}(x; \delta), \;   \;   r > 0. \label{836}
\end{multline}
and
\begin{multline}
(2-a\delta)(x+\delta)e^{-\sqrt{m}{\cal H}(x)}{\cal Q}_0(x; \delta) = (1-a\delta)(x+\delta)e^{-\sqrt{m}{\cal H}(x-\delta)}{\cal Q}_0 (x-\delta; \delta) \\
 +(x+\delta)e^{-\sqrt{m}{\cal H}(x+\delta)}{\cal Q}_0 (x+\delta; \delta)
 +\delta e^{-\sqrt{m}{\cal H}(x)}{\cal Q}_1 (x; \delta), \;  \;  r=0. \label{837}
\end{multline}
Expanding ${\cal Q}_r(x; \delta)$  and ${\cal Q}_0(x; \delta)$  as
\begin{equation*}
{\cal Q}_r(x; \delta) =  {\cal Q}_r(x)+\delta {\cal Q}_r^{(1)}(x)+O(\delta^2), \; \; r\geq 0,
\end{equation*}
we obtain from (\ref{836}) at the first two orders ($O(\delta)$ and $O(\delta^{3/2})$) the following
equations
\begin{eqnarray}
[{\cal H}'(x)]^2 &=& \frac{1}{x}, \label{838}\\
2{\cal H}'(x) {\cal Q}_r'(x) &=& \left[\; \frac{r}{x}{\cal H}'(x)-a{\cal H}'(x)-{\cal H}''(x)\right]{\cal Q}_r(x). \label{839}
\end{eqnarray}
Taking  ${\cal H}'(x)<0$, from (\ref{838}) we obtain ${\cal H}(x)=-2\sqrt{x}+constant$.  By matching to the ray solution at $x=1$ (corresponding to $t=0$), where $\Psi=0$, we must have  ${\cal H}(1)=0$ so that
\begin{equation}\label{840}
{\cal H}(x) = 2\left(1-\sqrt{x}\right).
\end{equation}
Using (\ref{840}) we solve (\ref{839}) to obtain
\begin{equation}
{\cal Q}_r(x)= {\cal Q}_*(r)x^{r/2+1/4}e^{-ax/2},  \;  r>0.
\end{equation}
Here ${\cal Q}_*(r)$ is a function of $r$. Thus, for $r>0$, the boundary solution near $y=0$ is
\begin{equation}\label{849}
\pi(k,r) \sim {\cal Q}_*(r) \; m^{\nu_1+r/2}x^{r/2+1/4}e^{-ax/2}e^{-2\sqrt{mr}(1-\sqrt{x})}.
\end{equation}

We determine ${\cal Q}_*(r)$ and $\nu_1$  by matching this boundary layer solution to the ray solution ${\cal P}(x, y)$.  We expand the ray solution ${\cal P}(x, y)$  as $y\rightarrow 0$.  We let  $s\rightarrow 0$, $\; t\rightarrow \infty$ with $st$  fixed and  $0<st<1$.  Then (\ref{Psi}), (\ref{817}), (\ref{818}), and (\ref{833})  give
\begin{eqnarray}
x &\sim& \frac{(1-st)^2}{s+1}+\frac{2s(1-st)}{s+1}, \label{841}\\
y &\sim& \frac{s^2}{s+1}, \label{842}\\
{\cal K} &\sim&  \frac{a}{2}\sqrt{1-st} \exp \left[-a\left(1-st +\frac{s^2t^2}{2}\right)\right], \label{843} \\
\Psi &\sim& (1-st)^2\log\left(1+\frac{s}{1-st}\right)+s^2\log\left(\frac{1-st}{se^t}\right)-\log(1+s) \nonumber\\
&\sim& s^2\left[\log(1-st)-\log s-2t\right]. \label{844}
\end{eqnarray}
We can invert (\ref{841}) and (\ref{842}) and write $s$ and $t$ in terms of $x$ and $y$, for $y\rightarrow 0$, as
\begin{eqnarray}
s &=& \sqrt{y}+\frac{y}{2}+O(y^{3/2}), \label{845}\\
1-st &=& \sqrt{x}+\left(\frac{\sqrt{x}}{2}-1\right)\sqrt{y}+O(y), \label{846}\\
t &=& \frac{1-\sqrt{x}}{\sqrt{y}}+\frac{1}{2}+O(\sqrt{y}). \label{847}
\end{eqnarray}
Using (\ref{845}) - (\ref{847}) in (\ref{843}) and (\ref{844}) we then obtain
\begin{eqnarray}
\delta {\cal K} &\sim& \frac{\delta a}{2}x^{1/4}e^{-a(x+1)/2}, \\
\Psi &\sim& -2\sqrt{y}(1-\sqrt{x})+\frac{y}{2}\log\left(\frac{x}{y}\right),
\end{eqnarray}
and thus near $y=0$ the ray expansion becomes
\begin{equation}\label{848}
{\cal P}(x, y) \sim \frac{a}{2}\; m^{r/2-1}r^{-r/2}x^{r/2+1/4}e^{-a(x+1)/2}e^{-2\sqrt{mr}(1-\sqrt{x})}.
\end{equation}
Therefore, comparing (\ref{849}) to (\ref{848}) we conclude that ${\cal Q}_*(r) = \frac{a}{2}e^{-a/2}r^{-r/2}$ and $\nu_1=-1$.
We have thus obtained (\ref{861}).

For $r=0$, we rewrite (\ref{835}) as
\begin{equation}
\pi(k,0) \sim (1-\rho)\rho^k+\frac{e^{-2\sqrt{m}(1-\sqrt{x})}}{\sqrt{m}}{\cal Q}_0(x).
\end{equation}
Then, from (\ref{837}), we obtain to the leading order
\begin{equation}
{\cal Q}_0(x) = -{\cal Q}_1 (x) = -{\cal Q}_*(1)x^{3/4}e^{-ax/2} = -\frac{a}{2}x^{3/4}e^{-a(x+1)/2},
\end{equation}
and thus obtain (\ref{876}).

\subsection{Corner layer near $(x,y)=(0,0)$}
We consider  $(x,y)$  near  $(0,0)$ and go back to the original discrete variables $(k,r)$, with $k, r = O(1)$.  We also set
\begin{eqnarray}
\pi(k,r) &=& m^{-5/4}e^{-2\sqrt{mr}}\;  {\cal T}(k,r;\; m), \; \; r>0,  \label{850}\\
\pi(k,0) &=& (1-\rho)\rho^k+ m^{-5/4}e^{-2\sqrt{m}}\; {\cal T}(k,0;\; m), \; \; r=0, \label{851}
\end{eqnarray}
and  note that ${\cal T}(0, 0;\; m)=0$.  The scaling in (\ref{850}) and (\ref{851}) can be inferred by expanding (\ref{861}) and (\ref{876}) as $x=\delta k=k/m \rightarrow 0$.
Then the balance equation (\ref{bl2}) becomes
\begin{multline}
(2-a\delta)(k+r+1)e^{-2\sqrt{mr}}\; {\cal T}(k,r;\; m) = (1-a\delta)(k+r+1)e^{-2\sqrt{mr}}\; {\cal T}(k-1,r;\; m)\\
+(k+1) e^{-2\sqrt{mr}}\; {\cal T}(k+1,r;\; m)
 + (r+1)e^{-2\sqrt{m(r+1)}}\; {\cal T}(k,r+1;\; m), \label{852}
\end{multline}
for $\;  k>0, \; \; r > 0$, and
\begin{multline}
(2-a\delta)(k+1)\; {\cal T}(k,0;\; m) = (1-a\delta)(k+1)\; {\cal T}(k-1,0;\; m)\\
+ (k+1)\; {\cal T}(k+1,0;\; m) + {\cal T}(k,1;\; m), \;   \;  \;  k >0, \;  \; r = 0. \label{853}
\end{multline}
The boundary condition at $k=0$ is
\begin{equation}\label{856}
(2-a\delta)\; e^{-2\sqrt{mr}}\;  {\cal T}(0,r;\; m) = \frac{1}{r+1}\; e^{-2\sqrt{mr}}\; {\cal T}(1,r;\; m)
 + e^{-2\sqrt{m(r+1)}}\; {\cal T}(0,r+1;\; m),
\end{equation}
for $r>0$, and the corner condition at $(k, r)=(0, 0)$ is
\begin{equation}\label{873}
0 = {\cal T}(1,0;\; m)+{\cal T}(0,1;\; m).
\end{equation}

Expanding ${\cal T}(k,r;\; m)$ as ${\cal T}(k,r;\; m)={\cal T}(k,r)+o(1)$ we obtain from (\ref{852}) - (\ref{856}) the following equations for the leading term
\begin{multline}
2(k+r+1)\; {\cal T}(k,r) = (k+r+1)\; {\cal T}(k-1,r)+(k+1)\; {\cal T}(k+1,r), \\     k>0, \;   \;   r>0,
\label{857}
\end{multline}
\begin{multline}
2(k+1)\; {\cal T}(k,0) = (k+1)\; {\cal T}(k-1,0) + (k+1)\; {\cal T}(k+1,0) + {\cal T}(k,1),\\  k>0,  \;  \; r=0,\label{858}
\end{multline}
\begin{equation}
2(r+1)\; {\cal T}(0,r) = \; {\cal T}(1,r), \;  \;   \;  r> 0. \label{859}
\end{equation}
We shall first analyze (\ref{857}) and (\ref{859}) for $r>0$, and then solve (\ref{858}) for ${\cal T}(k,0)$.

For $r>0$ we use the generating function ${\cal G}(z)= {\cal G}(z;r)\equiv \sum_{k=0}^\infty {\cal T}(k,r) z^k$.
Multiplying (\ref{857}) by $z^k$ and summing over  $k\geq 0$ gives
\begin{equation} \label{860}
(z-1)^2{\cal G}'(z) = [-(r+2)z+2(r+1)]{\cal G}(z).
\end{equation}
The general solution to (\ref{860}) is
\begin{equation*}
{\cal G}(z;r) = {\cal G}_*(r)(1-z)^{-r-2}e^{r/(1-z)}
\end{equation*}
where ${\cal G}_*(r)$ is a function of $r$.
Inverting the generating function, we can write
\begin{equation}\label{862}
{\cal T}(k,r) = {\cal G}_*(r) \frac{1}{2\pi i}\oint \frac{z^{-k-1}}{(1-z)^{r+2}}\; e^{r/(1-z)} dz, \; \; r>0
\end{equation}
where the integral is a complex contour integral along a small loop around $z=0$.
We obtain ${\cal G}_*(r)$ by asymptotic matching. Letting $x\rightarrow 0$ ($x=k/m$) in (\ref{861}), and invoking the asymptotic matching condition shows that
\begin{equation}\label{863}
{\cal T}(k,r) \sim  \frac{a}{2}\; r^{-r/2}k^{r/2+1/4}e^{-a/2}e^{2\sqrt{kr}}, \; \; \; k\rightarrow \infty.
\end{equation}
We now expand (\ref{862}) as $k\rightarrow \infty$.  This will verify that (\ref{863}) is satisfied and also determine the function ${\cal G}_*(r)$.
As $k\rightarrow \infty$ the behavior of the integral in (\ref{862}) will be determined by a saddle point near $z=1$.  We write the integrand in (\ref{862}) as
\begin{equation}
\frac{\exp[\hat{H}(z)]}{z(1-z)^{r+2}} \equiv \frac{1}{z(1-z)^{r+2}}\exp\left[-k\log z+\frac{r}{1-z}\right].
\end{equation}
The saddle point(s) satisfy $\hat{H}'(z)=0$ so that
\begin{equation}\label{864}
\frac{k}{z} = \frac{r}{(1-z)^2}.
\end{equation}
(\ref{864}) defines the saddle $\hat{z}$ as a function of $k/r$, and it is given explicitly by
\begin{equation}
\hat{z} = 1+\frac{r}{2k}-\frac{\sqrt{r(4k+r)}}{2k} \sim 1- \sqrt{\frac{r}{k}}, \;  k\rightarrow \infty.
\end{equation}
Evaluating the integral along the steepest descent path through the saddle point we obtain, as $k\rightarrow \infty$,
\begin{equation}\label{865}
{\cal T}(k,r) \sim \frac{1}{2\sqrt{\pi}} r^{-r/2-3/4}k^{r/2+1/4} e^{r/2}e^{2\sqrt{kr}}{\cal G}_*(r).
\end{equation}
Therefore, (\ref{863}) and (\ref{865}) suggest that ${\cal G}_*(r) = a\sqrt{\pi}\; r^{3/4}e^{-(a+r)/2}$, and thus, for $r\geq 1$, we obtain (\ref{866}).

Next we show that (\ref{866}), when expanded for $r\rightarrow \infty$ asymptotically matches to the expansion in the boundary solution near $x=0$, as $y\rightarrow 0$.
Letting $y\rightarrow 0$ in (\ref{831}) - (\ref{868}) we obtain
\begin{eqnarray*}
{\cal S}_0(y) &=&   a\sqrt{\pi}y^{3/4} e^{-a/2}+O(y^{5/4}) \\
\Psi(0,y) &=&-2\sqrt{y} + \frac{y}{2}+O(y^{3/2})\\
\Psi_y(0,y) &=&   -\frac{1}{\sqrt{y}}+O(1).
\end{eqnarray*}
Thus as $y \rightarrow 0$ ($y=r/m$) the asymptotic behavior of (\ref{830}) is
\begin{equation}\label{869}
\pi(k, r) \sim a\sqrt{\pi}\; m^{-5/4}r^{3/4} e^{-a/2} e^{r/2} e^{-2\sqrt{mr}}\; \frac{2^k r^k}{k!}.
\end{equation}
This is an approximation to $\pi(k,r)$ that applies in the matching region where $k=O(1), \; \; r\rightarrow \infty, \; \;  r/m\rightarrow 0.$

Next we expand the corner approximation (\ref{866}) as $r\rightarrow \infty$.  Expanding the integrand in (\ref{866}) around $z=0$ by setting $z=u/r$ and evaluating the resulting integral, using
\begin{equation*}
\frac{r^k e^r}{2\pi i}\oint u^{-k-1} e^{2u} du = \frac{2^k r^k e^r}{k!},
\end{equation*}
shows that (\ref{866}), as $r\rightarrow \infty$ with $k=O(1)$, agrees with (\ref{869}).

We verify that the corner approximation (\ref{866}) matches to the ray solution $\delta {\cal K}(x, y)\exp[m \Psi(x, y)]$, in an intermediate limit where $k, r \rightarrow \infty$ but $x=k/m, \; y=r/m \rightarrow 0$.  We expand (\ref{866}) for $k$ and $r$ simultaneously large, writing the integrand as
\begin{equation*}
\frac{\exp[\tilde{H}(z)]}{z(1-z)^2} \equiv \frac{1}{z(1-z)^2} \exp\left[-k \log z-r\log(1- z)+\frac{r}{1-z}\right]
\end{equation*}
and again using the saddle point method.  The saddle point(s) now satisfy $\tilde{H}'(z)=0$ so that
\begin{equation} \label{870}
\frac{k}{r} = \frac{z}{1-z}\left(1+\frac{1}{1-z}\right).
\end{equation}
(\ref{870}) defines the saddle $\tilde{z}$ as a function of $k/r$, and it is given explicitly by
\begin{equation*}
\tilde{z} = 1-\sqrt{\frac{r}{k+r}}.
\end{equation*}
We note that $\tilde{z} \rightarrow 0$ as $k/r \rightarrow 0$ and $\tilde{z} \rightarrow 1$ as $k/r \rightarrow \infty$.
Evaluating the integral along the steepest descent path through the saddle point we obtain
\begin{eqnarray}
\pi(k,r) &\sim& \frac{a}{2}\; m^{-5/4}r^{1/4}e^{-(a+r)/2}e^{-2\sqrt{mr}}e^{r/(1-\tilde{z})}\frac{\tilde{z}^{-k}(1-\tilde{z})^{-r}}{\sqrt{\tilde{z}(1-\tilde{z})}} \nonumber\\
&=& \frac{a}{2}\; m^{-5/4}r^{1/4}e^{-(a+r)/2}e^{-2\sqrt{mr}}e^{\sqrt{r(k+r)}}\nonumber\\
&\times&\left(1-\sqrt{\frac{r}{k+r}}\right)^{-k-1/2}\left(\frac{k+r}{r}\right)^{r/2+1/4}.\label{871}
\end{eqnarray}
This approximation applies in the matching region where $k, r \rightarrow \infty$, but $k/m, r/m \rightarrow 0$.

In the ray solution we let $x \rightarrow 0$ and $y\rightarrow 0$, which corresponds to $s \rightarrow 0$ and $\; t\rightarrow \infty$, with $x/y$ fixed (corresponding to $1-st \approx s$).  We find that in this limit
$$1-st \sim \sqrt{x+y}-\sqrt{y}, \; \;  s\sim \sqrt{y}, \; \; t \sim \frac{1}{\sqrt{y}}-\frac{\sqrt{x+y}}{\sqrt{y}}+1,$$
and thus, from (\ref{Psi}) - (\ref{Psiy}) and (\ref{833}),  obtain
\begin{eqnarray*}
\Psi_x &\sim& \log\left(\frac{\sqrt{x+y}}{\sqrt{x+y}-\sqrt{y}}\right),\\
\Psi_y &\sim& \log\left(\frac{\sqrt{x+y}}{\sqrt{y}}\right)-\frac{1}{\sqrt{y}}+\frac{\sqrt{x+y}}{\sqrt{y}}-1, \\
\Psi &\sim& x\log\left(\frac{\sqrt{x+y}}{\sqrt{x+y}-\sqrt{y}}\right)-2\sqrt{y}+\sqrt{y(x+y)}-\frac{y}{2}+ y\log\left(\frac{\sqrt{x+y}}{\sqrt{y}}\right),\\
{\cal K} &\sim& \frac{a}{2}e^{-a/2}\frac{\sqrt{x+y}}{\left(\sqrt{x+y}-\sqrt{y}\right)^{1/2}}.
\end{eqnarray*}
Therefore, as $(x, y)\rightarrow (0,0)$, in terms of $x=k/m$ and $y=r/m$,
\begin{eqnarray}\label{872}
{\cal P}(x,y) &\sim& \frac{a}{2}\; m^{-1}e^{-(a+r)/2}e^{-2\sqrt{mr}}e^{\sqrt{r(k+r)}}\nonumber\\
&\times&\frac{\sqrt{k+r}}{\left(\sqrt{k+r}-\sqrt{r}\right)^{1/2}}\left(\frac{\sqrt{k+r}}{\sqrt{k+r}-\sqrt{r}}\right)^k\left(\frac{\sqrt{k+r}}{\sqrt{r}}\right)^r.
\end{eqnarray}
We can easily check that (\ref{872}) agrees precisely with (\ref{871}).

It remains to compute ${\cal T}(k,0)$, by solving (\ref{858}) and (\ref{873}).  We use the generating function
$${\cal G}_0(z) \equiv \sum_{k=0}^\infty {\cal T}(k,0) z^k,$$
multiply (\ref{858}) by $z^k$, and sum over $k\geq 0$, to obtain
\begin{equation}\label{874}
(z-1)^2{\cal G}_0'(z) +2(z-1){\cal G}_0(z)+ {\cal G}_*(1)(1-z)^{-3}e^{1/(1-z)}=0.
\end{equation}
The solution to (\ref{874}) is
\begin{equation}
{\cal G}_0(z) = {\cal G}_*(1)\frac{z e^{1/(1-z)}}{(z-1)^3} =-a\sqrt{\pi}\; e^{-(a+1)/2}z (1-z)^{-3} e^{1/(1-z)},
\end{equation}
and thus,
\begin{eqnarray}\label{875}
{\cal T}(k,0) &=& \frac{1}{2\pi i}\oint z^{-k-1}{\cal G}_0(z)dz \nonumber\\
&=& -a\sqrt{\pi}e^{-(a+1)/2}  \frac{1}{2\pi i}\oint \frac{z^{-k}}{(1-z)^3}\; e^{1/(1-z)}dz.
\end{eqnarray}

To verify the matching between (\ref{851}) and (\ref{876}) we evaluate the integral in (\ref{875}) as $k\rightarrow \infty$.  We rewrite the integrand in (\ref{875}) as
\begin{equation*}
\frac{\exp[\tilde{H}_0(z)]}{(1-z)^3} \equiv \frac{1}{(1-z)^3} \exp\left[-k \log z+\frac{1}{1-z}\right],
\end{equation*}
and use the saddle point method. The saddle point(s) $\hat{z}_0$ satisfy $\tilde{H}_0'(z)=0$ so that
$k=z/(1-z)^2$  and hence
\begin{equation*}
\hat{z}_0=1-\frac{\sqrt{4k+1}}{2k}+\frac{1}{2k} \sim 1-\frac{1}{\sqrt{k}}.
\end{equation*}
Evaluating the integral along the steepest descent path through the saddle point and using the result in (\ref{851}) we obtain
\begin{equation}\label{877}
\pi(k,0) \sim (1-\rho)\rho^k- \frac{a}{2}e^{-a/2}m^{-5/4}k^{3/4}e^{-2\sqrt{m}}e^{2\sqrt{k}}.
\end{equation}
As $x \rightarrow 0$ ($x=k/m$) the asymptotic behavior of (\ref{876}) agrees with (\ref{877}).
This completes the analysis of the corner range where $(x, y)\approx(0,0)$.

\subsection{Analysis near the corner $(x, y)=(1, 0)\;$:\\ \; \; $x=1-O(m^{-1/2}), \; y=O(m^{-1/2})$}

We examine the problem for $x = 1-O(\sqrt{\delta})$ and $y=O(\sqrt{\delta})$, which corresponds to $k=m-O(\sqrt{m})$ and $r=O(\sqrt{m})$.
We let
\begin{equation}\label{8ap}
x= 1-\xi \sqrt{\delta}, \; \;  y=R\sqrt{\delta},
\end{equation}
and set
\begin{equation}\label{878}
\pi(k, r) = \delta^{3/2} \Omega(\xi, R;\; \delta).
\end{equation}
The scaling in (\ref{878}) can be inferred by expanding (\ref{854}) as $x=1-\xi \sqrt{\delta} \rightarrow 1$ and $y=R\sqrt{\delta} \rightarrow 0$.
Then, for $0<\xi<1$ and $R>0$,  (\ref{bl2}) becomes
\begin{multline}\label{879}
(2-a\delta)\Omega(\xi, R;\; \delta)=\Omega(\xi-\sqrt{\delta}, R;\; \delta)+(1-a\delta) \Omega(\xi+\sqrt{\delta}, R;\; \delta)\\
+\frac{(R\sqrt{\delta}+\delta)\Omega(\xi, R+\sqrt{\delta};\; \delta)}{1+(R-\xi)\sqrt{\delta}+\delta} -\frac{R\sqrt{\delta}\; \Omega(\xi-\sqrt{\delta}, R;\; \delta)}{1+(R-\xi)\sqrt{\delta}+\delta}.
\end{multline}
The boundary condition (\ref{ab}) becomes
\begin{equation}\label{880}
\frac{1+\delta}{1+R\sqrt{\delta}+\delta}\Omega(-\sqrt{\delta}, R;\; \delta) =
(1-\delta)\;\Omega(0, R-\sqrt{\delta};\; \delta), \; \; R>0.
\end{equation}
Expanding $\Omega(\xi, R;\; \delta)$ as
$$\Omega(\xi, R;\; \delta) = \Omega(\xi, R)+O(\sqrt{\delta}),$$
from (\ref{879}) and (\ref{880})  we obtain to the leading order
\begin{eqnarray}
\Omega + \Omega_{\xi\xi}+R(\Omega_\xi+\Omega_R)=0\; ; \; \; \xi, \; R >0 \label{881}\\
R\; \Omega+\Omega_\xi-\Omega_R=0 \; ; \; \; \xi=0, \; R>0. \label{882}
\end{eqnarray}

We must thus solve a parabolic PDE in the quarter plane, subject to an oblique derivative boundary condition along $\xi = 0$.  While we were not able to solve this problem exactly, we shall use asymptotic matching to infer various properties of $\Omega(\xi, R)$ as $\xi$ and/or $R$ become(s) large, or if $R\rightarrow 0$.  We shall also obtain an integral equation for the boundary values $\Omega(0,R)$, and show that the problem can be reduced to the heat equation with a moving boundary.  We comment that the PDE in (\ref{881}) is not separable due to the term $R\Omega_\xi$.  When considering the analogous infinite server model (see \cite{KB}) Knessl obtained, in a certain heavy traffic limit, a problem very similar to (\ref{881}) and (\ref{882}).  However, there the term $R\Omega_\xi$ was replaced by $\xi \Omega_\xi$, so that the PDE was separable.  Then despite the oblique derivative boundary condition Knessl could solve this problem explicitly, in terms of contour integrals of parabolic cylinder functions.  For the infinite server model the corresponding probability $\pi(k,r)$ on the scale (\ref{8ap}) was $O(m^{-1})$ and most of the probability mass accumulated on this scale.  In the present processor sharing model we have, in view of (\ref{878}),  $\pi(k,r)=O(m^{-3/2})$ for $k=m-O(\sqrt{m})$ and $r=O(\sqrt{m})$.  Thus the total mass is roughly $O(m^{-3/2})\times O(\sqrt{m}) \times O(\sqrt{m}) = O(m^{-1/2})$.  Indeed we already showed that for the PS model with $m\rightarrow \infty$ and $1-\rho = O(m^{-1})$, most mass occurs either along $r=0, \; k=O(m)$ or $r=O(m), \; k=m-O(1)$.  The PS model thus, on the $(\xi, R)$ scale, leads to a more difficult mathematical problem, but one whose analysis is perhaps less critical, due to the small probability mass.

We first examine (\ref{881}) and (\ref{882}) for $\xi$ and/or $R\rightarrow \infty$.  This will establish the asymptotic matching between the $(\xi, R)$ scale and the ray expansion on the $(x, y)$ scale.  We introduce the small parameter $\varepsilon >0$, setting
\begin{equation*}
\xi = \frac{\xi^*}{\varepsilon}, \; \; R=\frac{R^*}{\varepsilon},
\end{equation*}
and
\begin{equation}\label{892}
\Omega = \varepsilon^{\nu^*}e^{{\cal A}(\xi^*, R^*)/\varepsilon^2}{\cal B}(\xi^*, R^*)\left[1+O(\varepsilon^2)\right]
\end{equation}
From (\ref{881}) and (\ref{882}) at the first two orders ($O(\varepsilon^{-2})$ and $O(1)$) we obtain
\begin{eqnarray}
{\cal A}_{\xi^*}^2+R^*({\cal A}_{\xi^*}+{\cal A}_{R^*})=0, \label{883}\\
2{\cal A}_{\xi^*}{\cal B}_{\xi^*}+R^*({\cal B}_{\xi^*}+{\cal B}_{R^*})+({\cal A}_{\xi^*\xi^*}+1){\cal B}=0 \label{890}
\end{eqnarray}
and when $\xi^*=0$  (\ref{882}) leads to
\begin{eqnarray}
{\cal A}_{\xi^*}-{\cal A}_{R^*}+R^*=0, \label{886}\\
{\cal B}_{\xi^*}-{\cal B}_{R^*}=0. \label{891}
\end{eqnarray}
We are thus using a ray expansion on the PDE (\ref{881}).  We again use the method of characteristics, writing (\ref{883}) as $${\cal F}^*(\xi^*, R^*, {\cal A}, {\cal A}_{\xi^*}, {\cal A}_{R^*})=0$$ where
\begin{equation}
{\cal F}^* \equiv {\cal A}_{\xi^*}^2+R^*({\cal A}_{\xi^*}+{\cal A}_{R^*}).
\end{equation}
The characteristic equations are
\begin{eqnarray}
\frac{d \xi^*}{d \omega}&=&-\frac{\partial {\cal F}^*}{\partial {\cal A}_{\xi^*}} = -2 {\cal A}_{\xi^*}-R^*, \label{884}\\
\frac{d R^*}{d\omega} &=& - \frac{\partial {\cal F}^*}{\partial {\cal A}_{R^*}} = -R^*,\\
\frac{d {\cal F}^*}{d\omega} &=& -2{\cal A}_{\xi^*}^2-R^*({\cal A}_{\xi^*}+{\cal A}_{R^*})=-{\cal A}_{\xi^*}^2, \\
\frac{\partial {\cal A}_{\xi^*}}{d\omega} &=& -\frac{\partial {\cal F}^*}{\partial \xi^*} = 0,\\
\frac{\partial {\cal A}_{R^*}}{d\omega} &=& -\frac{\partial {\cal F}^*}{\partial R^*} = -{\cal A}_{\xi^*}-{\cal A}_{R^*}.\label{885}
\end{eqnarray}
Letting $R^*=u$ at $\omega=0$ (using rays start from ($R^*, \xi^*$) = ($u, 0$) at $\omega=0$) from (\ref{884})-(\ref{885}) and (\ref{886}) we obtain
\begin{eqnarray}
R^*&=& \varepsilon R =u e^{-\omega}, \label{887} \\
\xi^* &=& \varepsilon \xi = u \left[ 2\omega + e^{-\omega}-1\right], \label{888}\\
{\cal A} &=& -u^2 \omega. \label{889}
\end{eqnarray}
Here we also used ${\cal A}(0, 0)=0$ to determine an integration constant.
We rewrite (\ref{890}) as
\begin{equation}\label{893}
\frac{d{\cal B}}{d\omega}=({\cal A}_{\xi^*\xi^*}+1){\cal B}
\end{equation}
The general solution to (\ref{891}) and (\ref{893}) is
\begin{equation}\label{894}
{\cal B} = {\cal B}_0 \frac{u e^{\omega}}{\sqrt{2\omega+1}},
\end{equation}
where ${\cal B}_0$ is a constant.

To determine ${\cal B}_0$  we examine the asymptotic behavior of the solution (\ref{892}) as $\xi \rightarrow 0$ and $R \rightarrow \infty$, which corresponds to $u\rightarrow \infty$ and $\omega \rightarrow 0$.  From (\ref{887})-(\ref{889}) and (\ref{894}) we obtain
\begin{equation*}
\varepsilon R \sim u, \; \; \varepsilon \xi \sim u \omega, \; \;
{\cal A} \sim -\varepsilon^2 R\xi, \; \; {\cal B} \sim \varepsilon {\cal B}_0 R.
\end{equation*}
Thus, from (\ref{878}) and (\ref{892})
\begin{equation}
\pi(k, r) \sim \delta^{3/2}\varepsilon^{\nu^*+1}{\cal B}_0 R e^{-R\xi}.
\end{equation}
Letting $y \rightarrow 0$ ($y=\sqrt{\delta}R$) and $x=1-\sqrt{\delta}\xi$ in (\ref{854}) gives
\begin{equation}\label{8aq}
\pi(k, r) \sim \delta^{3/2} aR e^{-a}e^{-\xi R},
\end{equation}
and matching thus forces $\nu^* = -1$ and ${\cal B}_0 = a e^{-a}$.
Therefore,
\begin{eqnarray}
\delta^{3/2} \Omega(\xi, R) &\sim& \delta^{3/2}a e^{-a}\frac{ \varepsilon^{-1} u e^\omega}{\sqrt{2\omega+1}}\exp\left(-\varepsilon^{-2} u^2 \omega\right)\label{84}\\
&=& \delta^{3/2}ae^{-a}\frac{ R e^{2\omega}}{\sqrt{2\omega+1}}\exp\left(-R^2\omega e^{2\omega}\right). \label{895}
\end{eqnarray}
Here we used $u=\varepsilon R e^\omega$.  This gives the behavior of $\Omega(\xi, R)$ as $\xi, R \rightarrow \infty$ with $0\leq \xi/R <\infty$.

We can rewrite (\ref{895}) in terms of the Lambert $W$-function.  From (\ref{887}) and (\ref{888}) we obtain
$\xi/R = e^{\omega} \left(2\omega+e^{-\omega}-1\right)= e^\omega (2\omega-1)+1$, so that
\begin{equation}
\frac{\xi}{R}-1 = 2\sqrt{e}\left(\omega-\frac{1}{2}\right)e^{\omega-1/2}. \label{899}
\end{equation}
Setting  $W\equiv\omega -1/2$  (\ref{899}) becomes  $$We^W = \frac{1}{2\sqrt{e}}\left(\frac{\xi}{R}-1\right),$$
so that $\omega = 1/2+W$, where $W$ is the Lambert function.  Thus (\ref{895}) becomes
\begin{equation}\label{82}
\delta^{3/2}\Omega(\xi, R) \sim \delta^{3/2}ae^{-a}\frac{(\xi- R)^2 }{4R\sqrt{2(W+1)}W^2}\exp\left[-(\xi-R)^2\left(\frac{W+1/2}{4W^2}\right)\right].
\end{equation}

We next show that (\ref{82}) or (\ref{895}) also follows by expanding the ray expansion on the $(x, y)$ scale, as $(x, y) \rightarrow (0,0)$ along lines where $y/x$ is fixed.
We expand the ray solution ${\cal P}(x, y)$ as $s\rightarrow 0$ with $t$ fixed.
Letting $s\rightarrow 0$ with $t$ fixed in (\ref{817}) and (\ref{818}) gives
\begin{eqnarray*}
x(s,t)&=& 1+s\left(1-2t-e^{-t}\right) + O(s^2), \\
y(s,t) &=& s e^{-t} +O(s^2),
\end{eqnarray*}
and thus $ye^t \sim s$ and
\begin{equation}\label{896}
s \sim \frac{1-x-y}{2t-1}, \;  \; \; \frac{1-x}{y}-1 \sim 2\sqrt{e}\left(t-\frac{1}{2}\right)e^{t-1/2}.
\end{equation}
From (\ref{896}), (\ref{Psi}) - (\ref{Psiy}), and (\ref{823}) we obtain
\begin{eqnarray}
\Psi &\sim& x\left[s+s^2\left(t-\frac{1}{2}\right)\right]+y\left[s\left(1-e^t\right)+s^2\left(t-\frac{1}{2}\right)\right]-s+\frac{s^2}{2}\nonumber\\
&\sim& \left(\frac{s^2}{2}-s\right)(1-x-y) + s^2[(x+y)t- 1] \nonumber\\
&\sim& -\frac{(1-x-y)^2}{(2t-1)^2}\left[t+(1-x-y)t-\frac{1-x-y}{2}\right],\label{897}
\end{eqnarray}
and
\begin{equation}
{\cal K} \sim  \frac{a e^{-a} s  e^t}{\sqrt{2t+1}} \sim \frac{a e^{-a}}{\sqrt{2t+1}}\frac{(1-x-y)^2}{y(2t-1)^2}. \label{898}
\end{equation}
Since $(1-x)/y = \xi/R$, in view of (\ref{899}) and (\ref{896})  we can identify $t$, as $t\sim \omega$ in this limit.
Therefore, letting $x=1-\sqrt{\delta}\xi$  and $y=\sqrt{\delta}R$ in (\ref{897}) and (\ref{898}) we find that ${\cal P}(x, y)$ agrees with  (\ref{895}).

We next examine the asymptotic behavior of (\ref{895}) as $R\rightarrow 0$ and $\xi\rightarrow \infty$ (corresponding to $u \rightarrow 0$ and $\omega\rightarrow \infty$ with $u\omega$ fixed).  From (\ref{887}) and (\ref{888}) we obtain
\begin{equation*}
u \sim  \frac{\varepsilon\xi}{2\omega}, \; \; \;  e^\omega \sim \frac{\xi}{2\omega R}, \; \; \; \omega \sim \log\left(\frac{\xi}{R}\right).
\end{equation*}
Then (\ref{84}) becomes
\begin{equation}\label{8as}
\delta^{3/2}\Omega(\xi, R) \sim \delta^{3/2}\left(\frac{a e^{-a}}{4\sqrt{2}}\right)\frac{ \xi^2}{R}|\log R|^{-5/2}\exp\left[-\frac{(\xi-R)^2}{4\log(\xi/R)}\right].
\end{equation}
If we fix $R$ and let $\xi \rightarrow \infty$ (corresponding to $u, \omega \rightarrow \infty$), similarly we obtain
\begin{equation}\label{8ae}
\pi(k, r) \sim \delta^{3/2}\left(\frac{a e^{-a}}{4\sqrt{2}}\right)\frac{ \xi^2}{R}(\log \xi)^{-5/2}\exp\left[-\frac{(\xi-R)^2}{4\log(\xi/R)}\right].
\end{equation}
A more uniform result can be obtained by replacing $\log \xi$ in (\ref{8ae}) by $\log(\xi/R)$, which will contain also (\ref{8as}) as a special case.

We next examine (\ref{881}) and (\ref{882}) as $\xi \downarrow 0,\;  R\rightarrow \infty$ with $\xi R$ fixed.  We shall thus obtain some of the higher order terms in (\ref{8aq}).  We set
$$R=\frac{R^*}{\varepsilon}, \; \; \xi = \varepsilon \xi^\diamond,$$
and expand $\Omega (\xi, R)$ as
$$\Omega(\xi, R) = ae^{-a}\left[\frac{1}{\varepsilon}\Omega^{(0)}(\xi^\diamond, R^*)+\varepsilon\Omega^{(1)}(\xi^\diamond, R^*)+\varepsilon^3 \Omega^{(2)}(\xi^\diamond, R^*)+O(\varepsilon^5)\right].$$
Then, from (\ref{881})  we obtain for the first three orders ($O(\varepsilon^{-2})$, $O(1)$, $O(\varepsilon^2)$)
\begin{eqnarray}
\Omega^{(0)}_{\xi^\diamond\xi^\diamond}+R^*\Omega^{(0)}_{\xi^\diamond} &=& 0, \label{8k}\\
\Omega^{(1)}_{\xi^\diamond\xi^\diamond}+R^*\Omega^{(1)}_{\xi^\diamond} &=& -\Omega^{(0)}-R^*\Omega^{(0)}_{R^*}, \\
\Omega^{(2)}_{\xi^\diamond\xi^\diamond}+R^*\Omega^{(2)}_{\xi^\diamond} &=& -\Omega^{(1)}-R^*\Omega^{(1)}_{R^*}, \label{8au}
\end{eqnarray}
and from (\ref{882}) we obtain the following boundary conditions at $\xi^\diamond = 0$
\begin{eqnarray}
\Omega^{(0)}_{\xi^\diamond}+R^*\Omega^{(0)} &=& 0, \label{8at}\\
\Omega^{(1)}_{\xi^\diamond}+R^*\Omega^{(1)} &=& \Omega^{(0)}_{R^*},\\
\Omega^{(2)}_{\xi^\diamond}+R^*\Omega^{(2)} &=& \Omega^{(1)}_{R^*}. \label{8l}
\end{eqnarray}
Solving (\ref{8k}) - (\ref{8l}) recursively we obtain after some calculation
\begin{eqnarray}
\Omega^{(0)} &=& R^*e^{-R^*\xi^\diamond}, \\
\Omega^{(1)} &=& \left[-\frac{R^*(\xi^\diamond)^2}{2} + \xi^\diamond +\frac{1}{R^*}\right] e^{-R^*\xi^\diamond},\\
\Omega^{(2)} &=& \left[\frac{R^*(\xi^\diamond)^4}{8}-\frac{(\xi^\diamond)^3}{6}-\frac{(\xi^\diamond)^2}{2R^*}-\frac{\xi^\diamond}{(R^*)^2} -\frac{6}{(R^*)^3}\right]e^{-R^*\xi^\diamond}.
\end{eqnarray}
We note that (\ref{8k}) and (\ref{8at}) imply that $\Omega^{(0)} = e^{-\xi^\diamondsuit R^*} f_0(R^*)$, and we ultimately determine $f_0$ from the solvability condition for (\ref{8au}).  Then to obtain $\Omega^{(1)}$ we must also consider the equation for $\Omega^{(3)}$, etc.
Therefore, as $\xi \downarrow 0$ and $R\rightarrow \infty$ with $\xi R=O(1)$,
\begin{multline}\label{8bb}
\Omega (\xi, R) = ae^{-a}e^{-\xi R}[\; R  +\left(\frac{1}{R} + \xi-\frac{R\xi^2}{2}\right) \\ +\left(-\frac{6}{R^3}-\frac{\xi}{R^2}-\frac{\xi^2}{2R}-\frac{\xi^3}{6}+\frac{R\xi^4}{8}\right)+O(R^{-5})\; ].
\end{multline}

Next we analyze (\ref{881}) and (\ref{882}) using a Laplace transform.  We set $\Omega(\xi, R) = {\cal D}(\xi, R)/R$.  Then we rewrite (\ref{881}) and (\ref{882}) as
\begin{eqnarray}
{\cal D}_{\xi\xi}+R({\cal D}_\xi+{\cal D}_R) &=& 0\; ; \; \; \xi, \; R>0,\label{85}\\
{\cal D}_\xi - {\cal D}_R+\left(\frac{1}{R}+R\right){\cal D} &=& 0\; ; \; \; \xi=0, \; \; R>0. \label{86}
\end{eqnarray}
We assume that ${\cal D} \rightarrow 0$ as $R\rightarrow 0$ and use a double Laplace transform, with
\begin{equation}
{\cal U}^*(R; \alpha) \equiv \int_0^\infty {\cal D}(\xi, R)e^{-\alpha \xi}\; d\xi,
\end{equation}
and
\begin{eqnarray}
{\cal U}(\alpha, \beta) &\equiv& \int_0^\infty  {\cal U}^*(R, \alpha)e^{-\beta R} \; dR, \nonumber\\
&=& \int_0^\infty \int_0^\infty {\cal D}(\xi, R)e^{-\alpha \xi}e^{-\beta R}\; d\xi \; dR.
\end{eqnarray}
Taking the Laplace transform of (\ref{85}) over $\xi$ gives
\begin{equation}\label{87}
\alpha^2 {\cal U}^*(R;\alpha)+\alpha R {\cal U}^*(R; \alpha)+R{\cal U}^*_R(R;\alpha)
=\alpha {\cal D}(0, R)+{\cal D}_\xi(0, R)+R{\cal D}(0, R).
\end{equation}
Taking the Laplace transform of (\ref{87}) over $R$ we obtain
\begin{equation}
\alpha^2 {\cal U}-(\alpha+\beta) {\cal U}_\beta -{\cal U} = \int_0^\infty e^{-\beta R}\left\{\alpha {\cal D}(0,R)+R\left[\frac{d}{dR}\left(\frac{{\cal D}(0,R)}{R}\right)\right]\right\}dR. \label{88}
\end{equation}
Here we also used (\ref{86}) to eliminate ${\cal D}_\xi(0, R)$.
Integrating by parts in the right-hand side and dividing both sides of (\ref{88}) by $-(\alpha+\beta)$ we obtain
\begin{equation}\label{89}
{\cal U}_\beta +\left(\frac{1-\alpha^2}{\alpha+\beta}\right){\cal U} = \int_0^\infty e^{-\beta R}\left(\frac{1}{R(\alpha+\beta)}-1\right) {\cal D}(0,R) dR,
\end{equation}
and we can rewrite (\ref{89}) as
\begin{equation}\label{8a}
\left[(\alpha+\beta)^{1-\alpha^2}{\cal U}\right]_\beta = \int_0^\infty \left(\frac{1}{R}-\alpha-\beta\right){\cal D}(0,R) e^{-\beta R}(\alpha+\beta)^{-\alpha^2} dR.
\end{equation}
We assume that the real part of $\alpha^2$ is negative, i.e. $\Re (\alpha^2)<0$, integrate  (\ref{8a}) over $(-\alpha, \beta)$, and obtain
\begin{multline}
{\cal U}(\alpha, \beta) = \frac{1}{\alpha+\beta}\int_{-\alpha}^\beta\int_0^\infty \left(\frac{\alpha+\beta}{\alpha+z}\right)^{\alpha^2}\left(\frac{1}{R}-\alpha-z\right){\cal D}(0, R) e^{-zR}\; dR \; dz \\
= \int_{-1}^0\int_0^\infty \left(\frac{1}{1+\theta}\right)^{-\alpha^2}\left[\frac{1}{\tilde{\chi}}-(\alpha+\beta)(1+\theta)\right]{\cal D}(0,\tilde{\chi}) e^{-\beta\tilde{\chi}}e^{-(\alpha+\beta)\tilde{\chi}\theta}\; d\tilde{\chi} \; d\theta. \label{8b}
\end{multline}
Here we changed variables with $z=(\alpha+\beta)\theta+\beta$ and $R=\tilde{\chi}$.  Inverting (\ref{8b}) over $\beta$ gives
\begin{eqnarray}
{\cal U}^*(R;\alpha) &=& \int_{-1}^0 \left(\frac{1}{1+\theta}\right)^{\alpha^2}e^{-\alpha\theta R/(1+\theta)}\nonumber\\
&\times& \left[\left(\frac{1}{R}-\alpha+\frac{\alpha \theta}{1+\theta}\right){\cal D}\left(0,\frac{R}{1+\theta}\right)-\frac{1}{1+\theta} {\cal D}_{\tilde{\chi}}\left(0, \frac{R}{1+\theta}\right)\right]d\theta \nonumber\\
&=& \int_0^1 \tilde{\theta}^{-\alpha^2} e^{-\alpha R (1-1/\tilde{\theta})}\left[\left(\frac{1}{R}-\frac{\alpha}{\tilde{\theta}}\right) {\cal D}\left(0, \frac{R}{\tilde{\theta}}\right)-\frac{1}{\tilde{\theta}}{\cal D}_{\tilde{\chi}}\left(0, \frac{R}{\tilde{\theta}}\right)\right] d\tilde{\theta}\nonumber\\
&=& \int_R^\infty \left(\frac{\chi}{R}\right)^{\alpha^2} e^{\alpha(\chi-R)}\left[\left(\frac{1}{\chi^2}-\frac{\alpha}{\chi}\right){\cal D}(0, \chi)-\frac{{\cal D}_{\chi}(0, \chi)}{\chi}\right]d\chi, \label{8c} \nonumber\\
\end{eqnarray}
where we set $\theta = \tilde{\theta}-1$ and $\tilde{\theta}=R/\chi$.

Inverting (\ref{8c}) over $\alpha$ we obtain
\begin{eqnarray}
{\cal D}(\xi, R) &=& \int_R^\infty \{\; \frac{1}{2\pi i}\int_{Br} \exp\left[\alpha (\chi-R+\xi)+\alpha^2\log\chi-\alpha^2\log R\right]\nonumber\\
&\times& \left[(1-\alpha\chi){\cal D}(0, \chi)-\chi{\cal D}_{\chi}(0, \chi)\right]\; d\alpha \} \; \frac{d\chi}{\chi^2} \nonumber\\
&=&  \frac{1}{2\sqrt{\pi}}\int_R^\infty \exp\left[-\frac{(\chi-R+\xi)^2}{4\log(\chi/R)}\right]\nonumber\\
&\times& \left[\left(\frac{1}{\chi}+\frac{\chi-R+\xi}{2\log(\chi/R)}\right)\frac{{\cal D}(0, \chi)}{\chi\sqrt{\log(\chi/R)}}-\frac{{\cal D}_\chi(0,\chi)}{\chi\sqrt{\log(\chi/R)}}\right] d\chi\nonumber \label{8d}\\
\\
&=& \frac{1}{2\sqrt{\pi}}\int_R^\infty \frac{(\chi-R+\xi)^2}{4[\log(\chi/R)]^{5/2}}\frac{{\cal D}(0, \chi)}{\chi^2} \exp\left[-\frac{(\chi-R+\xi)^2}{4\log(\chi/R)}\right] d\chi \nonumber\\
&-& \frac{1}{2\sqrt{\pi}}\int_R^\infty \frac{d}{d\chi}\left\{\frac{{\cal D}(0,\chi)}{\chi}\exp\left[-\frac{(\chi-R+\xi)^2}{4\log(\chi/R)}\right]\right\}\frac{d\chi}{\sqrt{\log(\chi/R)}}.\nonumber\\
\end{eqnarray}
Integrating the second integral by parts, we obtain an alternate expression for ${\cal D}(\xi, R)$ as
\begin{multline}\label{8be}
{\cal D}(\xi, R) = \frac{1}{2\sqrt{\pi}}\int_R^\infty \left\{\frac{(\chi-R+\xi)^2}{4[\log(\chi/R)]^{5/2}}-\frac{1}{2[\log(\chi/R)]^{3/2}}\right\}\frac{{\cal D}(0, \chi)}{\chi^2} \\
\times \exp\left[-\frac{(\chi-R+\xi)^2}{4\log(\chi/R)}\right] d\chi.
\end{multline}

Next we derive an integral equation for ${\cal D}(0, R)$ by
letting $\xi \rightarrow 0$ in (\ref{8d}). We can rewrite (\ref{8d}) as (\ref{8g}).
If $\xi=0$, the second integral has a non-integrable singularity at the boundary $\chi=R$, due to the factor $[\log(\chi/R)]^{-3/2}$.  If we set $\chi = Re^v$, the second integral becomes
\begin{multline}
\frac{1}{2\sqrt{\pi}}\int_0^\infty
\frac{\xi{\cal D}(0, Re^v)}{2v^{3/2}}\exp\left[-\frac{R^2(e^v-1+\xi/R)^2}{4v}\right]dv\\
= \frac{1}{4\sqrt{\pi}}\int_0^\infty \frac{{\cal D}(0, Re^{w\xi^2})}{w^{3/2}} \exp\left[-\frac{R^2(e^{w\xi^2}-1+\xi/R)^2}{4w\xi^2}\right] dw \label{8e}
\end{multline}
where $v=w\xi^2$.
Then letting $\xi\downarrow 0$ in (\ref{8e}) we obtain
\begin{equation}\label{8f}
\frac{{\cal D}(0, R)}{4\sqrt{\pi}}\int_0^\infty w^{-3/2} e^{-1/4w} dw = \frac{{\cal D}(0, R)}{2\sqrt{\pi}}\int_0^\infty u^{-1/2} e^{-u} du = \frac{{\cal D}(0, R)}{2}.
\end{equation}
Therefore, from (\ref{8g}) and (\ref{8f}) we obtain (\ref{8av}).
Since ${\cal D}(0,R) = R\Omega(0,R)$ we thus also have the following integral equation for $\Omega(0,R)$
\begin{multline}\label{8aw}
R\Omega(0, R) =\frac{1}{\sqrt{\pi}}\int_R^\infty \frac{1}{\sqrt{\log(\chi/R)}} \; \exp\left[-\frac{(\chi-R)^2}{4\log(\chi/R)}\right] \\
\times \left[\frac{\chi-R}{2\log(\chi/R)}\; \Omega(0, \chi)-\Omega_\chi(0,\chi)\right] d\chi.
\end{multline}

We can also analyze the integral equation (\ref{8aw}) asymptotically.  Consider the limit $R\rightarrow \infty$.  Then assuming that $\Omega(0, R)$ has mild (e.g. algebraic) growth, the function $\exp\{-(\chi-R)^2/[4\log(\chi/R)]\}\Omega(0, \chi)$ becomes sharply concentrated at $\chi =R$, which is the lower limit on the integral in (\ref{8aw}).  We thus set $\chi =R+{\cal V}/R$  in (\ref{8aw}) and note that
\begin{equation}
\frac{(\chi-R)^2}{\log(\chi/R)} = {\cal V}+\frac{{\cal V}^2}{2R^2}-\frac{{\cal V}^3}{12 R^4} +O(R^{-6}).
\end{equation}
Then if we evaluate (\ref{8aw}) by an implicit form of Watson's lemma, to leading order the right-hand side becomes $R\Omega(0,R)$, which is the same as the left-hand side, but this argument does not determine $\Omega(0,R)$.  However, by considering higher order terms in the expansion of the integral, we find after some calculation that the right-hand side of (\ref{8aw}) has the form, with $\Omega(0, R)=\Omega_0(R)$,
\begin{multline}\label{8ba}
R\Omega_0(R)+\frac{1}{R^3}\left[2R^2\Omega_0''(R)-4R\Omega_0'(R)+4\Omega_0(R)\right]\\
+\frac{1}{R^5}\left[8R^3\Omega_0'''(R)-36R^2\Omega_0''(R)+88R\Omega_0'(R)-100 \Omega_0(R)\right]+O(R^{-7}\Omega_0).
\end{multline}
We note that if $\Omega_0$ has algebraic behavior as $R\rightarrow \infty$ then $\Omega_0,\; R\Omega_0', \; R^2\Omega_0'', \; \dots $ all have the same order of magnitude.  The expansion in (\ref{8ba}) is in powers of $R^{-2}$, but the coefficient of $R^{-1}\Omega_0$ turned out to be zero. Setting (\ref{8ba}) equal to $R\Omega_0(R)$  we see that asymptotically the integral equation (\ref{8aw}) can be approximated by $2R^2\Omega_0''(R)-4R\Omega_0'(R)+4\Omega_0(R) =0.$  This is an ODE of Cauchy-Euler type that admits the solutions $\Omega_0(R)=R$ and $\Omega_0(R)=R^2$.  But by asymptotic matching to the ray expansion ($x=1$, $\; y>0$) we know that $\Omega_0(R) \sim ae^{-a}R$,  which precludes the second solution.  Then from (\ref{8ba}) we can conclude that if $\Omega_0(R)=ae^{-a}[R+cR^{-1}+O(R^{-3})]$ as $R\rightarrow \infty$, then $c=1$.  This result agrees precisely with (\ref{8bb}) with $\xi=0$.  The third ($-6R^{-3}$) term in the expansion of $\Omega_0(R)$ can also be obtained by explicitly evaluating the $O(R^{-7}\Omega_0)$  terms in (\ref{8ba}).

Now consider the limit $\xi, R \rightarrow \infty$ with $\xi/R$ fixed, and equation (\ref{8g}).  We evaluate the integral by an implicit form of the Laplace method.  Scaling $\chi = Rv$ the integrand will be maximal where
$$\frac{d}{dv}\left\{\frac{[R(v-1)+\xi]^2}{\log v}\right\}=0,$$
or $$ 2v\log v = v-1+\frac{\xi}{R}.$$
Letting $v=\sqrt{e}e^u$ this equation becomes equivalent to (\ref{899}), with $u=W((\xi/R-1)/(2\sqrt{e}))$ being the Lambert-$W$ function.  Then if we evaluate (\ref{8g}) by the Laplace method we obtain, after multiplying by $R$,
\begin{multline}\label{8bc}
R\Omega(\xi, R) \sim \exp \left[-R^2 w_* v_*^2\right]  \sqrt{\frac{2\pi}{\bigstar}} \\
\times\left\{ \frac{[R(v_*-1)+\xi]R\Omega_0(Rv_*)}{4\sqrt{\pi}(\log v_*)^{3/2}}
-\frac{1}{2\sqrt{\pi}\sqrt{\log v_*}}\left[\frac{d}{dv} \Omega_0(Rv)\right]\Bigl |_{v=v_*} \right\}
\end{multline}
where $v_*$ is the location of the maximum, $w_*=\log v_*$ and $\bigstar = H''(v_*)$ where $H(v) = [R(v-1)+\xi]^2/(4\log v).$  We see that the first term in the right-hand side of (\ref{8bc}) dominates the second, and after some calculation we find that $\bigstar = R^2(1+2w_*)/(2w_*)$ so that (\ref{8bc}) becomes
\begin{equation}\label{8bd}
\Omega(\xi, R) \sim \Omega_0(Rv_*) \frac{v_*}{\sqrt{1+2w_*}} \exp\left(-R^2 w_* v_*^2\right).
\end{equation}
Now by matching $\Omega_0(Rv_*) \sim ae^{-a}Rv_*$, so that (\ref{8bd}) agrees precisely with (\ref{895}), after we identify $w_* \leftrightarrow \omega$,  $\; v_* \leftrightarrow e^\omega$.

Thus (\ref{8be}) or (\ref{8g}) yields quite a bit of information as $R\rightarrow \infty$, since we can then localize the integral operator.  However, the opposite limit $R\rightarrow 0^+$ seems much more difficult as then the global nature of the operator persists, and the asymptotic evaluation of the integral seems to require considering separately the contributions from different ranges (such as $\chi \sim R$, $\; \chi =O(1)$ and $\chi =O(R)$).  If  $\chi = O(1)$ the contribution will of course involve ${\cal D}(0, \chi)$ (or $\Omega_0(\chi)$) for $\chi = O(1)$, which we do not have explicitly.

We can also derive a heat equation from (\ref{881}) and (\ref{882}).  We set $\xi-R=\zeta$ and $\Omega(\xi, R) = \tilde{{\cal D}}(\zeta, R)/R$.  Then (\ref{881}) and (\ref{882}) can be rewritten as
\begin{eqnarray}
\tilde{{\cal D}}_{\zeta\zeta}+R\tilde{{\cal D}}_R &=& 0\; ; \; \; R>0,  \;  \zeta+R>0, \label{8h}\\
2\tilde{{\cal D}}_\zeta-\tilde{{\cal D}}_R +\left(R+\frac{1}{R}\right)\tilde{{\cal D}} &=& 0\; ; \; \; R>0,\; \;  \zeta+R=0. \label{8i}
\end{eqnarray}
Changing variables from $R$ to $R=e^{-\Lambda}$ and letting $\tilde{{\cal D}}(\zeta, R)={\cal D}^*(\zeta, \Lambda)$, from (\ref{8h}) and (\ref{8i}), for $-\infty<\Lambda<\infty$,  we obtain
\begin{eqnarray}
{\cal D}^*_\Lambda = {\cal D}^*_{\zeta\zeta} &;& \zeta+e^{-\Lambda}>0, \label{8j}\\
2e^{-\Lambda}{\cal D}^*_\zeta+{\cal D}^*_\Lambda +\left(1+e^{-2\Lambda}\right){\cal D}^* = 0 &;&  \zeta+e^{-\Lambda}=0.\label{8ar}
\end{eqnarray}
We have thus reduced the problem to the heat equation, with $\Lambda$ corresponding to the time variable, on a domain with a moving boundary $\zeta=-e^{-\Lambda}$.  On this boundary the condition (\ref{8ar}) must be satisfied for all $\Lambda \in (-\infty, \infty)$.  We could convert this problem to an integral equation, but this would not seem to have any advantages over (\ref{8av}) or (\ref{8aw}).

\subsection{Analysis near the corner $(x, y)=(1, 0)\;$:\\ \; \;  $x=1-O(m^{-1/4}), \; y=O(m^{-1/2})$}
Next we examine the behavior of $\pi(k, r)$ in a region that is further away from the boundary of $x=1$, than the $\xi$ scale.
In this range the asymptotics of $\pi(k,r)$ can be obtained by expanding the ray solution ${\cal P}(x, y) \sim \delta {\cal K} \exp(\Psi/\delta)$ as $(x, y)\rightarrow (0,0)$.  We set
$$1-x=\delta^{1/4}\eta, \; \; y=R\sqrt{\delta}$$
and expand the ray solution ${\cal P}(x, y)$  as $s\rightarrow 0$ and $t\rightarrow \infty$ letting
\begin{eqnarray}
t &=& -\log s + T_*, \; \; \; T_* = T_0 + \delta^{1/4}T_1+\delta^{1/2}T_2 +O(\delta^{3/4}), \label{8m}\\
s &=& \delta^{1/4}S_*, \; \; \; S_* =  S_0 + \delta^{1/4}S_1+\delta^{1/2}S_2 +O(\delta^{3/4}),\label{8ay}\\
R &=& \bar{Y}(S_*, T_*; \delta) = \bar{Y}_0+ \delta^{1/4} \bar{Y}_1 + \delta^{1/2} \bar{Y}_2 +O(\delta^{3/4}), \label{8ax}\\
\eta &=& \bar{Z}(S_*, T_*; \delta) = \bar{Z}_0 +\delta^{1/4} \bar{Z}_1 + \delta^{1/2} \bar{Z}_2 +O(\delta^{3/4}). \label{8n}
\end{eqnarray}
On the ($\eta, R$) scale the ray expansion simplifies considerably, and with this simplified form we will be able to relate this scale to the ($\xi, R$) scale in subsection \textbf{5.4}, as well as the ($x, r$) scale in subsection \textbf{5.2}.
In (\ref{8ax}) and (\ref{8n}) it is understood that $S_*$ and $T_*$ are replaced by their expansions in powers of $\delta^{1/4}$, as in (\ref{8m}) and (\ref{8ay}).  Then we have $\bar{Y}_n = \bar{Y}_n(S_0, T_0)$, $\bar{Z}_n = \bar{Z}_n(S_0, T_0)$ and thus $\bar{Y}_0=R$, $\; \bar{Z}_0=\eta$ and $\; \bar{Y}_n = \bar{Z}_n = 0$ for $n>0$.  Using (\ref{8m}) - (\ref{8n}) in (\ref{Psi}) - (\ref{833}) we obtain, after a lengthy calculation
\begin{equation}\label{8o}
\pi(k, r) \sim  \frac{\delta a e^{-a} e^{T_0/2}}{\sqrt{1+e^{-T_0}}}\left(4e^{T_0}+2+\eta\sqrt{\frac{1+e^{-T_0}}{R}}\right)^{-1/2} \exp \left[ \frac{\Psi^{(0)}}{\sqrt{\delta}}+\frac{\Psi^{(1)}}{\delta^{1/4}}+\Psi^{(2)}\right]
\end{equation}
where
\begin{eqnarray}
\Psi^{(0)} &=& -\frac{\eta}{2}\sqrt{\frac{R}{1+e^{-T_0}}}-R\log \left(1+e^{T_0}\right)+\frac{1}{2}\left(\frac{R}{1+e^{-T_0}}\right), \label{8s}\\
\Psi^{(1)} &=& -\frac{\eta^2}{4}\sqrt{\frac{R}{1+e^{-T_0}}}-\frac{\eta}{2}\left(\frac{R}{1+e^{-T_0}}\right)+R\sqrt{\frac{R}{1+e^{-T_0}}}-\frac{3}{4}\left(\frac{R}{1+e^{-T_0}}\right)^{3/2}, \label{8z} \nonumber\\
\end{eqnarray}
and
\begin{multline}
\Psi^{(2)} = -\left(4e^{T_0}+2+\eta\sqrt{\frac{1+e^{-T_0}}{R}}\right)^{-1}[\; \frac{5}{32}\eta^4+\eta^3\left(\frac{e^{T_0}}{2}+\frac{3}{4}\right)\sqrt{\frac{R}{1+e^{-T_0}}}\\
+\eta^2\frac{R}{1+e^{-T_0}}\left(e^{T_0}+\frac{55}{48}-\frac{3}{4}e^{-T_0}\right)-\eta\left(\frac{R}{1+e^{-T_0}}\right)^{3/2}\left(\frac{2}{3}e^{T_0}+\frac{4}{3}+\frac{5}{2}e^{-T_0}\right)\\
-\left(\frac{R}{1+e^{-T_0}}\right)^2\left(e^{T_0}+\frac{55}{32}+\frac{9}{4}e^{-T_0}-\frac{1}{2}e^{-2T_0}\right) \; ]. \label{8aa}
\end{multline}
Here $T_0$ is defined implicitly by
\begin{equation}\label{8r}
\eta \sqrt{\frac{1+e^{-T_0}}{R}}= 2T_0-\frac{1}{2}\log\delta-\log R+\log\left(1+e^{-T_0}\right)-1.
\end{equation}
Since $\bar{Y}_0=R$ and $\; \bar{Z}_0=\eta\; $ imply that $S_0 = \sqrt{R/(1+e^{-T_0})}$ and $\eta= S_0[2T_0-2\log S_0-(\log\delta)/2-1]$ , we obtain (\ref{8r}).  Note that $T_0$ depends on $\eta$ and $R$, and also weakly upon $m$, due to the term $-\log \delta = \log m$.

To check that (\ref{8o}) satisfies (\ref{801}), we rewrite (\ref{801}) in terms of $\eta$ and $R$ as
\begin{equation}\label{8p}
\left(1-\delta^{1/4}\eta+\delta^{1/2}R\right)\left(2-e^{\Psi_\eta/\delta^{1/4}}\right)-\left(1-\delta^{1/4}\eta\right)e^{-\Psi_\eta/\delta^{1/4}}-\delta^{1/2}Re^{\Psi_R/\delta^{1/2}} = 0.
\end{equation}
Expanding $\Psi$ as
\begin{equation}\label{8q}
\Psi = \delta^{1/2} \Psi^{(0)} + \delta^{3/4}\Psi^{(1)} +\delta \Psi^{(2)} +O(\delta^{5/4}),
\end{equation}
and using this in (\ref{8p}), we obtain at the first three orders ($O(\delta^{1/2})$, $\;O(\delta^{3/4})$, $\; O(\delta)$)
\begin{eqnarray}
R\left[1-e^{\Psi^{(0)}_R}\right] &=& \left[\Psi^{(0)}_\eta\right]^2, \label{8x}\\
\eta \left[\Psi^{(0)}_\eta\right]^2-R \Psi^{(0)}_\eta &=& 2 \Psi^{(0)}_\eta \Psi^{(1)}_\eta + R e^{\Psi^{(0)}_R}\Psi^{(1)}_R,  \label{8y}
\end{eqnarray}
and
\begin{multline}\label{8ab}
R e^{\Psi^{(0)}_R}\Psi^{(2)}_R+2\Psi^{(0)}_\eta \Psi^{(2)}_\eta = -\left[\Psi^{(1)}_\eta\right]^2-R\Psi^{(1)}_\eta-\frac{R}{2}e^{\Psi^{(0)}_R}\left[\Psi^{(1)}_R\right]^2\\
+2\eta\Psi^{(0)}_\eta\Psi^{(1)}_\eta -\frac{R}{2}\left[\Psi^{(0)}_\eta\right]^2-\frac{1}{12}\left[\Psi^{(0)}_\eta\right]^4.
\end{multline}
Differentiating both sides of (\ref{8r}) with respect to $\eta$ gives
\begin{equation}\label{8t}
\frac{\partial T_0}{\partial\eta} = \left(\frac{\eta}{2}\frac{e^{-T_0}}{\sqrt{R(1+e^{-T_0})}}+1+\frac{1}{1+e^{-T_0}}\right)^{-1}\sqrt{\frac{1+e^{-T_0}}{R}}.
\end{equation}
From (\ref{8s}) we have
\begin{equation}
\Psi^{(0)}_\eta = -\frac{1}{2}\sqrt{\frac{R}{1+e^{-T_0}}}-\left[ \frac{\eta \sqrt{R} e^{-T_0}}{4 \left(1+e^{T_0}\right)^{3/2}}+\frac{R}{2\left(1+e^{-T_0}\right)}+\frac{R}{2\left(1+e^{-T_0}\right)^2}\right]\frac{\partial T_0}{\partial\eta},
\end{equation}
and thus, using (\ref{8t}), we obtain
\begin{equation}\label{8v}
\Psi^{(0)}_\eta = -\sqrt{\frac{R}{1+e^{-T_0}}}.
\end{equation}
Similarly, differentiating both sides of (\ref{8r}) with respect to $R$ gives
\begin{equation}\label{8u}
\frac{\partial T_0}{\partial R} = \left(\frac{\eta}{2}\frac{e^{-T_0}}{\sqrt{R(1+e^{-T_0})}}+1+\frac{1}{1+e^{-T_0}}\right)^{-1}\frac{1}{R}\left(1-\frac{\eta}{2}\sqrt{\frac{1+e^{-T_0}}{R}}\right).
\end{equation}
From (\ref{8s}) and (\ref{8u}) we obtain
\begin{equation}\label{8w}
\Psi^{(0)}_R = -\log\left(1+e^{T_0}\right).
\end{equation}
We can easily check that (\ref{8v}) and (\ref{8w}) satisfy (\ref{8x}).  In the same manner, we can check that (\ref{8s})-(\ref{8aa}) satisfy (\ref{8y}) and (\ref{8ab}).

We examine the asymptotic behavior of (\ref{8o}) as $R\rightarrow 0$ with a fixed $\eta >0$ (corresponding to $S_0 \rightarrow 0$ and $T_0\rightarrow \infty$), which should match to (\ref{861}) when it is expanded as $r\rightarrow \infty$.
Letting $T_0\rightarrow \infty$ (and $R\rightarrow 0$) in (\ref{8r}) we obtain
\begin{equation}
T_0 \sim \frac{1}{4}\log\delta+\frac{1}{2}\log R+\frac{1}{2}\left(\frac{\eta}{\sqrt{R}}+1\right),
\end{equation}
which when used in (\ref{8o}) gives
\begin{eqnarray}
\pi(k, r) &\sim& \frac{\delta a e^{-a}}{2}\; \exp [\; -\frac{1}{\sqrt{\delta}}\left(\frac{\eta\sqrt{R}}{2}+\frac{R}{4}\log\delta+\frac{R\log R}{2}+\frac{\eta \sqrt{R}}{2}\right)\nonumber\\
&-& \frac{1}{\delta^{1/4}}\left(\frac{\eta^2 \sqrt{R}}{4}+\frac{\eta R}{2}- R\sqrt{R}+\frac{3}{4}R^{3/2}\right)
-\frac{\eta^3 \sqrt{R}}{8}-\frac{\eta^2 R}{4} +\frac{\eta R^{3/2}}{6}+\frac{R^2}{4} \; ]\nonumber\\
&\sim& \frac{ a e^{-a}}{2}\; m^{r/2-1}r^{-r/2}\exp \left(-\eta\sqrt{r} m^{1/4}-\frac{\eta^2 \sqrt{r}}{4}\right). \label{8ac}
\end{eqnarray}
Expanding (\ref{861}) as $x\rightarrow 1$ and letting $x = 1-\eta m^{-1/4}$ also lead to (\ref{8ac}), which verifies the matching between the ($\eta, R$) and ($x, r$) scales.

Next we expand (\ref{8o}) and (\ref{8r}) as $\eta \rightarrow 0$ with a fixed $R$, by setting $\eta = \xi \delta^{-1/4}$, which corresponds to $T_0\rightarrow -\infty$ and $S_0\rightarrow 0$ with $S_0 T_0 << 1$.  We note that in this limit,  $\; \xi/R\rightarrow \infty$.  Then we obtain from (\ref{8o})
\begin{multline}\label{8ad}
\pi(k, r) \sim \delta a e^{-a} e^{T_0}\left(2+\xi\delta^{1/4}\frac{e^{-T_0/2}}{\sqrt{R}}\right)^{-1/2}
\exp [\;-\delta^{-1/2}e^{T_0}\frac{R}{2}-\delta^{-1/4}e^{T_0/2}\frac{\xi \sqrt{R}}{2}\\
+\delta^{-1/4}e^{T_0/2}R\sqrt{R}
- \frac{R^2}{2}\left(2+\delta^{1/4}e^{-T_0/2}\frac{\xi}{\sqrt{R}}\right)^{-1}\; ],
\end{multline}
and (\ref{8r}) becomes
\begin{equation}
T_0 \sim \delta^{1/4}e^{-T_0/2}\frac{\xi}{\sqrt{R}}+\frac{1}{2}\log\delta +\log R +1.
\end{equation}
We let $\nabla = \delta^{1/4}e^{-T_0/2}\xi/\sqrt{R}$ and using $T_0 \sim \nabla +(\log\delta)/2 + \log R +1$ obtain
\begin{equation}\label{8af}
\frac{\xi}{R} \sim \sqrt{e}\nabla e^{\nabla/2}.
\end{equation}
We also note that $e^{T_0} = \sqrt{\delta}\xi^2/(R\nabla^2)$.
Then we can rewrite (\ref{8ad}) as
\begin{equation}\label{8ag}
\pi(k, r) \sim \delta^{3/2} a e^{-a} \frac{\xi^2}{R}\frac{1}{\nabla^2\sqrt{2+\nabla}}
\exp \left(-\frac{\xi^2}{2\nabla^2}-\frac{\xi^2}{2\nabla}+\frac{\xi R}{\nabla}- \frac{R^2}{2(2+\nabla)}\right).
\end{equation}
Since $\xi/R \rightarrow \infty$, from (\ref{8af}) we obtain  $\nabla \sim 2\log(\xi/R)$  and (\ref{8ag}) thus agrees with (\ref{8ae}).
This verifies the matching between the ($\eta, R$) and ($\xi, R$) scales, for $\xi \rightarrow \infty$ and $\eta \rightarrow 0$.

\subsection{Analysis near the corner $(x, y)=(1, 0)\;$: $\; r=O(1)$}

We consider $r=O(1)$, which corresponds to there being only a few occupied secondary spaces, and $x \sim 1$.  We shall discuss the scales $k=m-O(1)$, $k=m-O(\sqrt{m})$, $\; k=m-O(\sqrt{m\log m})$ and $k=m-O(\sqrt{m}\log m) = m-\sqrt{m}(\log m)\bar{\xi}$.  To get an idea of the forms of the expansions for $\pi(k,r)$ on these scales, we first expand the results of subsection \textbf{5.5}, which apply on the $(\eta, R)$ scale, in the limit $\eta\rightarrow 0$, $\; R\rightarrow 0$, and then rewrite the results in terms of $\bar{\xi}$ and $r$.  We note that
\begin{equation}
1-x = \delta^{1/4}\eta = \delta^{1/2}\xi = -\delta^{1/2}(\log\delta)\; \bar{\xi}, \; \; R=\delta^{1/2} r,
\end{equation}
which relates the variables $x$, $\eta$, $\bar{\xi}$ and $\xi$.

We rewrite (\ref{8r}), which defines $T_0$, as
\begin{equation}\label{8ai}
\frac{\eta}{\sqrt{R}}=\frac{\xi}{\sqrt{r}} = -\frac{\bar{\xi}\log\delta}{\sqrt{r}} =  \frac{1}{\sqrt{1+e^{-T_0}}}\left[ 2T_0-\log\delta-\log r+\log\left(1+e^{-T_0}\right)-1\right].
\end{equation}
Thus, if $T_0 = O(1)$,
\begin{equation} \label{8ah}
\frac{\bar{\xi}}{\sqrt{r}} \sim \frac{1}{\sqrt{1+e^{-T_0}}}, \; \; \; T_0 \sim -\log\left(\frac{r}{\bar{\xi}^2}-1\right).
\end{equation}
There is a singularity in this approximation to $T_0$ when $r=\bar{\xi}^2$, and this will lead to a transition in the asymptotics along $r=\bar{\xi}^2$.  We also note that
\begin{center}
$T_0 = 0 \; \; $if $ \; \;  r=2\bar{\xi}^2, \; \; \;
T_0 < 0 \; \;$ if $ \; \;  r>2\bar{\xi}^2, \; \; \;
T_0 > 0 \; \;$ if $ \; \;  \bar{\xi}^2<r<2\bar{\xi}^2,$
\end{center}
and (\ref{8ai}) implies that if $r<\bar{\xi}^2$,  $T_0 \rightarrow \infty$ with
\begin{equation}\label{31}
T_0 \sim -\frac{\log\delta}{2}\left(\frac{\bar{\xi}}{\sqrt{r}}-1\right) = \frac{\log m}{2}\left(\frac{\bar{\xi}}{\sqrt{r}}-1\right).
\end{equation}

We first consider the case $r<\bar{\xi}^2$.  Letting $\eta = m^{-1/4}(\log m) \bar{\xi}$, $\; R=m^{-1/2}r$ and $T_0\rightarrow \infty$ in (\ref{8o}) yields
\begin{equation}\label{8aj}
\pi(k, r) \sim \delta \frac{a e^{-a}}{2} \exp \left[-\frac{\bar{\xi}\sqrt{r}}{2}\log m-r T_0+\frac{r}{2}\right].
\end{equation}
As $T_0\rightarrow \infty$ (\ref{8ai}) implies that
\begin{equation}
T_0 \sim \frac{\bar{\xi}\log m}{2\sqrt{r}} -\frac{\log m}{2}+\frac{\log r}{2}+\frac{1}{2},
\end{equation}
which refines the approximation in (\ref{31}).
Thus (\ref{8aj}) becomes
\begin{eqnarray} \label{8ak}
\pi(k, r) &\sim& \delta \frac{a e^{-a}}{2} \exp \left[-\bar{\xi}\sqrt{r}\log m +\frac{r}{2}\log m-\frac{r\log r}{2}\right] \nonumber\\
&=& \frac{a e^{-a}}{2} m^{r/2-\bar{\xi}\sqrt{r}-1}r^{-r/2}, \; \; \bar{\xi}>\sqrt{r}.
\end{eqnarray}
We note that the expression in (\ref{8ak}) is algebraically small in $m$, and becomes $O(m^{-r/2-1})$ when $\bar{\xi}=\sqrt{r}.$
It is easy to see that (\ref{8ak}) is precisely the expansion of (\ref{8ac}) as $\eta\rightarrow 0$ ($\eta = m^{-1/4}(\log m)\; \bar{\xi}$) and thus it also matches to the boundary layer solution (\ref{861}), when it is expanded as $x\rightarrow 1$.

When $r>\bar{\xi}^2$ we refine (\ref{8ah}) by setting
\begin{equation}\label{32}
 T_0 = -\log\left(\frac{r}{\bar{\xi}^2}-1\right)+\frac{\clubsuit}{\log m}.
\end{equation}
Using (\ref{32}) in (\ref{8ai}) we obtain
\begin{equation}
\clubsuit \sim \frac{4\log(r/\bar{\xi}-\bar{\xi})+2}{1-\bar{\xi}^2/r}
\end{equation}
and thus
\begin{equation}\label{8al}
T_0 \sim -\log\left(\frac{r}{\bar{\xi}^2}-1\right)+\frac{4\log(r/\bar{\xi}-\bar{\xi})+2}{(1-\bar{\xi}^2/r)\log m}.
 \end{equation}
Using (\ref{8al}) in (\ref{8o}) leads to
\begin{equation}\label{8ao}
\pi(k, r) \sim ae^{-a}\frac{m^{-\bar{\xi}^2/2}}{m\sqrt{\log m}}\frac{\bar{\xi}^2}{\sqrt{r(r-\bar{\xi}^2)}}\left(\frac{\bar{\xi}}{r-\bar{\xi}^2}\right)^{\bar{\xi}^2}\left(\frac{r-\bar{\xi}^2}{r}\right)^r, \; \; \bar{\xi}<\sqrt{r}.
\end{equation}

Next we examine the behavior of (\ref{8o}) in the transition region where $\bar{\xi} \sim \sqrt{r}$. We set
\begin{equation}
T_0 = \log\log m + \log \diamondsuit
\end{equation}
and
\begin{equation}\label{8bm}
\bar{\xi} = \sqrt{r}+\frac{\spadesuit}{\log m}
\end{equation}
for $\diamondsuit, \; \spadesuit  \ll \log m$.
Then we obtain from (\ref{8ai})
\begin{equation}
\spadesuit \sim 2\sqrt{r} \log\log m,
\end{equation}
and then defining $\bar{\xi}^*$ by
\begin{equation}\label{8an}
\bar{\xi} = \sqrt{r}+2\sqrt{r}\left(\frac{\log\log m}{\log m}\right)+\frac{\bar{\xi} ^*}{\log m}
\end{equation}
we find that (\ref{8ai}) becomes
\begin{equation}\label{33}
\frac{\bar{\xi}^*}{\sqrt{r}}=2\log \diamondsuit-\frac{1}{2\diamondsuit}-\log r-1+o(1).
\end{equation}
This defines $\diamondsuit$ implicitly in terms of $\bar{\xi}^*$ and $r$.  Then on the ($\bar{\xi}^*, r$) scale we obtain from (\ref{8o})
\begin{equation}\label{8am}
\pi(k, r) \sim ae^{-a}\sqrt{\frac{\diamondsuit}{4\diamondsuit+1}} m^{-r/2-1}(\log m)^{-2r} r^{-r/2}e^{-\bar{\xi} ^*\sqrt{r}}.
\end{equation}

As $\bar{\xi} ^* \rightarrow \infty$  (\ref{8am}) should match to (\ref{8ak}).  From (\ref{33}) we see that $\diamondsuit\rightarrow \infty$ and then using (\ref{8an}), (\ref{8am}) becomes the same as (\ref{8ak}).
We next check the matching between (\ref{8am}) and (\ref{8ao}).  As $\bar{\xi} ^* \rightarrow -\infty$, $\; \diamondsuit\rightarrow 0$ and (\ref{8an}) becomes
\begin{equation}\label{34}
\bar{\xi} ^2 \sim r\left(1+\frac{2\bar{\xi} ^*}{\sqrt{r}\log m}+\frac{4\log\log m}{\log m}\right).
\end{equation}
Using (\ref{34}) in (\ref{8ao}) we obtain
\begin{equation}\label{8az}
\pi(k, r) \sim \frac{ae^{-a}}{m}\left(-\frac{\sqrt{r}}{2\bar{\xi} ^*}\right)^{1/2}
\exp\left[-\frac{r}{2}\log m-2r\log\log m -\bar{\xi} ^*\sqrt{r} -\frac{r}{2}\log r\right].
\end{equation}
Here we used $-\bar{\xi} ^* \gg \log\log m$.
The exponential part of (\ref{8az}) is the same as the exponential part of (\ref{8am}).  From (\ref{33}) we obtain as $\diamondsuit\rightarrow 0$
\begin{equation}
-\frac{\bar{\xi} ^*}{\sqrt{r}} \sim \frac{1}{2\diamondsuit},
\end{equation}
and thus
\begin{equation}
\sqrt{\frac{\diamondsuit}{4\diamondsuit+1}} \sim \sqrt{\diamondsuit} \sim \left(-\frac{\sqrt{r}}{2\bar{\xi} ^*}\right)^{1/2}.
\end{equation}
Therefore, as $\bar{\xi} ^* \rightarrow -\infty$ and $\diamondsuit\rightarrow 0$ (\ref{8am}) matches to (\ref{8ao}).

Next we examine  $\pi(k, r)$ on the ($\bar{\xi}, r$) scale by analyzing
the balance equation (\ref{bl2}).  We set
\begin{equation}\label{35}
\pi \sim \frac{1}{m \sqrt{\log m}}\bar{K}(\bar{\xi}, r)\exp[f(\bar{\xi})\log m]
\end{equation}
and rewrite
(\ref{bl2}) as
\begin{multline}\label{8bf}
\left(2-\frac{a}{m}\right)\; \bar{K}(\bar{\xi},
r)e^{f(\bar{\xi})\log m} \\= \left(1-\frac{a}{m}\right) \;
\bar{K}\left(\bar{\xi}+\frac{1}{\sqrt{m}\log m},
r\right)e^{f(\bar{\xi}+1/(\sqrt{m}\log m))\log
m}\\
+\frac{m-\bar{\xi}\sqrt{m}\log m +1}{m-\bar{\xi}\sqrt{m}\log
m+r+1} \bar{K}\left(\bar{\xi}-\frac{1}{\sqrt{m}\log m},
r\right)e^{f(\bar{\xi}-1/(\sqrt{m}\log m))\log m} \\
 +\frac{r+1}{m-\bar{\xi}\sqrt{m}\log m +r+1}\;
\bar{K}(\bar{\xi}, r+1)e^{f(\bar{\xi})\log m}.
\end{multline}
From (\ref{8bf}) to leading order ($O(1/m)$) we obtain the limiting equation
\begin{equation}\label{8bg}
-\bar{K}(\bar{\xi}, r) [f'({\bar{\xi}})]^2+r \bar{K}(\bar{\xi}, r)-(r+1)\bar{K}(\bar{\xi}, r+1) = 0.
\end{equation}
By matching (\ref{35}) to (\ref{8ao}), when $r>\bar{\xi} ^2$ we must have $f(\bar{\xi}) = -\bar{\xi}^2/2$.  Using this in (\ref{8bg}) we obtain
\begin{equation}\label{8bh}
\frac{\bar{K}(\bar{\xi},r+1 )}{\bar{K}(\bar{\xi},r)} =\frac{r- \bar{\xi} ^2}{r+1}.
\end{equation}
Solving (\ref{8bh}) gives
\begin{equation}
\bar{K}(\bar{\xi}, r) = \left[\prod_{j=1}^{r-1} (j-\bar{\xi}^2)\right]\frac{\bar{K}(\bar{\xi},1)}{r!} = \frac{\Gamma(r-\bar{\xi}^2)}{\Gamma(1-\bar{\xi}^2)}\frac{\bar{K}(\bar{\xi},1)}{r!},
\end{equation}
and thus
\begin{equation}\label{8bi}
\pi(k,r) \sim \frac{m^{-\bar{\xi}^2/2}}{m\sqrt{\log m}}
\frac{\Gamma(r-\bar{\xi}^2)}{r!}\frac{\bar{K}(\bar{\xi},1)}{\Gamma(1-\bar{\xi}^2)}.
\end{equation}

We rewrite (\ref{8bi}) as
\begin{equation}\label{8bn}
\pi(k,r) \sim \frac{m^{-\bar{\xi}^2/2}}{m\sqrt{\log m}}
\frac{\Gamma(r-\bar{\xi}^2)}{r!}\; \bar{H}(\bar{\xi}), \; \; \bar{H}(\bar{\xi})=\frac{\bar{K}(\bar{\xi},1)}{\Gamma(1-\bar{\xi}^2)}.
\end{equation}
We see that as $\bar{\xi} \uparrow \sqrt{r}$, the factor $\Gamma(r-\bar{\xi}^2)$ becomes singular, which indicates a transition in the asymptotics.  Analysis of the range $\bar{\xi}\sim \sqrt{r}$ will also determine $\bar{H}(\bar{\xi})$ in (\ref{8bn}), and thus $\bar{K}(\bar{\xi}, 1)$ in (\ref{8bi}).  To study the transition we scale $\bar{\xi}-\sqrt{r} = O\left[(\log m)^{-1/2}\right]$ with
$$\bar{\xi} = \sqrt{r}+\frac{\circledast}{\sqrt{\log m}}, \;   \;    \; \pi(k,r) = {\cal H}_r(\circledast; m).$$
Then (\ref{bl2}) becomes
\begin{multline}
\left(2-\frac{a}{m}\right)\; {\cal H}_r(\circledast; m) =
\left[1-\frac{r}{m}+O\left(\frac{\log m}{m\sqrt{m}}\right)\right]\; {\cal H}_r\left(\circledast-\frac{1}{\sqrt{m\log m}}; m\right)\\
 +\left[\frac{r+1}{m}+O\left(\frac{\log m}{m\sqrt{m}}\right)\right]\; {\cal H}_{r+1}(\circledast; m)+\left(1-\frac{a}{m}\right) {\cal H}_r\left(\circledast+\frac{1}{\sqrt{m\log m}}; m\right),
\end{multline}
which we further approximate by
\begin{equation}\label{8bp}
(r+1){\cal H}_{r+1}(\circledast; m)-r{\cal H}_r(\circledast; m)+\frac{1}{\log m}{\cal H}_r''(\circledast; m)+o\left(\frac{1}{\log m}\right) = 0.
\end{equation}
Before analyzing (\ref{8bp}) we derive matching conditions for ${\cal H}_r(\circledast; m)$ as $\circledast\rightarrow \pm \infty$.
On the $(\bar{\xi}, r)$ scale for $\bar{\xi}>\sqrt{r}$, (\ref{8ak}) applies.  This can also be obtained by analyzing (\ref{bl2}) on the ($\bar{\xi}, r$) scale for $\bar{\xi}>\sqrt{r}$ and using asymptotic matching.  In either case we ultimately conclude that (\ref{861}) remains valid for $x\rightarrow 1$ as long as $\bar{\xi} = (1-x)\sqrt{m}/\log m > \sqrt{r}.$  By setting $\bar{\xi} = \sqrt{r}+\circledast/\sqrt{\log m}$ in (\ref{8ak}) we obtain
$$\frac{a e^{-a}}{2} m^{-r/2-1}r^{-r/2}e^{-\circledast\sqrt{r\log m}}\left[1+o(1)\right]$$
so that as $\circledast \rightarrow +\infty$ the matching condition for ${\cal H}_r(\circledast; m)$ is
\begin{equation}\label{8bs}
{\cal H}_r(\circledast; m) \sim \frac{a e^{-a}}{2} m^{-r/2-1}r^{-r/2}e^{-\circledast\sqrt{r\log m}}, \; \; \circledast \rightarrow +\infty.
\end{equation}
By examining (\ref{8bn}) as $\bar{\xi}\rightarrow \sqrt{r}$ and using
\begin{equation}
\Gamma(r-\bar{\xi}^2) \sim \frac{1}{\sqrt{r}-\bar{\xi}}\cdot \frac{1}{\sqrt{r}+\bar{\xi}} = -\frac{\sqrt{\log m}}{\circledast}\left(2\sqrt{r}+\frac{\circledast}{\sqrt{\log m}}\right)^{-1} \sim -\frac{\sqrt{\log m}}{2\circledast\sqrt{r}}
\end{equation}
we obtain from (\ref{8bn})
\begin{equation}\label{8bo}
{\cal H}_r(\circledast; m) \sim m^{-r/2-1}e^{-\circledast\sqrt{r\log m}}
\frac{e^{-\circledast^2/2}}{-2\circledast}\frac{\bar{H}(\sqrt{r})}{r!\sqrt{r}}, \; \; \circledast \rightarrow -\infty.
\end{equation}
In view of the matching conditions we set
\begin{equation}\label{8bt}
{\cal H}_r(\circledast; m) = m^{-r/2-1}e^{-\circledast\sqrt{r\log m}} [\hbar_r(\circledast)+o(1)].
\end{equation}
With (\ref{8bt}) we see that
\begin{equation}\label{8bj}
{\cal H}_r''(\circledast; m)-(\log m)r{\cal H}_r'(\circledast;m) \sim m^{-r/2-1}e^{-\circledast\sqrt{r\log m}}(-2)\sqrt{r\log m}\; \hbar_r(\circledast).
\end{equation}
Also, since $\bar{\xi}\sim \sqrt{r}$ we have $\bar{\xi}<\sqrt{r+1}$ and thus (\ref{8bn}) may be used to approximate ${\cal H}_{r+1}(\circledast; m)$ in (\ref{8bp}), as
\begin{equation}\label{8bq}
{\cal H}_{r+1}(\circledast; m) \sim \frac{m^{-r/2-1}}{\sqrt{\log m}}e^{-\circledast\sqrt{r\log m}}e^{-\circledast^2/2}
\frac{\bar{H}(\sqrt{r})}{(r+1)!}.
\end{equation}
Using (\ref{8bj}) and (\ref{8bq}) in (\ref{8bp}) we obtain to leading order
\begin{equation}\label{8bk}
2\sqrt{r}\hbar_r '(\circledast) = \frac{e^{-\circledast^2/2}}{r!} \bar{H}(\sqrt{r})
\end{equation}
The solution to (\ref{8bk}) that decays as $\circledast \rightarrow -\infty$ is
\begin{equation}\label{8br}
\hbar_r (\circledast) = \frac{\bar{H}(\sqrt{r})}{2r!\sqrt{r}}\int_{-\infty}^\circledast e^{-u^2/2}du.
\end{equation}

To determine $\bar{H}(\sqrt{r})$ we use the matching condition in (\ref{8bs}).
Letting $\circledast \rightarrow \infty$, we evaluate the integral in (\ref{8br}) as $\sqrt{2\pi}$ and rewrite (\ref{8bt}), for $\circledast\rightarrow \infty$, as
\begin{equation}\label{8bu}
{\cal H}_r(\circledast; m) \sim m^{-r/2-1}e^{-\circledast\sqrt{r\log m}} \frac{\bar{H}(\sqrt{r})}{r!\sqrt{r}}\sqrt{\frac{\pi}{2}}.
\end{equation}
Comparing (\ref{8bs}) to (\ref{8bu}) we obtain
\begin{equation}
\bar{H}(\sqrt{r})=\frac{a e^{-a}}{\sqrt{2\pi}} r^{-r/2}\sqrt{r}\; \Gamma(r+1),
\end{equation}
so that
\begin{equation}\label{8bv}
\bar{H}(z) = \frac{a e^{-a}}{\sqrt{2\pi}} z^{-z^2}z\Gamma(1+z^2).
\end{equation}
Thus on the ($\circledast, r$) transition scale we have
\begin{equation}
\pi(k,r) = {\cal H}_r(\circledast; m) \sim \frac{a e^{-a}}{2\sqrt{2\pi}}m^{-r/2-1}r^{-r/2}e^{-\circledast\sqrt{r\log m}} \int_{-\infty}^\circledast e^{-u^2/2}du,  \; \; r\geq 1,
\end{equation}
and when $\bar{\xi}<\sqrt{r}$, from (\ref{8bn}) and (\ref{8bv}), we obtain
\begin{equation}\label{8bw}
\pi(k,r) \sim \frac{a e^{-a}}{\sqrt{2\pi}}\frac{m^{-\bar{\xi}^2/2}}{m\sqrt{\log m}}
\frac{\Gamma(1+\bar{\xi}^2)\Gamma(r-\bar{\xi}^2)}{r!}\;  \bar{\xi}^{-\bar{\xi}^2}\bar{\xi}, \;  \;  r\geq 1.
\end{equation}
We note that letting $\bar{\xi}, r \rightarrow \infty$ in (\ref{8bw}) and using Stirling's formula for the Gamma functions lead to (\ref{8ao}).
This verifies the matching between the ($\bar{\xi},r$) scale and the ray expansion on the ($x, y$) scale.

We next examine $r=0$.  For $\bar{\xi}>1$ we can obtain the expansion of $\pi(k,0)$ as a limiting case of (\ref{876}), which leads to
\begin{equation}\label{37}
\pi(k,0)-(1-\rho)\rho^k \sim -\frac{ae^{-a}}{2}m^{-1/2-\bar{\xi}}, \; \; \bar{\xi}>1.
\end{equation}
For $0<\bar{\xi}<1$ we find that (\ref{8bg}) still holds at $r=0$, if we write
$$\pi(k,0)-(1-\rho)\rho^k \sim \frac{m^{-\bar{\xi}^2/2}}{m\sqrt{\log m}} \bar{K}(\bar{\xi}, 0).$$  Hence,
\begin{equation}\label{36}
\pi(k,0)-(1-\rho)\rho^k \sim \frac{a e^{-a}}{\sqrt{2\pi}}\frac{m^{-\bar{\xi}^2/2}}{m\sqrt{\log m}}
\bar{\xi}^{1-\bar{\xi}^2}\Gamma(1+\bar{\xi}^2)\Gamma(-\bar{\xi}^2), \; \; 0<\bar{\xi}<1.
\end{equation}
We note that $\pi(k,0)-(1-\rho)\rho^k$ is negative, since $\Gamma(-\bar{\xi}^2)<0$.  There is a transition in the asymptotics when $\bar{\xi}=1$, and (\ref{36}) also breaks down as $\bar{\xi}\rightarrow 0$, due to the singularities of $\Gamma(-\bar{\xi}^2)$.  We also note that on the $\bar{\xi}$ scale
\begin{equation}
(1-\rho)\rho^k = \frac{ae^{-a}}{m}\left[1+\frac{a\log m}{\sqrt{m}}\bar{\xi}+O\left(\frac{\log^2 m}{m}\right)\right] =  O(m^{-1}),
\end{equation}
so that the geometric part of $\pi(k, 0)$ is still asymptotically dominant for $\bar{\xi}>0.$

To study the transition range $\bar{\xi}\sim 1$ we let $\pi(k,0)-(1-\rho)\rho^k \sim {\cal H}_0(\circledast; m)$ where now $\circledast=\sqrt{\log m}(\bar{\xi}-1)$.  Note that the transition point is the same for $\pi(k,0)$ and $\pi(k,1)$.  By examining the matching condition to (\ref{37}) as $\circledast\rightarrow +\infty$ and to (\ref{36}) as $\circledast\rightarrow -\infty$, we again set
\begin{equation}\label{38}
{\cal H}_0(\circledast; m) \sim m^{-3/2}e^{-\circledast\sqrt{\log m}}\hbar_0(\circledast)
\end{equation}
and from (\ref{8bp}) with $r=0$ we find that $(\log m)^{-1}{\cal H}_0''(\circledast; m)\sim -{\cal H}_1(\circledast; m)$ so that $\hbar_0(\circledast)=-\hbar_1(\circledast)$ and hence
\begin{equation}
\pi(k,0) - (1-\rho)\rho^k \sim -\frac{ae^{-a}}{2\sqrt{2\pi}}m^{-3/2}e^{-\circledast\sqrt{\log m}}\int_{-\infty}^\circledast e^{-u^2/2} du, \; \; \circledast = \sqrt{\log m}(\bar{\xi}-1).
\end{equation}
We note that on the $\bar{\xi}$ scale the identity $\sum_{r=0}^m \pi(N-r, r) = (1-\rho)\rho^N$ is asymptotically satisfied for $\bar{\xi}<1$, in view of the fact that
$$\sum_{r=0}^\infty \frac{\Gamma(r-\bar{\xi}^2)}{r!} = 0,  \; \; 0<\bar{\xi}<1.$$

Now we consider scales that have $\bar{\xi}\rightarrow 0$, so that $m-k=o(\sqrt{m}\log m)$.  If we define $\tilde{\xi}$ by
$$m-k = \sqrt{m\log m}\; \tilde{\xi}, \; \; \; \tilde{\xi} = \sqrt{\log m}\; \bar{\xi}$$
then the factor $m^{-\bar{\xi}^2/2} = \exp(-\bar{\xi}^2\log m/2) = e^{-\tilde{\xi}^2/2}$  becomes $O(1)$ and (\ref{8bn}) simplifies to
\begin{equation}\label{8cd}
\pi(k,r) \sim \frac{1}{m\log m}\frac{a e^{-a}}{\sqrt{2\pi}}\; \frac{\tilde{\xi}}{r}\; e^{-\tilde{\xi}^2/2}, \; \; r\geq 1.
\end{equation}
Here we used $\Gamma(r-\bar{\xi}^2)\sim (r-1)!$ and $\bar{H}(\bar{\xi}) \sim ae^{-a}\bar{\xi}/\sqrt{2\pi}$ as $\bar{\xi}\rightarrow 0$.  We have verified, after a lengthy calculation which we omit, that (\ref{8cd}) can also be obtained by expanding $m^{-3/2}\Omega(\xi, R)$ for $\xi\rightarrow \infty$, $\; R\rightarrow 0$ with $\xi=O(\sqrt{-\log R})$, using the integral in (\ref{8g}) and our knowledge of ${\cal D}(0, R)$ (or $\Omega(0, R)$) as $R\rightarrow \infty$.  If we expand $\pi(k,0)$ on the $\bar{\xi}$ scale for $\bar{\xi}\rightarrow 0$, we obtain from (\ref{36})
\begin{equation}\label{302}
\pi(k,0)\sim ae^{-a}\left[\frac{1}{m}-\frac{1}{m}\frac{e^{-\tilde{\xi}^2/2}}{\sqrt{2\pi}\tilde{\xi}}\right]
\end{equation}
so that on the $\tilde{\xi}$ scale the geometric part of $\pi(k,0)$ no longer dominates the right side of (\ref{36}).  We thus re-examine the balance equation (\ref{bl2}) along $r=0$, i.e.,
\begin{equation}\label{39}
\left(2-\frac{a}{m}\right)\pi(k, 0) = \left(1-\frac{a}{m}\right)\pi(k-1, 0)+\pi(k+1, 0)+\frac{1}{k+1}\pi(k,1).
\end{equation}
On the $\tilde{\xi}$ scale we let $\pi(k,0)\sim m^{-1}\heartsuit(\tilde{\xi})$ and use (\ref{8cd}) to approximate $\pi(k,1)$.  Then (\ref{39}) becomes asymptotically (after we multiply by $m^2\log m$)
\begin{equation}\label{300}
-\heartsuit''(\tilde{\xi}) \sim \frac{ae^{-a}}{\sqrt{2\pi}}\; \tilde{\xi} e^{-\tilde{\xi}^2/2}.
\end{equation}
The solution of (\ref{300}) that matches, as $\tilde{\xi}\rightarrow \infty$, to (\ref{36}) is given by
\begin{equation}
\heartsuit(\tilde{\xi}) = ae^{-a}\left[1-\frac{1}{\sqrt{2\pi}}\int_{\tilde{\xi}}^\infty e^{-u^2/2} du\right]
\end{equation}
so that on the $\tilde{\xi}$ scale with $r=0$ we have
\begin{equation}\label{301}
\pi(k, 0)\sim \frac{ae^{-a}}{m}\left[1-\frac{1}{\sqrt{2\pi}}\int_{\tilde{\xi}}^\infty e^{-u^2/2} du\right].
\end{equation}
For $\tilde{\xi}\rightarrow \infty$  (\ref{301}) reduces to (\ref{302}).  This shows also that on the $\tilde{\xi}$ scale ($k=m-O(\sqrt{m\log m})$) $\pi(k,r) \; (r\geq 1)$ is smaller than $\pi(k, 0)$ by a factor of $(\log m)^{-1}$.

Next we discuss the problem on the ($n, r$) scale where $k=m-n$ and $n=O(1)$. Setting $\pi(k,r) = {\cal L}(n,r;\; m)$  we write (\ref{bl2}) and (\ref{bcm}) as
\begin{multline}\label{8by}
\left(2-\frac{a}{m}\right){\cal L}(n,r;\; m) = \left(1-\frac{a}{m}\right){\cal L}(n+1,r;\; m) + \frac{m-n+1}{m+r-n+1}{\cal L}(n-1,r;\; m)\\
+\frac{r+1}{m+r-n+1}{\cal L}(n,r+1;\; m), \; \; n\geq1
\end{multline}
and
\begin{multline}\label{8bz}
\left(2-\frac{a}{m}\right){\cal  L}(0,r;\; m) = \left(1-\frac{a}{m}\right){\cal L}(1,
r;\; m)+\frac{r+1}{m+r+1}{\cal L}(0, r+1; m)\\+\left(1-\frac{a}{m}\right){\cal L}(0, r-1;\; m), \; \; r\geq1.
\end{multline}
The corner condition (\ref{27}) becomes
\begin{equation}\label{304}
\left(2-\frac{1}{m}\right){\cal L}(0,0;\; m) = \left(1-\frac{1}{m}\right){\cal L}(1,0;\; m)+\frac{1}{m+1}{\cal L}(0,1;\; m).
\end{equation}
We expand ${\cal L}$ as
\begin{equation}\label{303}
{\cal L}(n, r;\; m)= \frac{1}{m^{3/2}\sqrt{\log m}}\left[{\cal L}(n,r)+ \frac{1}{\log m}{\cal L}^{(1)}(n,r)+O(\log^{-2}m)\right].
\end{equation}
The scale factor $m^{-3/2}(\log m)^{-1/2}$ must be included in view of matching considerations, which we discuss shortly.  Using (\ref{303}) in
(\ref{8by})-(\ref{304}) we find that the leading term ${\cal L}(n,r)$ satisfies
\begin{eqnarray}
2{\cal L}(n,r) - {\cal L}(n+1, r)-{\cal L}(n-1, r) &=& 0, \; \; r\geq 0, \label{8ca}\\
2{\cal L}(0,r) -{\cal L}(1,r)-{\cal L}(0,r-1) &=& 0, \; \; r\geq1,\label{8cb}
\end{eqnarray}
with the corner condition
\begin{equation}\label{8cc}
2{\cal L}(0,0) - {\cal L}(1,0) = 0.
\end{equation}
The correction term ${\cal L}^{(1)}(n,r)$ also satisfies (\ref{8ca})-(\ref{8cc}).

The general solution to (\ref{8ca}) is ${\cal L}(n,r) = {\cal V}(r)+n{\cal W}(r)$ and we then obtain from (\ref{8cb}) and (\ref{8cc})
\begin{equation}\label{8ce}
{\cal V}(r)-{\cal V}(r-1) = {\cal W}(r), \; \; r\geq1,
\end{equation}
and
\begin{equation}
{\cal V}(0)={\cal W}(0),
\end{equation}
which implies that ${\cal L}(n,0) = (n+1){\cal V}(0)$.
A similar argument can be used to determine ${\cal L}^{(1)}(n, r)$ as ${\cal V}^{(1)}(r)+n{\cal W}^{(1)}(r)$ so we write the expansion on the $(n,r)$ scale as
\begin{multline}\label{305}
{\cal L}(n,r;\; m) \sim \frac{1}{m^{3/2}\sqrt{\log m}}\; \{\; {\cal V}(0)+\sum_{l=1}^r {\cal W}(l)+n{\cal W}(r) \\
+\frac{1}{\log m}\left[{\cal V}^{(1)}(0)+\sum_{l=1}^r {\cal W}^{(1)}(l)+n{\cal W}^{(1)}(r)\right]\; \}
\end{multline}
for $r\geq1$, and
\begin{equation}\label{306}
{\cal L}(n,0;\; m) \sim  \frac{n+1}{m^{3/2}\sqrt{\log m}}\left[{\cal V}(0)+\frac{1}{\log m}{\cal V}^{(1)}(0)\right].
\end{equation}

Now we try to match (\ref{305}) and (\ref{306}) to (\ref{8cd}) and (\ref{301}), noting that $n=\tilde{\xi}\sqrt{m\log m}$.  On the $\tilde{\xi}$ scale ${\cal L}(n,0)$ in (\ref{306}) becomes $O(m^{-1})$ which is of the same order as (\ref{301}), but the expansions cannot match since (\ref{306}) will be linear in $n$ (or $\tilde{\xi}$), while (\ref{301}) does not vanish as $\tilde{\xi}\rightarrow 0 $.  Problems also arise in matching the $n$ and $\tilde{\xi}$ scales for $r\geq1$.  For a fixed $r$ and $n\rightarrow \infty$,  the leading term in (\ref{305}) becomes $m^{-3/2}(\log m)^{-1/2}n{\cal W}(r) = m^{-1}\tilde{\xi}{\cal W}(r)$ which is larger than (\ref{8cd}) by a factor of $\log m$.  We must thus set ${\cal W}(r)=0$.  We can match (\ref{8cd}) to the correction term in (\ref{305}) by setting ${\cal W}^{(1)}(r)=ae^{-a}/(\sqrt{2\pi}r)$.  Then we would have, for $r\geq1$,
\begin{equation}\label{307}
{\cal L}(n,r;\; m) \sim \frac{1}{m^{3/2}\sqrt{\log m}}\left\{{\cal V}(0)+\frac{1}{\log m}\left[\left(\frac{n}{r}+\sum_{l=1}^r\frac{1}{l}\right)\frac{ae^{-a}}{\sqrt{2\pi}}+{\cal V}^{(1)}(0)\right]\right\},
\end{equation}
and
\begin{equation}\label{308}
{\cal L}(n,0;\; m) \sim \frac{n+1}{m^{3/2}\sqrt{\log m}}\left[{\cal V}(0)+\frac{1}{\log m}{\cal V}^{(1)}(0)\right].
\end{equation}
We have verified, by numerical computations, that the order of magnitude of $\pi(k,r)$ on the $(n,r)$ scale does seem to be $O(m^{-3/2}(\log m)^{-1/2})$, that $\pi(k,0)$ is approximately proportional to $n+1$, and that $\pi(k,r)$ is approximately constant for $r\geq1$ and $k=m-O(1)$.  However, the problems with the matching suggest that there is yet another scale in the problem, which corresponds to $n\rightarrow \infty$ and $\tilde{\xi}\rightarrow 0$.  We have not been able to identify this new scale.

To summarize this subsection, we obtained results for $\pi(k,r)$ on the $(\bar{\xi}, r)$  scale, treating separately the cases $\bar{\xi}>\sqrt{r}$, $\; \bar{\xi}\sim \sqrt{r}$, and $\; \bar{\xi}<\sqrt{r}$,  for $r\geq1$.  For $r=0$ we gave results for $\bar{\xi}>1$, $\; \bar{\xi}\sim 1$, and $\bar{\xi}<1$.  For $\tilde{\xi}=\bar{\xi}\sqrt{\log m}=O(1)$ we obtained the simplified result in (\ref{8cd}) for $\pi(k,r)$ for $r\geq1$.  For $\pi(k,0)$, (\ref{301}) applies on the $\tilde{\xi}$ scale.  On the $(n, r)$ scale we obtained (\ref{307}) for $r\geq1$ and (\ref{308}) for $r=0$.  However, there is still a "gap" in the asymptotics between the $n$ and $\tilde{\xi}$ scale, a gap which we have not been able to fill.

\section{Numerical Studies}

We assess the accuracy of some of the asymptotic formulas we obtained.

In \textbf{Table 1} and \textbf{Table 2} we test our asymptotic results for $m\rightarrow \infty$ with a fixed $\rho<1$.  \textbf{Table 1} has $\pi(0,1)$, where (\ref{763}) applies with $k=0$ and $r=1$, $\; \pi(m/2, 1)$, where (\ref{738}) applies with $x=1/2$ and $r=1$, and $\pi(m, 0)$, where (\ref{749}) applies with $n=0$.  We consider $\rho=0.5$ and increase $m$ from 10 to 30.  The agreement is not particularly good for $(k, r)=(0, 1)$ but does improve significantly as $m$ increases.  For $(k,r) = (m/2,1)$ (boundary layer near $y=0$ in subsection \textbf{4.2}) and $(k,r)=(m, 0)$ (corner layer in subsection \textbf{4.3}) the agreement is good even for $m=10$, with errors of at most $10 \%$.  In \textbf{Table 4}, we consider $(k, r) = (0,m)$, ($m/2, m$), ($m, m$), so that $y=1$ and $x=$0, 1/2, 1.  The asymptotic formulas that apply are now (\ref{3}) (with (\ref{760})-(\ref{734})), (\ref{30}) (with (\ref{phi}) - (\ref{769})), and (\ref{727}) (with $x=y=1$).  We again consider $\rho = 0.5$ and $m=$ 10, 20, 30.  Now we obtain generally excellent agreement, with the worst error in \textbf{Table 2} being about $3 \%$, which occurs for $k=0$ and $r=m=10$.

These comparisons show that the asymptotics agree  reasonably well with the exact numerical values of $\pi(k, r)$.  The comparisons also clearly demonstrate the necessity of analyzing separately the different ranges of $(k, r)$, when $m\rightarrow \infty$.

\newpage

\begin{center}
\textbf{Table 1}
\begin{equation*}
\rho =0.5\; ; \;  \;  (k,r) = (0,1), \; (m/2, 1), \; (m,0)\; ; \;  \;  10\leq m\leq 30.
\end{equation*}
\end{center}

\begin{center}
 \begin{tabular} {|c|c||c|c|}       \hline
$(k,r)$          & $m$  & exact                & asymptotic \\
\hline $(0,1)$   & $10$ & $5.83\times10^{-6}$  & $9.76\times10^{-6}$\\
\hline $   $     & $20$ & $1.78\times10^{-9}$  & $2.38\times10^{-9}$\\
\hline $   $     & $30$ & $8.43\times10^{-13}$ & $1.03\times10^{-12}$\\
\hline $(m/2,1)$ & $10$ & $8.57\times10^{-5}$  & $7.81\times10^{-5}$\\
\hline $   $     & $20$ & $7.92\times10^{-8}$  & $7.62\times10^{-8}$\\
\hline $   $     & $30$ & $7.61\times10^{-11}$ & $7.45\times10^{-11}$\\
\hline $(m,0)$   & $10$ & $2.71\times10^{-4}$  & $2.44\times10^{-4}$\\
\hline $   $     & $20$ & $2.51\times10^{-7}$  & $2.38\times10^{-7}$\\
\hline $   $     & $30$ & $2.41\times10^{-10}$ & $2.32\times10^{-10}$\\
\hline
\end{tabular}\\
\end{center}

\begin{center}
\textbf{Table 2}
\begin{equation*}
\rho =0.5\; ; \;  \;  (k,r) = (0,m), \; (m/2, m), \; (m,m)\; ; \;  \;  10\leq m\leq 30.
\end{equation*}
\end{center}

\begin{center}
 \begin{tabular} {|c|c||c|c|}       \hline
$(k,r)$          & $m$  & exact                & asymptotic \\
\hline $(0,m)$   & $10$ & $7.59\times10^{-13}$ & $7.40\times10^{-13}$\\
\hline $   $     & $20$ & $3.64\times10^{-25}$ & $3.59\times10^{-25}$\\
\hline $   $     & $30$ & $1.52\times10^{-37}$ & $1.51\times10^{-37}$\\
\hline $(m/2,m)$ & $10$ & $4.86\times10^{-9}$  & $4.81\times10^{-9}$\\
\hline $   $     & $20$ & $4.93\times10^{-17}$ & $4.90\times10^{-17}$\\
\hline $   $     & $30$ & $5.02\times10^{-25}$ & $5.00\times10^{-25}$\\
\hline $(m,m)$   & $10$ & $3.67\times10^{-7}$  & $3.57\times10^{-7}$\\
\hline $   $     & $20$ & $3.45\times10^{-13}$ & $3.41\times10^{-13}$\\
\hline $   $     & $30$ & $3.28\times10^{-19}$ & $3.25\times10^{-19}$\\
\hline
\end{tabular}\\
\end{center}

\end{document}